\documentclass{article}
\usepackage{amsmath,amsfonts,amssymb,bm,hyperref}
\numberwithin{equation}{section}
\newtheorem{theorem}{Theorem}[section]
\newtheorem{lemma}{Lemma}[section]

\newtheorem{example}{Example}[section]
\newtheorem{remark}{Remark}[section]
\newtheorem{definition}{Definition}[section]
\usepackage[T1]{fontenc}
\usepackage[utf8]{inputenc}
\usepackage{authblk}

\author{Olga Chervova\thanks{O.Chervova@ucl.ac.uk, \url{http://www.homepages.ucl.ac.uk/\~ucahoch/}} }
\author{Robert Downes\thanks{R.Downes@ucl.ac.uk, \url{http://www.homepages.ucl.ac.uk/\~zcahc37/}} }
\author{Dmitri Vassiliev\thanks{D.Vassiliev@ucl.ac.uk, \url{http://www.homepages.ucl.ac.uk/\~ucahdva/}}}
\affil{Department of Mathematics,
University College London,\\ Gower Street, London WC1E 6BT, UK}

\begin{document}

\title{The spectral function of a first order system}
\maketitle
\begin{abstract}
We consider an elliptic self-adjoint first order pseudodifferential
operator acting on columns of $m$ complex-valued half-densities over
a connected compact $n$-dimensional manifold without boundary. The eigenvalues of the principal
symbol are assumed to be simple but no assumptions are made on their
sign, so the operator is not necessarily semi-bounded. We study the
spectral function, i.e.~the sum of squares of Euclidean norms of
eigen\-functions evaluated at a given point of the manifold, with
summation carried out over all eigenvalues between zero and a
positive $\lambda$. We derive a two-term asymptotic formula for the
spectral function as $\lambda$ tends to plus infinity.
We then restrict our study to the case when $m=2$, $n=3$, the operator is
differential and has trace-free principal symbol, and address the
question: is our operator a massless Dirac operator? We prove that
it is a massless Dirac operator if and only if the following two
conditions are satisfied at every point of the manifold: a) the
subprincipal symbol is proportional to the identity matrix and b)
the second asymptotic coefficient of the spectral function is zero.
\end{abstract}


\


\tableofcontents

\section{Main results}
\label{Main results}

The aim of the paper is to extend the classical results of \cite{DuiGui}
to systems. We are motivated by the following two observations.
\begin{itemize}
\item
To our knowledge, all previous publications on systems give
formulae for the second asymptotic coefficient that are either
incorrect or incomplete (i.e.~an algorithm for the calculation
of the second asymptotic coefficient rather than an actual formula).
The appropriate bibliographic review is presented in
Section~\ref{Bibliographic review}.
\item
Systems are fundamentally different from scalar operators in that
spectral analysis of systems reveals a very rich geometric
structure. An important example of an elliptic system
is the massless Dirac operator which is examined in detail in our paper.
\end{itemize}

Consider a first order classical pseudodifferential
operator $A$ acting on columns
$v=\begin{pmatrix}v_1&\ldots&v_m\end{pmatrix}^T$
of complex-valued half-densities
over a connected compact $n$-dimensional manifold $M$.
Throughout this paper we assume that $m\ge2$ and $n\ge2$.

We assume the coefficients of the operator $A$ to be infinitely smooth. We also
assume that the operator $A$ is formally self-adjoint (symmetric):
$\int_Mw^*Av\,dx=\int_M(Aw)^*v\,dx$ for all infinitely smooth
$v,w:M\to\mathbb{C}^m$. Here and further on
the superscript $\,{}^*\,$ in matrices, rows and columns
indicates Hermitian conjugation in $\mathbb{C}^m$
and $dx:=dx^1\ldots dx^n$, where $x=(x^1,\ldots,x^n)$ are local
coordinates on $M$.

Let $A_1(x,\xi)$ be the principal symbol of the operator $A$.
Here $\xi=(\xi_1,\ldots,\xi_n)$ is the variable dual to the position
variable $x$; in physics literature the $\xi$ would be referred to
as \emph{momentum}. Our principal symbol $A_1$ is an $m\times m$
Hermitian matrix-function on $T'M:=T^*M\setminus\{\xi=0\}$,
i.e.~on the cotangent bundle with the zero section removed.

Let $h^{(j)}(x,\xi)$ be the eigenvalues of the principal symbol. We
assume these eigenvalues to be nonzero (this is a version of the
ellipticity condition) but do not make any assumptions on their
sign. We also assume that the eigenvalues $h^{(j)}(x,\xi)$ are
simple for all $(x,\xi)\in T'M$. The techniques developed in our
paper do not work in the case when eigenvalues of the principal
symbol have variable multiplicity, though they could probably be
adapted to the case of constant multiplicity different from
multiplicity 1. The use of the letter ``$h$'' for an eigenvalue of
the principal symbol is motivated by the fact that later it will
take on the role of a Hamiltonian, see formula (\ref{Hamiltonian system of equations}).

We enumerate the eigenvalues of the principal symbol
$h^{(j)}(x,\xi)$ in increasing order, using a positive index
$j=1,\ldots,m^+$ for positive $h^{(j)}(x,\xi)$ and a negative index
$j=-1,\ldots,-m^-$ for negative $h^{(j)}(x,\xi)$. Here $m^+$ is the
number of positive eigenvalues of the principal symbol and $m^-$ is
the number of negative ones. Of course, $m^++m^-=m$.

Under the above assumptions $A$ is a self-adjoint operator, in the
full functional analytic sense, in the Hilbert space
$L^2(M;\mathbb{C}^m)$ (Hilbert space of square integrable
complex-valued column ``functions'') with domain $H^1(M;\mathbb{C}^m)$
(Sobolev space of complex-valued column ``functions'' which are
square integrable together with their first partial derivatives) and
the spectrum of $A$ is discrete. These facts are easily established
by constructing the parametrix (approximate inverse) of the operator
$A+iI$. Note that for the special case of the massless Dirac operator
a detailed examination of relevant functional analytic properties was
performed in Chapter~4 of \cite{friedrich}.

Let $\lambda_k$ and
$v_k=\begin{pmatrix}v_{k1}(x)&\ldots&v_{km}(x)\end{pmatrix}^T$ be
the eigenvalues and eigenfunctions of the operator $A$. The
eigenvalues $\lambda_k$ are enumerated in increasing order with
account of multiplicity,
using a positive index $k=1,2,\ldots$ for positive $\lambda_k$
and a nonpositive index $k=0,-1,-2,\ldots$ for nonpositive $\lambda_k$.
If the operator $A$ is bounded from below (i.e.~if $m^-=0$)
then the index $k$ runs from some integer value to $+\infty$;
if the operator $A$ is bounded from above (i.e.~if $m^+=0$)
then the index $k$ runs from $-\infty$ to some integer value;
and if the operator $A$ is unbounded from above
and from below (i.e.~if $m^+\ne0$ and $m^-\ne0$)
then the index $k$ runs from $-\infty$ to $+\infty$.

\

We will be studying the following three objects.

\

\textbf{Object 1.}
Our first object of study is the \emph{propagator},
which is the one-parameter family of operators defined as
\begin{equation}
\label{definition of wave group}
U(t):=e^{-itA}
=\sum_k e^{-it\lambda_k}v_k(x)\int_M[v_k(y)]^*(\,\cdot\,)\,dy\,,
\end{equation}
$t\in\mathbb{R}$.
The propagator provides a solution to the Cauchy problem
\begin{equation}
\label{initial condition most basic}
\left.w\right|_{t=0}=v
\end{equation}
for the dynamic equation
\begin{equation}
\label{dynamic equation most basic}
D_tw+Aw=0\,,
\end{equation}
where $D_t:=-i\partial/\partial t$.
Namely, it is easy to see that if the column of half-densities $v=v(x)$
is infinitely smooth,
then, setting
$\,w:=U(t)\,v$, we get a time-dependent column of half-densities $w(t,x)$
which is also infinitely smooth
and which satisfies the equation
(\ref{dynamic equation most basic})
and the initial condition
(\ref{initial condition most basic}).
The use of the letter ``$U$'' for the propagator is motivated by the
fact that for each $t$ the operator $U(t)$ is unitary.

\

\textbf{Object 2.}
Our second object of study is the \emph{spectral function},
which is the real density defined as
\begin{equation}
\label{definition of spectral function}
e(\lambda,x,x):=\sum_{0<\lambda_k<\lambda}\|v_k(x)\|^2,
\end{equation}
where $\|v_k(x)\|^2:=[v_k(x)]^*v_k(x)$ is the square of the
Euclidean norm of the eigenfunction $v_k$ evaluated at the point
$x\in M$ and $\lambda$ is a positive parameter (spectral parameter).

\

\textbf{Object 3.}
Our third and final object of study is the \emph{counting function}
\begin{equation}
\label{definition of counting function}
N(\lambda):=\,\sum_{0<\lambda_k<\lambda}1\ =\int_Me(\lambda,x,x)\,dx\,.
\end{equation}
In other words, $N(\lambda)$ is the number of eigenvalues $\lambda_k$
between zero and $\lambda$.

\

It is natural to ask the question: why, in defining the spectral function
(\ref{definition of spectral function})
and the counting function
(\ref{definition of counting function}),
did we choose to perform summation
over all \emph{positive} eigenvalues up to a given positive $\lambda$
rather than
over all \emph{negative} eigenvalues up to a given negative $\lambda$?
There is no particular reason. One case reduces to the other by the change
of operator $A\mapsto-A$. This issue will be revisited in
Section~\ref{Spectral asymmetry}.

Further on we assume that $m^+>0$, i.e.~that the operator
$A$ is unbounded from above.

\

Our objectives are as follows.

\

\textbf{Objective 1.}
We aim to construct the propagator
(\ref{definition of wave group})
explicitly in terms of
oscillatory integrals, modulo an integral operator with an
infinitely smooth, in the variables $t$, $x$ and $y$, integral kernel.

\

\textbf{Objectives 2 and 3.}
We aim to derive, under appropriate assumptions on Hamiltonian
trajectories, two-term asymptotics for the spectral function
(\ref{definition of spectral function})
and the counting function
(\ref{definition of counting function}),
i.e.~formulae of the type
\begin{equation}
\label{two-term asymptotic formula for spectral function}
e(\lambda,x,x)=a(x)\,\lambda^n+b(x)\,\lambda^{n-1}+o(\lambda^{n-1}),
\end{equation}
\begin{equation}
\label{two-term asymptotic formula for counting function}
N(\lambda)=a\lambda^n+b\lambda^{n-1}+o(\lambda^{n-1}),
\end{equation}
as $\lambda\to+\infty$.
Obviously, here we expect the real constants $a$, $b$ and real densities
$a(x)$, $b(x)$ to be related in accordance with
\begin{equation}
\label{a via a(x)}
a=\int_Ma(x)\,dx,
\end{equation}
\begin{equation}
\label{b via b(x)}
b=\int_Mb(x)\,dx.
\end{equation}

\

It is well known that the above three objectives are closely
related: if one achieves Objective 1, then Objectives 2 and 3 follow via
Fourier Tauberian theorems \cite{DuiGui,mybook,ivrii_book,Safarov_Tauberian_Theorems}.

\

We are now in a position to state our main results.

\

\textbf{Result 1.}
We construct the propagator as a sum of $m$ oscillatory integrals
\begin{equation}
\label{wave group as a sum of oscillatory integrals}
U(t)\overset{\operatorname{mod}C^\infty}=\sum_j
U^{(j)}(t)\,,
\end{equation}
where the phase function of each oscillatory integral
$U^{(j)}(t)$ is associated with the corresponding
Hamiltonian $h^{(j)}(x,\xi)$. The symbol of the oscillatory integral
$U^{(j)}(t)$ is a complex-valued $m\times m$ matrix-function
$u^{(j)}(t;y,\eta)$, where $y=(y^1,\ldots,y^n)$ is the position of the
source of the wave (i.e.~this is the same $y$ that appears in formula
(\ref{definition of wave group})) and
$\eta=(\eta_1,\ldots,\eta_n)$ is the corresponding dual variable
(covector at the point $y$).
When $|\eta|\to+\infty$, the symbol admits an asymptotic expansion
\begin{equation}
\label{decomposition of symbol of OI into homogeneous components}
u^{(j)}(t;y,\eta)=u^{(j)}_0(t;y,\eta)+u^{(j)}_{-1}(t;y,\eta)+\ldots
\end{equation}
into components positively homogeneous in $\eta$, with the subscript
indicating degree of homogeneity.

The formula for the principal symbol of the oscillatory integral
$U^{(j)}(t)$ is known
\cite{SafarovDSc,NicollPhD}
and reads as follows:
\begin{multline}
\label{formula for principal symbol of oscillatory integral}
u^{(j)}_0(t;y,\eta)=
[v^{(j)}(x^{(j)}(t;y,\eta),\xi^{(j)}(t;y,\eta))]
\,[v^{(j)}(y,\eta)]^*
\\
\times\exp
\left(
-i\int_0^tq^{(j)}(x^{(j)}(\tau;y,\eta),\xi^{(j)}(\tau;y,\eta))\,d\tau
\right),
\end{multline}
where $v^{(j)}(z,\zeta)$ is the normalised eigenvector of the principal
symbol $A_1(z,\zeta)$ corresponding to the eigenvalue (Hamiltonian)
$h^{(j)}(z,\zeta)$,
\ $(x^{(j)}(t;y,\eta),\xi^{(j)}(t;y,\eta))$ is the Hamiltonian trajectory
originating from the point $(y,\eta)$, i.e.~solution of the system of
ordinary differential equations (the dot denotes differentiation in $t$)
\begin{equation}
\label{Hamiltonian system of equations}
\dot x^{(j)}=h^{(j)}_\xi(x^{(j)},\xi^{(j)}),
\qquad
\dot\xi^{(j)}=-h^{(j)}_x(x^{(j)},\xi^{(j)})
\end{equation}
subject to the initial condition $\left.(x^{(j)},\xi^{(j)})\right|_{t=0}=(y,\eta)$,
\ $q^{(j)}:T'M\to\mathbb{R}$ is the function
\begin{equation}
\label{phase appearing in principal symbol}
q^{(j)}:=[v^{(j)}]^*A_\mathrm{sub}v^{(j)}
-\frac i2
\{
[v^{(j)}]^*,A_1-h^{(j)},v^{(j)}
\}
-i[v^{(j)}]^*\{v^{(j)},h^{(j)}\}
\end{equation}
and
\begin{equation}
\label{definition of subprincipal symbol}
A_\mathrm{sub}(z,\zeta):=
A_0(z,\zeta)+\frac i2
(A_1)_{z^\alpha\zeta_\alpha}(z,\zeta)
\end{equation}
is the subprincipal symbol of the operator $A$,
with the subscripts $z^\alpha$ and $\zeta_\alpha$
indicating partial derivatives and
the repeated index $\alpha$ indicating summation over $\alpha=1,\ldots,n$.
Curly brackets in formula
(\ref{phase appearing in principal symbol})
denote the Poisson bracket on matrix-functions
\begin{equation}
\label{Poisson bracket on matrix-functions}
\{P,R\}:=P_{z^\alpha}R_{\zeta_\alpha}-P_{\zeta_\alpha}R_{z^\alpha}
\end{equation}
and its further generalisation
\begin{equation}
\label{generalised Poisson bracket on matrix-functions}
\{P,Q,R\}:=P_{z^\alpha}QR_{\zeta_\alpha}-P_{\zeta_\alpha}QR_{z^\alpha}\,.
\end{equation}

As the derivation of formula
(\ref{formula for principal symbol of oscillatory integral})
was previously performed only in theses \cite{SafarovDSc,NicollPhD},
we repeat it in Sections
\ref{Algorithm for the construction of the wave group}
and
\ref{Leading transport equations}
of our paper.
Our derivation differs slightly from that in \cite{SafarovDSc} and \cite{NicollPhD}.

Formula (\ref{formula for principal symbol of oscillatory integral})
is invariant under changes of local coordinates on the manifold $M$,
i.e.~elements of the $m\times m$ matrix-function
$u^{(j)}_0(t;y,\eta)$ are scalars on
$\mathbb{R}\times T'M$.
Moreover, formula (\ref{formula for principal symbol of oscillatory integral})
is invariant under the transformation of the eigenvector of
the principal symbol
\begin{equation}
\label{gauge transformation of the eigenvector}
v^{(j)}\mapsto e^{i\phi^{(j)}}v^{(j)},
\end{equation}
where
\begin{equation}
\label{phase appearing in gauge transformation}
\phi^{(j)}:T'M\to\mathbb{R}
\end{equation}
is an arbitrary smooth function.
When some quantity is defined up to the action of a certain
transformation, theoretical physicists refer to such a
transformation as a \emph{gauge transformation}. We follow this
tradition.
Note that our particular gauge
transformation
(\ref{gauge transformation of the eigenvector}),
(\ref{phase appearing in gauge transformation})
is quite common in quantum mechanics:
when $\phi^{(j)}$ is a function of the position variable $x$ only
(i.e.~when $\phi^{(j)}:M\to\mathbb{R}$) this gauge transformation
is associated with electromagnetism.

Both Y.~Safarov \cite{SafarovDSc} and W.J.~Nicoll \cite{NicollPhD} assumed
that the operator $A$ is semi-bounded from below but this assumption
is not essential and their formula
(\ref{formula for principal symbol of oscillatory integral})
remains true in the more general case that we are dealing with.

However, knowing the principal symbol
(\ref{formula for principal symbol of oscillatory integral})
of the oscillatory integral $U^{(j)}(t)$
is not enough if one wants to derive two-term
asymptotics
(\ref{two-term asymptotic formula for spectral function})
and
(\ref{two-term asymptotic formula for counting function}).
One needs information about $u^{(j)}_{-1}(t;y,\eta)$,
the component of the symbol of the
oscillatory integral $U^{(j)}(t)$
which is positively homogeneous in $\eta$ of degree~-1,
see formula (\ref{decomposition of symbol of OI into homogeneous components}),
but here the problem is that $u^{(j)}_{-1}(t;y,\eta)$
is not a true invariant in the sense that it depends on the choice of
phase function in the oscillatory integral. We overcome this difficulty
by observing that $U^{(j)}(0)$ is a pseudodifferential operator, hence,
it has a well-defined subprincipal symbol
$[U^{(j)}(0)]_\mathrm{sub}$. We prove that
\begin{equation}
\label{subprincipal symbol of OI at time zero}
\operatorname{tr}[U^{(j)}(0)]_\mathrm{sub}
=-i\{[v^{(j)}]^*,v^{(j)}\}
\end{equation}
and subsequently show that information contained in formulae
(\ref{formula for principal symbol of oscillatory integral})
and
(\ref{subprincipal symbol of OI at time zero})
is sufficient for the derivation of two-term
asymptotics
(\ref{two-term asymptotic formula for spectral function})
and
(\ref{two-term asymptotic formula for counting function}).

Note that the RHS of formula (\ref{subprincipal symbol of OI at time zero})
is invariant under the gauge transformation
(\ref{gauge transformation of the eigenvector}),
(\ref{phase appearing in gauge transformation}).

Formula~(\ref{subprincipal symbol of OI at time zero})
plays a central role in our paper.
Sections~\ref{Algorithm for the construction of the wave group}
and~\ref{Leading transport equations}
provide auxiliary material needed for the proof of
formula~(\ref{subprincipal symbol of OI at time zero}),
whereas the actual proof
of formula~(\ref{subprincipal symbol of OI at time zero})
is given in Section~\ref{Proof of formula}.

Let us elaborate briefly on the geometric meaning of the
RHS of (\ref{subprincipal symbol of OI at time zero})
(a more detailed exposition is presented in Section~\ref{U(1) connection}).
The eigenvector of the principal
symbol is defined up to a gauge transformation
(\ref{gauge transformation of the eigenvector}),
(\ref{phase appearing in gauge transformation})
so it is natural to introduce a $\mathrm{U}(1)$ connection on $T'M$ as
follows: when parallel transporting an eigenvector of the principal
symbol along a curve in $T'M$ we require that the derivative of the
eigenvector along the curve be orthogonal to the eigenvector itself.
This is equivalent to the introduction of an (intrinsic) electromagnetic
field on $T'M$, with the $2n$-component real quantity
\begin{equation}
\label{electromagnetic covector potential}
i\,(\,[v^{(j)}]^*v^{(j)}_{x^\alpha}\,,\,[v^{(j)}]^*v^{(j)}_{\xi_\gamma}\,)
\end{equation}
playing the role of the electromagnetic covector potential. Our
quantity (\ref{electromagnetic covector potential}) is a 1-form on
$T'M$, rather than on $M$ itself as is the case in ``traditional''
electromagnetism. The above $\mathrm{U}(1)$ connection generates
curvature which is a 2-form on $T'M$, an analogue of the
electromagnetic tensor. Out of this curvature 2-form one can
construct, by contraction of indices, a real scalar. This scalar
curvature is the expression appearing in the RHS of
formula (\ref{subprincipal symbol of OI at time zero}).

Observe now that $\sum_jU^{(j)}(0)$ is the identity operator
on half-densities. The subprincipal symbol of the identity operator
is zero, so formula (\ref{subprincipal symbol of OI at time zero})
implies
\begin{equation}
\label{sum of curvatures is zero}
\sum_j\{[v^{(j)}]^*,v^{(j)}\}=0.
\end{equation}
One can check the identity (\ref{sum of curvatures is zero})
directly, without constructing the oscillatory integrals
$U^{(j)}(t)$: it follows from the fact that the $v^{(j)}(x,\xi)$ form an
orthonormal basis, see end of Section \ref{U(1) connection} for details.
We mentioned the identity
(\ref{sum of curvatures is zero})
in order to highlight, once again, the fact that the curvature effects
we have identified are specific to systems and do not have an
analogue in the scalar case.

\

\textbf{Results 2 and 3.}
We prove, under appropriate assumptions on Hamiltonian trajectories
(see Theorems~\ref{theorem spectral function unmollified two term}
and \ref{theorem counting function unmollified two term} for details),
asymptotic formulae
(\ref{two-term asymptotic formula for spectral function})
and
(\ref{two-term asymptotic formula for counting function})
with
\begin{equation}
\label{formula for a(x)}
a(x)=\sum_{j=1}^{m^+}
\ \int\limits_{h^{(j)}(x,\xi)<1}{d{\hskip-1pt\bar{}}\hskip1pt}\xi\,,
\end{equation}
\begin{multline}
\label{formula for b(x)}
b(x)=-n\sum_{j=1}^{m^+}
\ \int\limits_{h^{(j)}(x,\xi)<1}
\Bigl(
[v^{(j)}]^*A_\mathrm{sub}v^{(j)}
\\
-\frac i2
\{
[v^{(j)}]^*,A_1-h^{(j)},v^{(j)}
\}
+\frac i{n-1}h^{(j)}\{[v^{(j)}]^*,v^{(j)}\}
\Bigr)(x,\xi)\,
{d{\hskip-1pt\bar{}}\hskip1pt}\xi\,,
\end{multline}
and $a$ and $b$ expressed via the above densities
(\ref{formula for a(x)})
and
(\ref{formula for b(x)})
as
(\ref{a via a(x)})
and
(\ref{b via b(x)}).
In
(\ref{formula for a(x)})
and
(\ref{formula for b(x)})
\,${d{\hskip-1pt\bar{}}\hskip1pt}\xi$ is shorthand for
${d{\hskip-1pt\bar{}}\hskip1pt}\xi:=(2\pi)^{-n}\,d\xi
=(2\pi)^{-n}\,d\xi_1\ldots d\xi_n$,
and the Poisson bracket on matrix-functions
$\{\,\cdot\,,\,\cdot\,\}$
and its further generalisation
$\{\,\cdot\,,\,\cdot\,,\,\cdot\,\}$
are defined by formulae
(\ref{Poisson bracket on matrix-functions})
and
(\ref{generalised Poisson bracket on matrix-functions})
respectively.

To our knowledge, formula (\ref{formula for b(x)}) is a new result.
Note that in \cite{SafarovDSc} this formula
(more precisely, its integrated over $M$ version (\ref{b via b(x)}))
was written incorrectly,
without the curvature terms
$\,-\frac{ni}{n-1}\int h^{(j)}\{[v^{(j)}]^*,v^{(j)}\}$.
See also Section~\ref{Bibliographic review}
where we give a more detailed bibliographic review.

It is easy to see that the right-hand sides of
(\ref{formula for a(x)})
and
(\ref{formula for b(x)})
behave as densities under changes of local coordinates
on the manifold $M$ and that these expressions are invariant
under gauge transformations
(\ref{gauge transformation of the eigenvector}),
(\ref{phase appearing in gauge transformation})
of the eigenvectors of the principal symbol.
Moreover, the right-hand sides of
(\ref{formula for a(x)})
and
(\ref{formula for b(x)})
are unitarily invariant,
i.e.~invariant under transformations of the operator
\begin{equation}
\label{unitary transformation of operator A}
A\mapsto RAR^*,
\end{equation}
where
\begin{equation}
\label{matrix appearing in unitary transformation of operator}
R:M\to\mathrm{U}(m)
\end{equation}
is an arbitrary smooth unitary matrix-function.
The fact that the RHS of
(\ref{formula for b(x)})
is unitarily invariant is non-trivial: the appropriate calculations
are presented in Section~\ref{U(m) invariance}.
The observation that without the curvature terms
$\,-\frac{ni}{n-1}\int h^{(j)}\{[v^{(j)}]^*,v^{(j)}\}$
(as in \cite{SafarovDSc}) the RHS of
(\ref{formula for b(x)}) is not unitarily invariant
was a major motivating factor in the writing of this paper.

\

We will now start making additional assumptions which will,
in the end, allow us to provide a simple
spectral theoretic characterisation of the massless Dirac operator.

\

\textbf{Additional assumption 1:}
\begin{equation}
\label{Assumption 1}
m=2\quad\text{and}\quad\operatorname{tr}A_1=0.
\end{equation}

In this case we can simplify notation by denoting the positive
eigenvalue of the principal symbol by
$h^+$, the corresponding eigenvector by
$v^+=\begin{pmatrix}v^+_1\\ v^+_2\end{pmatrix}$
and Hamiltonian
trajectories by $(x^+(t;y,\eta),\xi^+(t;y,\eta))$.
Obviously, the other eigenvalue of the principal
symbol is $-h^+$, the corresponding eigenvector is
$\begin{pmatrix}-\bar v^+_2\\\bar v^+_1\end{pmatrix}$
and Hamiltonian trajectories
are $(x^+(-t;y,\eta),\xi^+(-t;y,\eta))$ (time reversal).
Note that in theoretical physics the antilinear transformation
\begin{equation}
\label{charge conjugation at the level of principal symbol}
\begin{pmatrix}v^+_1\\ v^+_2\end{pmatrix}
\overset{\mathrm{C}}\mapsto
\begin{pmatrix}-\bar v^+_2\\\bar v^+_1\end{pmatrix}
\end{equation}
is referred to as \emph{charge conjugation} \cite{PCT}.

Moreover, in this case the two scalar invariants,
$\{
[v^+]^*,A_1-h^+,v^+
\}$
and
$h^+\{[v^+]^*,v^+\}$,
appearing in formula (\ref{formula for b(x)})
cease being independent and become related as
$\{
[v^+]^*,A_1-h^+,v^+
\}
=-2h^+\{[v^+]^*,v^+\}$.
Hence, formulae
(\ref{formula for a(x)})
and
(\ref{formula for b(x)})
simplify and now read
\begin{equation}
\label{formula for a(x) with assumption 1}
a(x)=
\int\limits_{h^+(x,\xi)<1}{d{\hskip-1pt\bar{}}\hskip1pt}\xi\,,
\end{equation}
\begin{equation}
\label{formula for b(x) with assumption 1}
b(x)=-n\int\limits_{h^+(x,\xi)<1}
\Bigl(
[v^+]^*A_\mathrm{sub}v^+
+\frac{n}{n-1}ih^+\{[v^+]^*,v^+\}
\Bigr)(x,\xi)\,
{d{\hskip-1pt\bar{}}\hskip1pt}\xi\,.
\end{equation}

\

\textbf{Additional assumption 2:}
\begin{equation}
\label{Assumption 2}
\text{the operator $A$ is differential}.
\end{equation}

In this case there are three further simplifications.

Firstly, the dimension of the manifold can only be $n=2$ or $n=3$.
This follows from the ellipticity condition and the fact that the
dimension of the real vector space of trace-free Hermitian $2\times2$
matrices is 3.

Secondly, the subprincipal symbol $A_\mathrm{sub}$ does not depend
on the dual variable~$\xi$ (momentum) and is a function of $x$
(position) only.

Thirdly, we acquire a geometric object, the \emph{metric}. Indeed, the
determinant of the principal symbol is a negative definite quadratic
form
\begin{equation}
\label{definition of metric}
\det A_1(x,\xi)=-g^{\alpha\beta}\xi_\alpha\xi_\beta
\end{equation}
and the coefficients $g^{\alpha\beta}(x)$,
$\alpha,\beta=1,\ldots,n$,
appearing in (\ref{definition of metric}) can be interpreted as the
components of a (contravariant) Riemannian metric. This implies, in
particular, that our Hamiltonian (positive eigenvalue of the
principal symbol) takes the form
\begin{equation}
\label{Hamiltonian expressed via metric}
h^+(x,\xi)=\sqrt{g^{\alpha\beta}(x)\,\xi_\alpha\xi_\beta}
\end{equation}
and the $x$-components of our Hamiltonian trajectories become
geodesics. Moreover, formulae
(\ref{formula for a(x) with assumption 1})
and
(\ref{a via a(x)})
simplify and now read
\begin{equation}
\label{formula for a(x) with assumptions 1 and 2}
a(x)=
(2\pi)^{-n}\,\omega_n\,\sqrt{\det g_{\alpha\beta}(x)}\,,
\end{equation}
\begin{equation}
\label{formula for a with assumptions 1 and 2}
a=
(2\pi)^{-n}\,\omega_n\operatorname{Vol}M\,,
\end{equation}
where $\omega_n$ is the volume of the unit ball in $\mathbb{R}^n$
and $\operatorname{Vol}M$ is the $n$-dimensional volume of the
Riemannian manifold $M$.

\

\textbf{Additional assumption 3:}
\begin{equation}
\label{Assumption 3}
n=3.
\end{equation}

In this case there are three more simplifications.

Firstly, the manifold $M$ is bound to be parallelizable (and, hence,
orientable). The relevant argument is presented in the beginning
of Section~\ref{Teleparallel connection}. From this point we work
only in local coordinates with prescribed orientation.

Secondly, we acquire the identity
\begin{equation}
\label{identity leading to definition of relative orientation}
\det g^{\alpha\beta}=
-\frac14\bigl[\operatorname{tr}\bigl((A_1)_{\xi_1}(A_1)_{\xi_2}(A_1)_{\xi_3}\bigr)\bigr]^2
\end{equation}
which allows us to define the topological invariant
\begin{equation}
\label{definition of relative orientation}
\mathbf{c}:=-\frac i2\sqrt{\det g_{\alpha\beta}}\,\operatorname{tr}
\bigl((A_1)_{\xi_1}(A_1)_{\xi_2}(A_1)_{\xi_3}\bigr).
\end{equation}
The number $\mathbf{c}$ defined by formula
(\ref{definition of relative orientation})
can take only two values, $+1$ or $-1$,
and describes the orientation of the principal symbol $A_1(x,\xi)$
relative to the chosen orientation of local coordinates,
see formula (\ref{definition of relative orientation more natural})
for a more natural geometric definition.
In calling the number $\mathbf{c}$ a topological
invariant we are referring to the topology of deformations of the elliptic trace-free principal
symbol $A_1(x,\xi)$ rather than the deformations of the manifold $M$ itself.

Thirdly, we acquire a new differential geometric object, namely, a
\emph{teleparallel connection}. This is an affine connection
defined as follows. Suppose we have a covector $\eta$ based at
the point $y\in M$ and we want to construct a parallel covector
$\xi$ based at the point $x\in M$. This is done by solving the
linear system of equations
\begin{equation}
\label{definition of parallel transport}
A_1(x,\xi)=A_1(y,\eta).
\end{equation}
Equation (\ref{definition of parallel transport}) is equivalent to a
system of three real linear algebraic equations for the three
real unknowns, components of the covector $\xi$, and it is easy to
see that this system has a unique solution. It is also easy to see that
the affine connection defined by formula (\ref{definition of parallel transport})
preserves the Riemannian norm of covectors,
i.e.~$g^{\alpha\beta}(x)\,\xi_\alpha\xi_\beta
=g^{\alpha\beta}(y)\,\eta_\alpha\eta_\beta$, hence, it is metric
compatible. The parallel transport defined by formula
(\ref{definition of parallel transport}) does not depend on the
curve along which we transport the (co)vector, so our connection has
zero curvature. The word ``teleparallel'' (parallel at a distance)
is used in theoretical physics to describe metric compatible affine
connections with zero curvature.
This terminology goes back to the works of A.~Einstein and
\'E.~Cartan \cite{unzicker-2005-,MR2276051,MR543192}, though Cartan preferred
to use the term ``absolute parallelism'' rather than ``teleparallelism''.

The teleparallel connection coefficients
$\Gamma^\alpha{}_{\beta\gamma}(x)$ can be written down explicitly in
terms of the principal symbol,
see formula (\ref{formula for teleparallel connection coefficients}),
and this allows us to define yet another geometric object --- the torsion tensor
\begin{equation}
\label{definition of torsion}
T^\alpha{}_{\beta\gamma}
:=\Gamma^\alpha{}_{\beta\gamma}-\Gamma^\alpha{}_{\gamma\beta}\,.
\end{equation}
Further on we raise and lower indices of the torsion tensor using
the metric. Torsion is a rank three tensor antisymmetric in the last
two indices. Because we are working in dimension three, it is
convenient, as in \cite{rotational_elasticity},
to apply the Hodge star in the last two indices
and deal with the rank two tensor
\begin{equation}
\label{definition of torsion with a star}
\overset{*}T{}^\alpha{}_\beta:=
\frac12\,T^{\alpha\gamma\delta}\,\varepsilon_{\gamma\delta\beta}
\,\sqrt{\det g_{\mu\nu}}
\end{equation}
rather than with the rank three tensor $T$. Here $\varepsilon$ is
the totally antisymmetric quantity, $\varepsilon_{123}:=+1$.

The teleparallel connection is a simpler geometric object than
the $\mathrm{U}(1)$ connection because the coefficients of
the teleparallel connection do not depend on the dual variable
(momentum), i.e.~they are ``functions'' on the base manifold $M$.
The relationship between the two connections is established
in Section~\ref{Teleparallel connection} where
we show that the scalar curvature of the $\mathrm{U}(1)$ connection is
expressed via the torsion of the teleparallel connection and the metric as
\begin{equation}
\label{curvature via torsion and metric}
-i\{[v^+]^*,v^+\}(x,\xi)
=\frac{\mathbf{c}}2\,
\frac
{\overset{*}T{}^{\alpha\beta}(x)\,\xi_\alpha\xi_\beta}
{(g^{\mu\nu}(x)\,\xi_\mu\xi_\nu)^{3/2}}\,.
\end{equation}

Integration of both terms appearing in formula
(\ref{formula for b(x) with assumption 1})
can now be carried out explicitly, giving
\begin{equation}
\label{integration in xi of term with subprincipal symbol}
\int\limits_{h^+(x,\xi)<1}
(
[v^+]^*A_\mathrm{sub}v^+
)(x,\xi)\,
{d{\hskip-1pt\bar{}}\hskip1pt}\xi
=\frac1{12\pi^2}
\bigl(\operatorname{tr}A_\mathrm{sub}\,\sqrt{\det g_{\alpha\beta}}\,\bigr)(x)\,,
\end{equation}
\begin{equation}
\label{integration in xi of term with curvature}
-i\int\limits_{h^+(x,\xi)<1}h^+\{[v^+]^*,v^+\}(x,\xi)\,
{d{\hskip-1pt\bar{}}\hskip1pt}\xi
=\frac{\mathbf{c}}{36\pi^2}
\bigl(\operatorname{tr}\overset{*}T\,\sqrt{\det g_{\alpha\beta}}\,\bigr)(x)\,,
\end{equation}
where $\operatorname{tr}\overset{*}T:=\overset{*}T{}^\alpha{}_\alpha$.
Note that $\operatorname{tr}\overset{*}T$ corresponds to one of the
three irreducible pieces of torsion, namely, the piece which is labelled
by theoretical physicists by the adjective ``axial'', see
\cite{rotational_elasticity,cartantorsionreview}
for details; it is interesting that this is exactly the
irreducible piece of torsion which is used when one models the
neutrino~\cite{MR2670535} or the electron \cite{mathematika} by means of Cosserat elasticity.
Formula (\ref{integration in xi of term with curvature})
follows immediately from (\ref{curvature via torsion and metric}),
whereas formula (\ref{integration in xi of term with subprincipal symbol})
is somewhat less obvious. In order to see where
formula (\ref{integration in xi of term with subprincipal symbol})
comes from one has to write the orthogonal projection
$v^+(x,\xi)\,[v^+(x,\xi)]^*$ as
$
v^+(x,\xi)\,[v^+(x,\xi)]^*=\frac1{2h^+(x,\xi)}(A_1(x,\xi)+h^+(x,\xi)\,I)
$
and use the fact that the principal symbol $A_1(x,\xi)$ is an odd function of $\xi$.

Substituting
(\ref{Assumption 3}),
(\ref{integration in xi of term with subprincipal symbol})
and
(\ref{integration in xi of term with curvature})
into
(\ref{formula for b(x) with assumption 1})
we get
\begin{equation}
\label{formula for b(x) with assumption 1, 2 and 3}
b(x)=\frac1{8\pi^2}
\bigl(
\bigl(
\mathbf{c}\operatorname{tr}\overset{*}T
-
2\operatorname{tr}A_\mathrm{sub}
\bigr)
\sqrt{\det g_{\alpha\beta}}
\,\bigr)(x)\,.
\end{equation}
An explicit self-contained expression for
$\operatorname{tr}\overset{*}T$
is given in formula
(\ref{explicit formula for the trace of torsion with a star}).

Note that the two traces appearing in formula
(\ref{formula for b(x) with assumption 1, 2 and 3}) have a different
meaning: $\operatorname{tr}\overset{*}T$ is the trace of a $3\times3$ tensor,
whereas $\operatorname{tr}A_\mathrm{sub}$ is the trace of a $2\times2$ matrix.

\

We now turn our attention to the massless Dirac operator. This
operator is defined in Appendix \ref{The massless Dirac operator},
see formula (\ref{definition of Weyl operator}),
and it does not fit into our scheme because this is an operator
acting on a 2-component complex-valued spinor (Weyl spinor) rather than
a pair of complex-valued half-densities. However, on a parallelizable manifold
components of a spinor can be identified with half-densities. We
call the resulting operator
\emph{the massless Dirac operator on half-densities}. The explicit
formula for the massless Dirac operator on half-densities is
(\ref{definition of Weyl operator on half-densities}).

The massless Dirac operator on half-densities is an operator of the
type described in this section (elliptic self-adjoint first order
operator acting on a column of complex-valued half-densities) which,
moreover, satisfies the additional assumptions (\ref{Assumption 1}),
(\ref{Assumption 2}) and (\ref{Assumption 3}). We address the
question: is a given operator $A$ a massless Dirac operator? The
answer is given by the following theorem which
we prove in Section~\ref{Proof of Theorem}.

\begin{theorem}
\label{main theorem} Let $A$ be an elliptic self-adjoint first order
pseudodifferential operator acting on columns of $m$ complex-valued
half-densities over a compact $n$-dimensional manifold. Suppose also
that this operator satisfies the additional assumptions
(\ref{Assumption 1}), (\ref{Assumption 2}) and (\ref{Assumption 3}).
Then $A$ is a massless Dirac operator on half-densities if and
only if the following two conditions are satisfied at every point of the
manifold $M$: a) the subprincipal symbol of the operator,
$A_\mathrm{sub}(x)$, is proportional to the identity matrix and b)
the second asymptotic coefficient of the spectral function, $b(x)$,
is zero.
\end{theorem}

The theorem stated above warrants the following remarks.
\begin{itemize}
\item
In stating Theorem \ref{main theorem} we did not make any
assumptions on Hamiltonian trajectories (loops). The
second asymptotic coefficient (\ref{formula for b(x) with assumption 1, 2 and 3})
is, in itself,
well-defined irrespective of how many loops we have. If one wishes
to reformulate the asymptotic formula
(\ref{two-term asymptotic formula for spectral function})
in such a way that it remains valid without assumptions on the
number of loops, this can easily be achieved, say, by
taking a convolution with a function
from Schwartz space $\mathcal{S}(\mathbb{R})$.
See Theorem~\ref{theorem spectral function mollified} for details.
\item
Conditions a) and b) in Theorem \ref{main theorem}
are invariant under special unitary transformations,
i.e.~transformations of the operator
(\ref{unitary transformation of operator A})
where $R=R(x)$ is an arbitrary smooth special unitary matrix-function.
This is not surprising as the massless Dirac operator is designed
around the concept of $\mathrm{SU}(2)$ invariance,
see Property 4 in Appendix \ref{The massless Dirac operator}.
\item
Condition b) in Theorem \ref{main theorem}
is actually invariant under the action of a broader group:
the unitary matrix-function appearing in formula
(\ref{unitary transformation of operator A})
does not have to be special.
\end{itemize}

\section{Algorithm for the construction of the propagator}
\label{Algorithm for the construction of the wave group}

We construct the propagator as a sum of $m$ oscillatory integrals
(\ref{wave group as a sum of oscillatory integrals}) where each
integral is of the form
\begin{equation}
\label{algorithm equation 1}
U^{(j)}(t)
=
\int e^{i\varphi^{(j)}(t,x;y,\eta)}
\,u^{(j)}(t;y,\eta)
\,\varsigma^{(j)}(t,x;y,\eta)\,d_{\varphi^{(j)}}(t,x;y,\eta)\,
(\ \cdot\ )\,dy\,{d{\hskip-1pt\bar{}}\hskip1pt}\eta\,.
\end{equation}
Here we use notation from the book \cite{mybook}, only adapted to
systems. Namely, the expressions appearing in formula
(\ref{algorithm equation 1}) have the following meaning.

\begin{itemize}
\item
The function
$\varphi^{(j)}$ is a phase function, i.e.~a function
$\mathbb{R}\times M\times T'M\to\mathbb{C}$
positively homogeneous in $\eta$ of degree 1 and satisfying
the conditions
\begin{equation}
\label{algorithm equation 2}
\varphi^{(j)}(t,x;y,\eta)
=(x-x^{(j)}(t;y,\eta))^\alpha\,\xi^{(j)}_\alpha(t;y,\eta)
+O(|x-x^{(j)}(t;y,\eta)|^2),
\end{equation}
\begin{equation}
\label{algorithm equation 3}
\operatorname{Im}\varphi^{(j)}(t,x;y,\eta)\ge0,
\end{equation}
\begin{equation}
\label{algorithm equation 4}
\det\varphi^{(j)}_{x^\alpha\eta_\beta}(t,x^{(j)}(t;y,\eta);y,\eta)\ne0.
\end{equation}
Recall that according to Corollary 2.4.5 from
\cite{mybook} we are guaranteed to have
(\ref{algorithm equation 4}) if we choose a phase function
\begin{multline}
\label{algorithm equation 5}
\varphi^{(j)}(t,x;y,\eta)
=(x-x^{(j)}(t;y,\eta))^\alpha\,\xi^{(j)}_\alpha(t;y,\eta)
\\
+\frac12C^{(j)}_{\alpha\beta}(t;y,\eta)
\,(x-x^{(j)}(t;y,\eta))^\alpha\,(x-x^{(j)}(t;y,\eta))^\beta
\\
+O(|x-x^{(j)}(t;y,\eta)|^3)
\end{multline}
with complex-valued symmetric matrix-function $C^{(j)}_{\alpha\beta}$ satisfying the strict inequality
$\operatorname{Im}C^{(j)}>0$
(our original requirement (\ref{algorithm equation 3}) implies only the
non-strict inequality $\operatorname{Im}C^{(j)}\ge0$).
Note that even though the matrix-function $C^{(j)}_{\alpha\beta}$ is
not a tensor, the inequalities $\operatorname{Im}C^{(j)}\ge0$ and
$\operatorname{Im}C^{(j)}>0$ are invariant under transformations of
local coordinates $x$; see Remark 2.4.9 in \cite{mybook} for details.

\item
The quantity $u^{(j)}$ is the symbol of our oscillatory integral,
i.e.~a complex-valued $m\times m$ matrix-function
$\mathbb{R}\times T'M\to\mathbb{C}^{m^2}$
which admits the asymptotic expansion
(\ref{decomposition of symbol of OI into homogeneous components}).
The symbol is the unknown quantity in our construction.

\item
The quantity $d_{\varphi^{(j)}}$ is defined in accordance with
formula (2.2.4) from \cite{mybook} as
\begin{equation}
\label{algorithm equation 6}
d_{\varphi^{(j)}}(t,x;y,\eta)
:=({\det}^2\varphi^{(j)}_{x^\alpha\eta_\beta})^{1/4}
=|\det\varphi^{(j)}_{x^\alpha\eta_\beta}|^{1/2}
\,e^{\,i\arg({\det}^2\varphi^{(j)}_{x^\alpha\eta_\beta})/4}.
\end{equation}
Note that in view of (\ref{algorithm equation 4}) our $d_{\varphi^{(j)}}$
is well-defined and smooth for $x$ close to $x^{(j)}(t;y,\eta)$. It is known
\cite{mybook} that under coordinate transformations
$d_{\varphi^{(j)}}$ behaves as a half-density in $x$ and
as a half-density to the power $-1$ in $y$.

In formula (\ref{algorithm equation 6}) we wrote
$({\det}^2\varphi^{(j)}_{x^\alpha\eta_\beta})^{1/4}$
rather than
$(\det\varphi^{(j)}_{x^\alpha\eta_\beta})^{1/2}$
in order to make this expression truly invariant under coordinate transformations.
Recall that local coordinates $x$ and $y$ are chosen independently
and that $\eta$ is a covector based at the point $y$.
Consequently,
$\det\varphi^{(j)}_{x^\alpha\eta_\beta}$ changes sign under inversions
of local coordinates $x$ or $y$,
whereas ${\det}^2\varphi^{(j)}_{x^\alpha\eta_\beta}$ retains sign under
inversions.

The choice of (smooth) branch of $\arg({\det}^2\varphi^{(j)}_{x^\alpha\eta_\beta})$ is assumed
to be fixed. Thus, for a given phase function
$\varphi^{(j)}$
formula (\ref{algorithm equation 6}) defines the quantity
$d_{\varphi^{(j)}}$
uniquely up to a factor $e^{ik\pi/2}$, $k=0,1,2,3$.
Observe now that if we set $t=0$ and choose the same local coordinates
for $x$ and $y$, we get $\varphi^{(j)}_{x^\alpha\eta_\beta}(0,y;y,\eta)=I$.
This implies that we can fully specify the choice of branch of
$\arg({\det}^2\varphi^{(j)}_{x^\alpha\eta_\beta})$
by requiring that
$d_{\varphi^{(j)}}(0,y;y,\eta)=1$.

The purpose of the introduction of the factor $d_{\varphi^{(j)}}$
in (\ref{algorithm equation 1}) is twofold.
\begin{itemize}
\item[(a)]
It ensures that the symbol $u^{(j)}$ is a function on
$\mathbb{R}\times T'M$ in the full differential geometric sense of the word,
i.e.~that it is invariant under transformations of local coordinates $x$ and $y$.
\item[(b)]
It ensures that the principal symbol $u^{(j)}_0$ does not depend
on the choice of phase function $\varphi^{(j)}$.
See Remark 2.2.8 in \cite{mybook} for more details.
\end{itemize}

\item
The quantity $\varsigma^{(j)}$ is a smooth cut-off function
$\mathbb{R}\times M\times T'M\to\mathbb{R}$
satisfying the following conditions.
\begin{itemize}
\item[(a)]
$\varsigma^{(j)}(t,x;y,\eta)=0$ on the set
$\{(t,x;y,\eta):\ |h^{(j)}(y,\eta)|\le1/2\}$.
\item[(b)]
$\varsigma^{(j)}(t,x;y,\eta)=1$ on the intersection
of a small conic neighbourhood of the set
\begin{equation}
\label{algorithm equation 7.1}
\{(t,x;y,\eta):\ x=x^{(j)}(t;y,\eta)\}
\end{equation}
with the set $\{(t,x;y,\eta):\ |h^{(j)}(y,\eta)|\ge1\}$.
\item[(c)]
$\varsigma^{(j)}(t,x;y,\lambda\eta)=\varsigma^{(j)}(t,x;y,\eta)$
for $\,|h^{(j)}(y,\eta)|\ge1$, $\,\lambda\ge1$.
\end{itemize}

\item
It is known (see Section 2.3 in \cite{mybook} for details)
that Hamiltonian trajectories generated by a Hamiltonian
$h^{(j)}(x,\xi)$ positively homogeneous in $\xi$ of degree~1
satisfy the identity
\begin{equation}
\label{algorithm equation 7.1.5}
(x^{(j)}_\eta)^{\alpha\beta}\xi^{(j)}_\alpha=0,
\end{equation}
where $(x^{(j)}_\eta)^{\alpha\beta}:=\partial(x^{(j)})^\alpha/\partial\eta_\beta$.
Formulae (\ref{algorithm equation 2}) and (\ref{algorithm equation 7.1.5})
imply
\begin{equation}
\label{algorithm equation 7.2}
\varphi^{(j)}_\eta(t,x^{(j)}(t;y,\eta);y,\eta)=0.
\end{equation}
This allows us to apply the stationary phase method in the neighbourhood
of the set (\ref{algorithm equation 7.1}) and disregard what happens
away from it.
\end{itemize}

\

Our task now is to construct the symbols  $u^{(j)}_0(t;y,\eta)$, $j=1,\ldots,m$,
so that our oscillatory integrals $U^{(j)}(t)$, $j=1,\ldots,m$,
satisfy the dynamic equations
\begin{equation}
\label{algorithm equation 8}
(D_t+A(x,D_x))\,U^{(j)}(t)\overset{\operatorname{mod}C^\infty}=0
\end{equation}
and initial condition
\begin{equation}
\label{algorithm equation 9}
\sum_jU^{(j)}(0)\overset{\operatorname{mod}C^\infty}=I\,,
\end{equation}
where $I$ is the identity operator on half-densities;
compare with formulae
(\ref{dynamic equation most basic}),
(\ref{initial condition most basic})
and (\ref{wave group as a sum of oscillatory integrals}).
Note that the pseudodifferential operator $A$ in formula
(\ref{algorithm equation 8}) acts on the oscillatory integral
$U(t)$ in the variable $x$; say, if $A$ is a differential
operator this means that in order to evaluate $A\,U^{(j)}(t)$
one has to perform the appropriate
differentiations of the oscillatory integral
(\ref{algorithm equation 1})
in the variable $x$.
Following the conventions of Section 3.3 of \cite{mybook},
we emphasise the fact that the pseudodifferential operator $A$ in formula
(\ref{algorithm equation 8}) acts on the oscillatory integral
$U(t)$ in the variable $x$ by writing this pseudodifferential operator
as $A(x,D_x)$, where
$D_{x^\alpha}:=-i\partial/\partial x^\alpha$.

We examine first the dynamic equation (\ref{algorithm equation 8}).
We have
\[
(D_t+A(x,D_x))\,U^{(j)}(t)=F^{(j)}(t)\,,
\]
where $F^{(j)}(t)$ is the oscillatory integral
\[
F^{(j)}(t)
=
\int e^{i\varphi^{(j)}(t,x;y,\eta)}
\,f^{(j)}(t,x;y,\eta)
\,\varsigma^{(j)}(t,x;y,\eta)\,d_{\varphi^{(j)}}(t,x;y,\eta)\,
(\ \cdot\ )\,dy\,{d{\hskip-1pt\bar{}}\hskip1pt}\eta
\]
whose matrix-valued amplitude $f^{(j)}$ is given by the formula
\begin{equation}
\label{algorithm equation 12}
f^{(j)}=D_tu^{(j)}+
\bigl(
\varphi^{(j)}_t+(d_{\varphi^{(j)}})^{-1}(D_t d_{\varphi^{(j)}})+s^{(j)}
\bigr)
\,u^{(j)},
\end{equation}
where the matrix-function $s^{(j)}(t,x;y,\eta)$ is defined as
\begin{equation}
\label{algorithm equation 13}
s^{(j)}=e^{-i\varphi^{(j)}}(d_{\varphi^{(j)}})^{-1}\,A(x,D_x)\,(e^{i\varphi^{(j)}}d_{\varphi^{(j)}})\,.
\end{equation}

Theorem 18.1 from \cite{shubin} gives us the following explicit asymptotic
(in inverse powers of $\eta$) formula for the
matrix-function (\ref{algorithm equation 13}):
\begin{equation}
\label{algorithm equation 14}
s^{(j)}=(d_{\varphi^{(j)}})^{-1}\sum_{\bm\alpha}
\frac1{{\bm\alpha}!}
\,A^{({\bm\alpha})}(x,\varphi^{(j)}_x)\,(D_z^{\bm\alpha}\chi^{(j)})\bigr|_{z=x}\ ,
\end{equation}
where
\begin{equation}
\label{algorithm equation 15}
\chi^{(j)}(t,z,x;y,\eta)
=e^{i\psi^{(j)}(t,z,x;y,\eta)}d_{\varphi^{(j)}}(t,z;y,\eta),
\end{equation}
\begin{equation}
\label{algorithm equation 16}
\psi^{(j)}(t,z,x;y,\eta)
=\varphi^{(j)}(t,z;y,\eta)
-\varphi^{(j)}(t,x;y,\eta)
-\varphi^{(j)}_{x^\beta}(t,x;y,\eta)\,(z-x)^\beta.
\end{equation}
In formula (\ref{algorithm equation 14})
\begin{itemize}
\item
${\bm\alpha}:=(\alpha_1,\ldots,\alpha_n)$ is a multi-index
(note the bold font which we use to distinguish
multi-indices and individual indices),
${\bm\alpha}!:=\alpha_1!\cdots\alpha_n!\,$,
$D_z^{\bm\alpha}:=D_{z^1}^{\alpha_1}\cdots D_{z^n}^{\alpha_n}$,
$D_{z^\beta}:=-i\partial/\partial z^\beta$,
\item
$A(x,\xi)$ is the full symbol of the pseudodifferential operator $A$
written in local coordinates~$x$,
\item
$A^{({\bm\alpha})}(x,\xi):=\partial_\xi^{\bm\alpha}A(x,\xi)$,
$\partial_\xi^{\bm\alpha}:=\partial_{\xi_1}^{\alpha_1}\cdots\partial_{\xi_n}^{\alpha_n}$
and $\partial_{\xi_\beta}:=\partial/\partial\xi_\beta\,$.
\end{itemize}

When $|\eta|\to+\infty$
the matrix-valued amplitude $f^{(j)}(t,x;y,\eta)$ defined by formula
(\ref{algorithm equation 12}) admits an asymptotic expansion
\begin{equation}
\label{algorithm equation 17}
f^{(j)}(t,x;y,\eta)=f^{(j)}_1(t,x;y,\eta)+f^{(j)}_0(t,x;y,\eta)+f^{(j)}_{-1}(t,x;y,\eta)+\ldots
\end{equation}
into components positively homogeneous in $\eta$, with the subscript
indicating degree of homogeneity. Note the following differences between formulae
(\ref{decomposition of symbol of OI into homogeneous components})
and (\ref{algorithm equation 17}).
\begin{itemize}
\item
The leading term in
(\ref{algorithm equation 17})
has degree of homogeneity 1, rather than 0 as in
(\ref{decomposition of symbol of OI into homogeneous components}).
In fact, the leading term in
(\ref{algorithm equation 17})
can be easily written out explicitly
\begin{equation}
\label{algorithm equation 18}
f^{(j)}_1(t,x;y,\eta)=
(\varphi^{(j)}_t(t,x;y,\eta)+A_1(x,\varphi^{(j)}_x(t,x;y,\eta)))\,u^{(j)}_0(t;y,\eta)\,,
\end{equation}
where $A_1(x,\xi)$ is the (matrix-valued) principal symbol of the pseudodifferential
operator $A$.
\item
Unlike the symbol $u^{(j)}(t;y,\eta)$, the amplitude
$f^{(j)}(t,x;y,\eta)$ depends on $x$.
\end{itemize}

We now need to exclude the dependence on $x$ from the amplitude
$f^{(j)}(t,x;y,\eta)$. This can be done by means of the algorithm
described in subsection 2.7.3 of \cite{mybook}.
We outline this algorithm below.

Working in local coordinates, define the matrix-function
$\varphi^{(j)}_{x\eta}$ in accordance with
$(\varphi^{(j)}_{x\eta})_\alpha{}^\beta:=\varphi^{(j)}_{x^\alpha\eta_\beta}$
and then define its inverse $(\varphi^{(j)}_{x\eta})^{-1}$ from the identity
$(\varphi^{(j)})_\alpha{}^\beta[(\varphi^{(j)}_{x\eta})^{-1}]_\beta{}^\gamma:=\delta_\alpha{}^\gamma$.
Define the ``scalar'' first order linear differential operators
\begin{equation}
\label{algorithm equation 19}
L^{(j)}_\alpha:=[(\varphi^{(j)}_{x\eta})^{-1}]_\alpha{}^\beta\,(\partial/\partial x^\beta),
\qquad\alpha=1,\ldots,n.
\end{equation}
Note that the coefficients of these differential operators are functions of the position
variable $x$ and the dual variable $\xi$. It is known, see part 2 of Appendix E in \cite{mybook},
that the operators (\ref{algorithm equation 19}) commute:
$\ L^{(j)}_\alpha L^{(j)}_\beta=L^{(j)}_\beta L^{(j)}_\alpha$,
$\ \alpha,\beta=1,\ldots,n$.

Denote
$\ L^{(j)}_{\bm\alpha}:=(L^{(j)}_1)^{\alpha_1}\cdots(L^{(j)}_n)^{\alpha_n}$,
$\ (-\varphi^{(j)}_\eta)^{\bm\alpha}:=(-\varphi^{(j)}_{\eta_1})^{\alpha_1}\cdots(-\varphi^{(j)}_{\eta_n})^{\alpha_n}$,
and, given an $r\in\mathbb{N}$, define the ``scalar'' linear differential operator
\begin{equation}
\label{algorithm equation 21}
\mathfrak{P}^{(j)}_{-1,r}:=
i(d_{\varphi^{(j)}})^{-1}
\,
\frac\partial{\partial\eta_\beta}
\,d_{\varphi^{(j)}}
\left(1+
\sum_{1\le|{\bm\alpha}|\le2r-1}
\frac{(-\varphi^{(j)}_\eta)^{\bm\alpha}}{{\bm\alpha}!\,(|{\bm\alpha}|+1)}
\,L^{(j)}_{\bm\alpha}
\right)
L^{(j)}_\beta\,,
\end{equation}
where $|{\bm\alpha}|:=\alpha_1+\ldots+\alpha_n$ and the repeated index $\beta$ indicates
summation over $\beta=1,\ldots,n$.

Recall Definition 2.7.8 from \cite{mybook}:
the linear operator $L$ is said to be
positively homogeneous in $\eta$ of degree $p\in\mathbb{R}$
if for any $q\in\mathbb{R}$ and any function $f$
positively homogeneous in $\eta$ of degree $q$
the function $Lf$ is
positively homogeneous in $\eta$ of degree $p+q$.
It is easy to see that the operator (\ref{algorithm equation 21}) is
positively homogeneous in $\eta$ of degree $-1$
and the first subscript in $\mathfrak{P}^{(j)}_{-1,r}$ emphasises this fact.

Let $\mathfrak{S}^{(j)}_0$ be the (linear) operator of restriction to $x=x^{(j)}(t;y,\eta)$,
\begin{equation}
\label{algorithm equation 22}
\mathfrak{S}^{(j)}_0:=\left.(\,\cdot\,)\right|_{x=x^{(j)}(t;y,\eta)}\,,
\end{equation}
and let
\begin{equation}
\label{algorithm equation 23}
\mathfrak{S}^{(j)}_{-r}:=\mathfrak{S}^{(j)}_0(\mathfrak{P}^{(j)}_{-1,r})^r
\end{equation}
for $r=1,2,\ldots$. Observe that our linear operators
$\mathfrak{S}^{(j)}_{-r}$, $r=0,1,2,\ldots$, are
positively homogeneous in $\eta$ of degree $-r$.
This observation allows us to define the linear operator
\begin{equation}
\label{algorithm equation 24}
\mathfrak{S}^{(j)}:=\sum_{r=0}^{+\infty}\mathfrak{S}^{(j)}_{-r}\ ,
\end{equation}
where the series is understood as an asymptotic series in inverse powers of $\eta$.

According to subsection 2.7.3 of \cite{mybook},
the dynamic equation (\ref{algorithm equation 8}) can now be rewritten in the equivalent form
\begin{equation}
\label{algorithm equation 25}
\mathfrak{S}^{(j)}f^{(j)}=0\,,
\end{equation}
where the equality is understood in the asymptotic sense, as
an asymptotic expansion in inverse powers of $\eta$.
Recall that the matrix-valued amplitude $f^{(j)}(t,x;y,\eta)$
appearing in (\ref{algorithm equation 25}) is defined
by formulae (\ref{algorithm equation 12})--(\ref{algorithm equation 16}).

Substituting (\ref{algorithm equation 24}) and (\ref{algorithm equation 17})
into (\ref{algorithm equation 25}) we obtain a hierarchy of equations
\begin{equation}
\label{algorithm equation 26}
\mathfrak{S}^{(j)}_0f^{(j)}_1=0,
\end{equation}
\begin{equation}
\label{algorithm equation 27}
\mathfrak{S}^{(j)}_{-1}f^{(j)}_1+\mathfrak{S}^{(j)}_0f^{(j)}_0=0,
\end{equation}
\[
\mathfrak{S}^{(j)}_{-2}f^{(j)}_1+\mathfrak{S}^{(j)}_{-1}f^{(j)}_0+\mathfrak{S}^{(j)}_0f^{(j)}_{-1}=0,
\]
\[
\ldots
\]
positively homogeneous in $\eta$ of degree 1, 0, $-1$, $\ldots$.
These are the \emph{transport} equations for the determination of the unknown
homogeneous components $u^{(j)}_0(t;y,\eta)$, $u^{(j)}_{-1}(t;y,\eta)$, $u^{(j)}_{-2}(t;y,\eta)$, $\ldots$,
of the symbol of the oscillatory integral (\ref{algorithm equation 1}).

Let us now examine the initial condition (\ref{algorithm equation 9}).
Each operator $U^{(j)}(0)$ is a pseudodifferential operator, only
written in a slightly nonstandard form. The issues here are as follows.
\begin{itemize}
\item
We use the invariantly defined phase function
$
\varphi^{(j)}(0,x;y,\eta)
=(x-y)^\alpha\,\eta_\alpha
+O(|x-y|^2)
$
rather than the linear phase function $(x-y)^\alpha\,\eta_\alpha$
written in local coordinates.
\item
When defining the (full) symbol of the operator $U^{(j)}(t)$ we excluded the variable
$x$ from the amplitude rather than the variable $y$. Note that when dealing
with pseudodifferential operators it is customary to exclude the variable $y$
from the amplitude; exclusion of the variable $x$ gives the dual symbol of
a pseudodifferential operator, see subsection 2.1.3 in \cite{mybook}.
Thus, at $t=0$, our symbol $u^{(j)}(0;y,\eta)$ resembles
the dual symbol of a pseudodifferential operator rather
than the ``normal'' symbol.
\item
We have the extra factor $d_{\varphi^{(j)}}(0,x;y,\eta)$ in our representation
of the operator $U^{(j)}(0)$ as an oscillatory integral.
\end{itemize}

The (full) dual symbol
of the pseudodifferential operator $U^{(j)}(0)$
can be calculated in local coordinates in accordance with the following
formula which addresses the issues highlighted above:
\begin{equation}
\label{algorithm equation 30}
\sum_{\bm\alpha}
\frac{(-1)^{|{\bm\alpha}|}}{{\bm\alpha}!}\,
\bigl(
D_x^{\bm\alpha}\,\partial_\eta^{\bm\alpha}\,
u^{(j)}(0;y,\eta)\,
e^{i\omega^{(j)}(x;y,\eta)}\,d_{\varphi^{(j)}}(0,x;y,\eta)
\bigr)
\bigr|_{x=y}\ ,
\end{equation}
where
$\omega^{(j)}(x;y,\eta)=\varphi^{(j)}(0,x;y,\eta)-(x-y)^\beta\,\eta_\beta\,$.
Formula (\ref{algorithm equation 30})
is a version of the formula from subsection 2.1.3 of \cite{mybook}, only with
the extra factor $(-1)^{|{\bm\alpha}|}$. The latter is needed because we are writing
down the dual symbol of the pseudodifferential operator $U^{(j)}(0)$ (no dependence on $x$)
rather than its ``normal'' symbol (no dependence on $y$).

The initial condition (\ref{algorithm equation 9}) can now be rewritten in explicit form as
\begin{equation}
\label{algorithm equation 32}
\sum_j
\sum_{\bm\alpha}
\frac{(-1)^{|{\bm\alpha}|}}{{\bm\alpha}!}\,
\bigl(
D_x^{\bm\alpha}\,\partial_\eta^{\bm\alpha}\,
u^{(j)}(0;y,\eta)\,
e^{i\omega^{(j)}(x;y,\eta)}\,d_{\varphi^{(j)}}(0,x;y,\eta)
\bigr)
\bigr|_{x=y}=I\,,
\end{equation}
where $I$ is the $m\times m$ identity matrix.
Condition (\ref{algorithm equation 32})
can be decomposed into components positively homogeneous in $\eta$
of degree $0,-1,-2,\ldots$, giving us a hierarchy of initial conditions.
The leading (of degree of homogeneity 0) initial condition reads
\begin{equation}
\label{algorithm equation 33}
\sum_j
u^{(j)}_0(0;y,\eta)=I\,,
\end{equation}
whereas lower order initial conditions are more complicated
and depend on the choice of our phase functions $\varphi^{(j)}$.

\section{Leading transport equations}
\label{Leading transport equations}

Formulae
(\ref{algorithm equation 22}),
(\ref{algorithm equation 18}),
(\ref{algorithm equation 2}),
(\ref{Hamiltonian system of equations})
and the identity $\xi_\alpha h^{(j)}_{\xi_\alpha}(x,\xi)=h^{(j)}(x,\xi)$
(consequence of the fact that $h^{(j)}(x,\xi)$ is positively homogeneous in $\xi$ of degree~1)
give us the following explicit representation
for the leading transport equation
(\ref{algorithm equation 26}):
\begin{equation}
\label{Leading transport equations equation 1}
\!\!
\bigl[
A_1\bigl(x^{(j)}(t;y,\eta),\xi^{(j)}(t;y,\eta)\bigr)
-
h^{(j)}\bigl(x^{(j)}(t;y,\eta),\xi^{(j)}(t;y,\eta)\bigr)
\bigr]
\,u^{(j)}_0(t;y,\eta)=0.
\end{equation}
Here, of course,
$h^{(j)}\bigl(x^{(j)}(t;y,\eta),\xi^{(j)}(t;y,\eta)\bigr)=h^{(j)}(y,\eta)$.

Equation (\ref{Leading transport equations equation 1}) implies that
\begin{equation}
\label{Leading transport equations equation 2}
u^{(j)}_0(t;y,\eta)=v^{(j)}(x^{(j)}(t;y,\eta),\xi^{(j)}(t;y,\eta))
\,[w^{(j)}(t;y,\eta)]^T,
\end{equation}
where $v^{(j)}(z,\zeta)$ is the normalised eigenvector of the principal
symbol $A_1(z,\zeta)$ corresponding to the eigenvalue $h^{(j)}(z,\zeta)$
and $w^{(j)}:\mathbb{R}\times T'M\to\mathbb{C}^m$ is a column-function,
positively homogeneous in $\eta$ of degree 0, that remains to be found.
Formulae
(\ref{algorithm equation 33})
and
(\ref{Leading transport equations equation 2})
imply the following initial condition for the unknown column-function $w^{(j)}$:
\begin{equation}
\label{Leading transport equations equation 3}
w^{(j)}(0;y,\eta)=\overline{v^{(j)}(y,\eta)}.
\end{equation}

We now consider the
next transport equation in our hierarchy,
equation (\ref{algorithm equation 27}).
We will write down the two terms appearing in
(\ref{algorithm equation 27}) separately.

In view of formulae
(\ref{algorithm equation 18})
and
(\ref{algorithm equation 21})--(\ref{algorithm equation 23}),
the first term in (\ref{algorithm equation 27}) reads
\begin{multline}
\label{Leading transport equations equation 4}
\mathfrak{S}^{(j)}_{-1}f^{(j)}_1=
\\
i
\left.
\left[
(d_{\varphi^{(j)}})^{-1}
\frac\partial{\partial\eta_\beta}
d_{\varphi^{(j)}}
\left(1-
\frac12
\varphi^{(j)}_{\eta_\alpha}
L^{(j)}_\alpha
\right)
\left(
L^{(j)}_\beta
\bigl(\varphi^{(j)}_t+A_1(x,\varphi^{(j)}_x)\bigr)
\right)
u^{(j)}_0
\right]
\right|_{x=x^{(j)}}\,,
\end{multline}
where we dropped, for the sake of brevity,
the arguments $(t;y,\eta)$ in $u^{(j)}_0$ and $x^{(j)}$,
and the arguments $(t,x;y,\eta)$
in $\varphi^{(j)}_t$, $\varphi^{(j)}_x$, $\varphi^{(j)}_\eta$ and $d_{\varphi^{(j)}}\,$.
Recall that the differential operators $L^{(j)}_\alpha$ are defined in accordance with
formula (\ref{algorithm equation 19})
and the coefficients of these operators depend on $(t,x;y,\eta)$.

In view of formulae
(\ref{algorithm equation 12})--(\ref{algorithm equation 17})
and
(\ref{algorithm equation 22}),
the second term in (\ref{algorithm equation 27}) reads
\begin{multline}
\label{Leading transport equations equation 5}
\mathfrak{S}^{(j)}_0f^{(j)}_0=
D_tu^{(j)}_0
\\
+\left.\left[
(d_{\varphi^{(j)}})^{-1}
\left(D_t+(A_1)_{\xi_\alpha}D_{x^\alpha}\right)
d_{\varphi^{(j)}}
+A_0
-\frac i2(A_1)_{\xi_\alpha\xi_\beta}C^{(j)}_{\alpha\beta}
\right]
\right|_{x=x^{(j)}}u^{(j)}_0
\\
+\bigl[A_1-h^{(j)}\bigr]u^{(j)}_{-1}\,,
\end{multline}
where
\begin{equation}
\label{Leading transport equations equation 6}
C^{(j)}_{\alpha\beta}:=\left.\varphi^{(j)}_{x^\alpha x^\beta}\right|_{x=x^{(j)}}
\end{equation}
is the matrix-function from
(\ref{algorithm equation 5}).
In formulae
(\ref{Leading transport equations equation 5})
and
(\ref{Leading transport equations equation 6})
we dropped, for the sake of brevity,
the arguments $(t;y,\eta)$ in $u^{(j)}_0$, $u^{(j)}_{-1}$, $C^{(j)}_{\alpha\beta}$ and $x^{(j)}$,
the arguments
$(x^{(j)}(t;y,\eta),\xi^{(j)}(t;y,\eta))$
in $A_0$, $A_1$, $(A_1)_{\xi_\alpha}$, $(A_1)_{\xi_\alpha\xi_\beta}$ and $h^{(j)}$,
and the arguments $(t,x;y,\eta)$
in $d_{\varphi^{(j)}}$ and $\varphi^{(j)}_{x^\alpha x^\beta}\,$.

Looking at
(\ref{Leading transport equations equation 4})
and
(\ref{Leading transport equations equation 5})
we see that the transport equation (\ref{algorithm equation 27}) has a complicated
structure.
Hence, in this section we choose not to perform the analysis
of the full equation (\ref{algorithm equation 27})
and analyse only one particular subequation of this equation.
Namely, observe that equation (\ref{algorithm equation 27})
is equivalent to $m$ subequations
\begin{equation}
\label{Leading transport equations equation 7}
\bigl[v^{(j)}\bigr]^*
\,
\bigl[\mathfrak{S}^{(j)}_{-1}f^{(j)}_1+\mathfrak{S}^{(j)}_0f^{(j)}_0\bigr]
=0,
\end{equation}
\begin{equation}
\label{Leading transport equations equation 8}
\bigl[v^{(l)}\bigr]^*
\,
\bigl[\mathfrak{S}^{(j)}_{-1}f^{(j)}_1+\mathfrak{S}^{(j)}_0f^{(j)}_0\bigr]
=0,
\qquad l\ne j,
\end{equation}
where we dropped, for the sake of brevity, the arguments
$(x^{(j)}(t;y,\eta),\xi^{(j)}(t;y,\eta))$
in $\bigl[v^{(j)}\bigr]^*$ and $\bigl[v^{(l)}\bigr]^*$.
In the remainder of this section we analyse (sub)equation
(\ref{Leading transport equations equation 7}) only.

Equation (\ref{Leading transport equations equation 7})
is simpler than each of the $m-1$ equations
(\ref{Leading transport equations equation 8})
for the following two reasons.

\begin{itemize}

\item
Firstly, the term
$\bigl[A_1-h^{(j)}\bigr]u^{(j)}_{-1}$
from
(\ref{Leading transport equations equation 5})
vanishes after multiplication by
$\bigl[v^{(j)}\bigr]^*$
from the left.
Hence, equation
(\ref{Leading transport equations equation 7})
does not contain $u^{(j)}_{-1}$.

\item
Secondly, if we substitute
(\ref{Leading transport equations equation 2})
into
(\ref{Leading transport equations equation 7}),
then the term with
\[
\partial[d_{\varphi^{(j)}}w^{(j)}(t;y,\eta)]^T/\partial\eta_\beta
\]
vanishes.
This follows from the fact that the scalar function
\[
\bigl[v^{(j)}\bigr]^*
\bigl(\varphi^{(j)}_t+A_1(x,\varphi^{(j)}_x)\bigr)
v^{(j)}
\]
has a second order zero, in the variable $x$, at $x=x^{(j)}(t;y,\eta)$.
Indeed, we have
\begin{multline*}
\left.
\left[
\frac\partial{\partial x^\alpha}
\bigl[v^{(j)}\bigr]^*
\bigl(\varphi^{(j)}_t+A_1(x,\varphi^{(j)}_x)\bigr)
v^{(j)}
\right]
\right|_{x=x^{(j)}}
\\
=
\bigl[v^{(j)}\bigr]^*
\left.
\left[
\bigl(\varphi^{(j)}_t+A_1(x,\varphi^{(j)}_x)\bigr)_{x^\alpha}
\right]
\right|_{x=x^{(j)}}
v^{(j)}
\\
=
\bigl[v^{(j)}\bigr]^*
\bigl(
-h^{(j)}_{x^\alpha}-C^{(j)}_{\alpha\beta}h^{(j)}_{\xi_\beta}
+(A_1)_{x^\alpha}+C^{(j)}_{\alpha\beta}(A_1)_{\xi_\beta}
\bigr)
v^{(j)}
\\
=
\bigl[v^{(j)}\bigr]^*(A_1)_{x^\alpha}v^{(j)}-h^{(j)}_{x^\alpha}
+
C^{(j)}_{\alpha\beta}
\bigl(
\bigl[v^{(j)}\bigr]^*(A_1)_{\xi_\beta}v^{(j)}-h^{(j)}_{\xi_\beta}
\bigr)
=0\,,
\end{multline*}
where in the last two lines we dropped,
for the sake of brevity,
the arguments
$(x^{(j)}(t;y,\eta),\xi^{(j)}(t;y,\eta))$
in $(A_1)_{x^\alpha}$, $(A_1)_{\xi_\beta}$,
$h^{(j)}_{x^\alpha}$, $h^{(j)}_{\xi_\beta}$,
and the argument $(t;y,\eta)$ in
$C^{(j)}_{\alpha\beta}$
(the latter is the
matrix-function
from
formulae (\ref{algorithm equation 5}) and (\ref{Leading transport equations equation 6})).
Throughout the above argument we used the fact that our
$\bigl[v^{(j)}\bigr]^*$ and $v^{(j)}$ do not depend on $x$:
their argument is
$(x^{(j)}(t;y,\eta),\xi^{(j)}(t;y,\eta))$.

\end{itemize}

Substituting
(\ref{Leading transport equations equation 4}),
(\ref{Leading transport equations equation 5})
and
(\ref{Leading transport equations equation 2})
into
(\ref{Leading transport equations equation 7})
we get
\begin{equation}
\label{Leading transport equations equation 9}
(D_t+p^{(j)}(t;y,\eta))\,[w^{(j)}(t;y,\eta)]^T=0\,,
\end{equation}
where
\begin{multline}
\label{Leading transport equations equation 10}
p^{(j)}=
i
\left.
[v^{(j)}]^*
\left[
\frac\partial{\partial\eta_\beta}
\left(1-
\frac12
\varphi^{(j)}_{\eta_\alpha}
L^{(j)}_\alpha
\right)
\left(
L^{(j)}_\beta
\bigl(\varphi^{(j)}_t+A_1(x,\varphi^{(j)}_x)\bigr)
\right)
v^{(j)}
\right]
\right|_{x=x^{(j)}}
\\
-i[v^{(j)}]^*\{v^{(j)},h^{(j)}\}
+\left.\left[
(d_{\varphi^{(j)}})^{-1}
\left(D_t+h^{(j)}_{\xi_\alpha}D_{x^\alpha}\right)
d_{\varphi^{(j)}}
\right]
\right|_{x=x^{(j)}}
\\
+[v^{(j)}]^*
\left(
A_0
-\frac i2(A_1)_{\xi_\alpha\xi_\beta}C^{(j)}_{\alpha\beta}
\right)
v^{(j)}.
\end{multline}

Note that the ordinary differential operator in the LHS of formula
(\ref{Leading transport equations equation 9}) is a scalar one, i.e.
it does not mix up the different components of the column-function
$w^{(j)}(t;y,\eta)$. The solution of the ordinary
differential equation
(\ref{Leading transport equations equation 9})
subject to the initial condition
(\ref{Leading transport equations equation 3})
is
\begin{equation}
\label{Leading transport equations equation 11}
w^{(j)}(t;y,\eta)=\overline{v^{(j)}(y,\eta)}
\exp\left(-i\int_0^tp^{(j)}(\tau;y,\eta)\,d\tau\right).
\end{equation}
Comparing formulae
(\ref{Leading transport equations equation 2}),
(\ref{Leading transport equations equation 11})
with formula
(\ref{formula for principal symbol of oscillatory integral})
we see that in order to prove the latter we need only to establish the
scalar identity
\begin{equation}
\label{Leading transport equations equation 12}
p^{(j)}(t;y,\eta)=q^{(j)}(x^{(j)}(t;y,\eta),\xi^{(j)}(t;y,\eta))\,,
\end{equation}
where $q^{(j)}$ is the function
(\ref{phase appearing in principal symbol}).
In view of the definitions of the quantities
$p^{(j)}$ and $q^{(j)}$,
see formulae
(\ref{Leading transport equations equation 10})
and
(\ref{phase appearing in principal symbol}),
and the definition of the subprincipal symbol
(\ref{definition of subprincipal symbol}),
proving the identity
(\ref{Leading transport equations equation 12})
reduces to proving the identity
\begin{multline}
\label{Leading transport equations equation 13}
\{
[v^{(j)}]^*,A_1-h^{(j)},v^{(j)}
\}
(x^{(j)},\xi^{(j)})
=
\\
-2
\left.
[v^{(j)}(x^{(j)},\xi^{(j)})]^*
\left[
\frac\partial{\partial\eta_\beta}
\left(1-
\frac12
\varphi^{(j)}_{\eta_\alpha}
L^{(j)}_\alpha
\right)
\left(
L^{(j)}_\beta
\bigl(\varphi^{(j)}_t+A_1(x,\varphi^{(j)}_x)\bigr)
\right)
v^{(j)}(x^{(j)},\xi^{(j)})
\right]
\right|_{x=x^{(j)}}
\\
+2\left.\left[
(d_{\varphi^{(j)}})^{-1}
\left(\partial_t+h^{(j)}_{\xi_\alpha}\partial_{x^\alpha}\right)
d_{\varphi^{(j)}}
\right]
\right|_{x=x^{(j)}}
\\
+[v^{(j)}(x^{(j)},\xi^{(j)})]^*
\left(
(A_1)_{x^\alpha\xi_\alpha}+(A_1)_{\xi_\alpha\xi_\beta}C^{(j)}_{\alpha\beta}
\right)
v^{(j)}(x^{(j)},\xi^{(j)}).
\end{multline}
Note that the expressions in the LHS and RHS of
(\ref{Leading transport equations equation 13}) have different
structure. The LHS of
(\ref{Leading transport equations equation 13})
is the generalised Poisson bracket
$\{
[v^{(j)}]^*,A_1-h^{(j)},v^{(j)}
\}$,
see (\ref{generalised Poisson bracket on matrix-functions}),
evaluated at $z=x^{(j)}(t;y,\eta)$, $\zeta=\xi^{(j)}(t;y,\eta)$,
whereas the RHS of
(\ref{Leading transport equations equation 13})
involves partial derivatives (in $\eta$) of
$v^{(j)}(x^{(j)}(t;y,\eta),\xi^{(j)}(t;y,\eta))$
(Chain Rule).
In writing
(\ref{Leading transport equations equation 13}) we also
dropped, for the sake of brevity,
the arguments $(t,x;y,\eta)$ in
$\varphi^{(j)}_t$, $\varphi^{(j)}_x$,
$\varphi^{(j)}_\eta$, $d_{\varphi^{(j)}}\,$
and the coefficients of the differential operators
$L^{(j)}_\alpha$ and $L^{(j)}_\beta$,
the arguments $(x^{(j)},\xi^{(j)})$
in $h^{(j)}_{\xi_\alpha}$,
$(A_1)_{x^\alpha\xi_\alpha}$ and $(A_1)_{\xi_\alpha\xi_\beta}$,
and the arguments $(t;y,\eta)$
in $x^{(j)}$, $\xi^{(j)}$ and $C^{(j)}_{\alpha\beta}$.

Before performing the calculations that will establish the identity
(\ref{Leading transport equations equation 13}) we make several
observations that will allow us to simplify these calculations
considerably.

Firstly, our function $p^{(j)}(t;y,\eta)$ does not depend on the choice
of the phase function $\varphi^{(j)}(t,x;y,\eta)$. Indeed, if
$p^{(j)}(t;y,\eta)$ did depend on the choice of phase function, then,
in view of formulae
(\ref{Leading transport equations equation 2})
and
(\ref{Leading transport equations equation 11})
the principal symbol of our oscillatory integral $U^{(j)}(t)$ would depend
on the choice of phase function, which would contradict Theorem
2.7.11 from \cite{mybook}. Here we use the fact that operators
$U^{(j)}(t)$ with different $j$ cannot compensate each other to give
an integral operator whose integral kernel is infinitely
smooth in $t$, $x$ and $y$ because all our $U^{(j)}(t)$ oscillate in
$t$ in a different way:
$\varphi^{(j)}_t(t,x^{(j)}(t;y,\eta);y,\eta)=-h^{(j)}(y,\eta)$
and we assumed the eigenvalues $h^{(j)}(y,\eta)$ of our principal
symbol $A_1(y,\eta)$ to be simple.

Secondly, the arguments (free variables) in
(\ref{Leading transport equations equation 13}) are
$(t;y,\eta)$. We fix an arbitrary point
$(\tilde t;\tilde y,\tilde\eta)\in\mathbb{R}\times T'M$
and prove formula
(\ref{Leading transport equations equation 13})
at this point.
Put
$(\xi^{(j)}_\eta)_\alpha{}^\beta:=\partial(\xi^{(j)})_\alpha/\partial\eta_\beta$.
According to Lemma 2.3.2 from \cite{mybook}
there exists a local coordinate system $x$ such that
$\det(\xi^{(j)}_\eta)_\alpha{}^\beta\ne0$.
This opens the way to the use of the linear phase function
\begin{equation}
\label{Leading transport equations equation 14}
\varphi^{(j)}(t,x;y,\eta)
=(x-x^{(j)}(t;y,\eta))^\alpha\,\xi^{(j)}_\alpha(t;y,\eta)
\end{equation}
which will simplify calculations to a great extent.
Moreover, we can choose a local coordinate system $y$ such that
\begin{equation}
\label{Leading transport equations equation 15}
(\xi^{(j)}_\eta)_\alpha{}^\beta(\tilde t;\tilde y,\tilde\eta)=\delta_\alpha{}^\beta
\end{equation}
which will simplify calculations even further.

The calculations we are about to perform will make use of the symmetry
\begin{equation}
\label{Leading transport equations equation 16}
(x^{(j)}_\eta)^{\gamma\alpha}(\xi^{(j)}_\eta)_\gamma{}^\beta
=
(x^{(j)}_\eta)^{\gamma\beta}(\xi^{(j)}_\eta)_\gamma{}^\alpha
\end{equation}
which is an immediate
consequence of formula (\ref{algorithm equation 7.1.5}).
Formula (\ref{Leading transport equations equation 16}) appears
as formula (2.3.3) in \cite{mybook} and the accompanying text
explains its geometric meaning. Note that at the point
$(\tilde t;\tilde y,\tilde\eta)$ formula
(\ref{Leading transport equations equation 16}) takes the
especially simple form
\begin{equation}
\label{Leading transport equations equation 17}
(x^{(j)}_\eta)^{\alpha\beta}(\tilde t;\tilde y,\tilde\eta)
=
(x^{(j)}_\eta)^{\beta\alpha}(\tilde t;\tilde y,\tilde\eta).
\end{equation}

Our calculations will also involve the quantity
$\varphi^{(j)}_{\eta_\alpha\eta_\beta}(\tilde t,\tilde x;\tilde y,\tilde\eta)$
where $\tilde x:=x^{(j)}(\tilde t;\tilde y,\tilde\eta)$.
Formulae
(\ref{Leading transport equations equation 14}),
(\ref{algorithm equation 7.1.5}),
(\ref{Leading transport equations equation 15})
and
(\ref{Leading transport equations equation 17})
imply
\begin{equation}
\label{Leading transport equations equation 18}
\varphi^{(j)}_{\eta_\alpha\eta_\beta}(\tilde t,\tilde x;\tilde y,\tilde\eta)
=
-(x^{(j)}_\eta)^{\alpha\beta}(\tilde t;\tilde y,\tilde\eta).
\end{equation}

Further on we denote $\tilde\xi:=\xi^{(j)}(\tilde t;\tilde y,\tilde\eta)$.

With account of all the simplifications listed above, we can rewrite
formula
(\ref{Leading transport equations equation 13}),
which is the identity that we are proving, as
\begin{multline}
\label{Leading transport equations equation 19}
\{
[v^{(j)}]^*,A_1-h^{(j)},v^{(j)}
\}
(\tilde x,\tilde\xi)
=
\\
\!\!\!\!\!\!\!\!\!\!\!\!\!\!\!\!\!\!\!\!\!\!\!\!\!\!\!\!
-2
[\tilde v^{(j)}]^*
\Bigl[
\frac{\partial^2}{\partial x^\alpha\partial\eta_\alpha}
\bigl(A_1(x,\xi^{(j)})-h^{(j)}(\tilde y,\eta)
\\
\qquad\qquad\qquad\qquad\qquad\qquad
-(x-x^{(j)})^\gamma h^{(j)}_{x^\gamma}(x^{(j)},\xi^{(j)})\bigr)
\,v^{(j)}(x^{(j)},\xi^{(j)})
\Bigr]
\Bigr|_{(x,\eta)=(\tilde x,\tilde\eta)}
\\
\!\!\!\!\!\!\!\!\!\!\!\!\!\!\!\!\!\!\!\!\!\!\!\!\!\!\!\!\!\!\!\!\!\!\!\!\!\!\!\!\!\!\!\!\!\!\!
-(\tilde x^{(j)}_\eta)^{\alpha\beta}\,
[\tilde v^{(j)}]^*
\Bigl[
\frac{\partial^2}{\partial x^\alpha\partial x^\beta}
\bigl(A_1(x,\xi^{(j)})-h^{(j)}(\tilde y,\eta)
\\
\qquad\qquad\qquad\qquad\qquad\qquad
-(x-x^{(j)})^\gamma h^{(j)}_{x^\gamma}(x^{(j)},\xi^{(j)})\bigr)
\,v^{(j)}(x^{(j)},\xi^{(j)})
\Bigr]
\Bigr|_{(x,\eta)=(\tilde x,\tilde\eta)}
\\
+[\tilde v^{(j)}]^*
(\tilde A_1)_{x^\alpha\xi_\alpha}
\tilde v^{(j)}
-\tilde h^{(j)}_{x^\alpha\xi_\alpha}
-\tilde h^{(j)}_{x^\alpha x^\beta}(\tilde
x^{(j)}_\eta)^{\alpha\beta}\,,
\qquad
\end{multline}

\noindent
where
$\tilde v^{(j)}=v^{(j)}(\tilde x,\tilde\xi)$,
$\tilde x^{(j)}_\eta=x^{(j)}_\eta(\tilde t;\tilde y,\tilde\eta)$,
$(\tilde A_1)_{x^\alpha\xi_\alpha}=(A_1)_{x^\alpha\xi_\alpha}(\tilde x,\tilde\xi)$,
$\tilde h^{(j)}_{x^\alpha\xi_\alpha}=h^{(j)}_{x^\alpha\xi_\alpha}(\tilde x,\tilde\xi)$,
$\tilde h^{(j)}_{x^\alpha x^\beta}=h^{(j)}_{x^\alpha x^\beta}(\tilde x,\tilde\xi)$,
$x^{(j)}=x^{(j)}(\tilde t;\tilde y,\eta)$
and
$\xi^{(j)}=\xi^{(j)}(\tilde t;\tilde y,\eta)$.

Note that the last two terms in the RHS of
(\ref{Leading transport equations equation 19})
originate from the term with
$d_{\varphi^{(j)}}$ in
(\ref{Leading transport equations equation 13}):
we used the fact that $d_{\varphi^{(j)}}$ does not depend on
$x$ and that
\begin{equation}
\label{Leading transport equations equation 20}
\left.
\left[
(d_{\varphi^{(j)}})^{-1}
\partial_t
d_{\varphi^{(j)}}
\right]
\right|_{(t,x;y,\eta)=(\tilde t,\tilde x;\tilde y,\tilde\eta)}
=-\frac12
\bigl(
\tilde h^{(j)}_{x^\alpha\xi_\alpha}
+\tilde h^{(j)}_{x^\alpha x^\beta}(\tilde x^{(j)}_\eta)^{\alpha\beta}
\bigr).
\end{equation}
Formula (\ref{Leading transport equations equation 20})
is a special case of formula (3.3.21) from \cite{mybook}.

Note also that the term $-h^{(j)}(\tilde y,\eta)$
appearing (twice) in the RHS of
(\ref{Leading transport equations equation 19})
will vanish after being acted upon with
the differential operators
$\frac{\partial^2}{\partial x^\alpha\partial\eta_\alpha}$
and
$\frac{\partial^2}{\partial x^\alpha\partial x^\beta}$
because it does not depend on $x$.

We have
\begin{multline}
\label{Leading transport equations equation 21}
[\tilde v^{(j)}]^*
\left.
\left[
\frac{\partial^2}{\partial x^\alpha\partial\eta_\alpha}
\bigl(A_1(x,\xi^{(j)})-(x-x^{(j)})^\gamma h^{(j)}_{x^\gamma}(x^{(j)},\xi^{(j)})\bigr)
\,v^{(j)}(x^{(j)},\xi^{(j)})
\right]
\right|_{(x,\eta)=(\tilde x,\tilde\eta)}
\\
=
[\tilde v^{(j)}]^*
(\tilde A_1)_{x^\alpha\xi_\alpha}
\tilde v^{(j)}
-\tilde h^{(j)}_{x^\alpha\xi_\alpha}
-\tilde h^{(j)}_{x^\alpha x^\beta}(\tilde x^{(j)}_\eta)^{\alpha\beta}
\\
+
[\tilde v^{(j)}]^*
\bigl(
(\tilde A_1)_{x^\alpha}
-\tilde h^{(j)}_{x^\alpha}
\bigr)
\bigl(
\tilde v^{(j)}_{\xi_\alpha}
+\tilde v^{(j)}_{x^\beta}(\tilde x^{(j)}_\eta)^{\alpha\beta}
\bigr),
\end{multline}
\begin{multline}
\label{Leading transport equations equation 22}
[\tilde v^{(j)}]^*
\left.
\left[
\frac{\partial^2}{\partial x^\alpha\partial x^\beta}
\bigl(A_1(x,\xi^{(j)})-(x-x^{(j)})^\gamma h^{(j)}_{x^\gamma}(x^{(j)},\xi^{(j)})\bigr)
\,v^{(j)}(x^{(j)},\xi^{(j)})
\right]
\right|_{(x,\eta)=(\tilde x,\tilde\eta)}
\\
=
[\tilde v^{(j)}]^*
(\tilde A_1)_{x^\alpha x^\beta}
\tilde v^{(j)}\,,
\end{multline}
where
$(\tilde A_1)_{x^\alpha}=(A_1)_{x^\alpha}(\tilde x,\tilde\xi)$,
$\tilde h^{(j)}_{x^\alpha}=h^{(j)}_{x^\alpha}(\tilde x,\tilde\xi)$,
$\tilde v^{(j)}_{\xi_\alpha}=v^{(j)}_{\xi_\alpha}(\tilde x,\tilde\xi)$
and
$\tilde v^{(j)}_{x^\beta}=v^{(j)}_{x^\beta}(\tilde x,\tilde\xi)$.
We also have
\begin{multline}
\label{Leading transport equations equation 23}
[\tilde v^{(j)}]^*
\bigl(
(\tilde A_1)_{x^\alpha}
-\tilde h^{(j)}_{x^\alpha}
\bigr)
\tilde v^{(j)}_{x^\beta}
+
[\tilde v^{(j)}]^*
\bigl(
(\tilde A_1)_{x^\beta}
-\tilde h^{(j)}_{x^\beta}
\bigr)
\tilde v^{(j)}_{x^\alpha}
\\
=
\tilde h^{(j)}_{x^\alpha x^\beta}
-
[\tilde v^{(j)}]^*
(\tilde A_1)_{x^\alpha x^\beta}
\tilde v^{(j)}.
\end{multline}
Using formulae
(\ref{Leading transport equations equation 23})
and
(\ref{Leading transport equations equation 17})
we can rewrite formula
(\ref{Leading transport equations equation 21})
as
\begin{multline}
\label{Leading transport equations equation 24}
[\tilde v^{(j)}]^*
\left.
\left[
\frac{\partial^2}{\partial x^\alpha\partial\eta_\alpha}
\bigl(A_1(x,\xi^{(j)})-(x-x^{(j)})^\gamma h^{(j)}_{x^\gamma}(x^{(j)},\xi^{(j)})\bigr)
\,v^{(j)}(x^{(j)},\xi^{(j)})
\right]
\right|_{(x,\eta)=(\tilde x,\tilde\eta)}
\\
=
[\tilde v^{(j)}]^*
(\tilde A_1)_{x^\alpha\xi_\alpha}
\tilde v^{(j)}
-\tilde h^{(j)}_{x^\alpha\xi_\alpha}
+
[\tilde v^{(j)}]^*
\bigl(
(\tilde A_1)_{x^\alpha}
-\tilde h^{(j)}_{x^\alpha}
\bigr)
\tilde v^{(j)}_{\xi_\alpha}
\\
-\frac12
\bigl(
[\tilde v^{(j)}]^*
(\tilde A_1)_{x^\alpha x^\beta}
\tilde v^{(j)}
+
\tilde h^{(j)}_{x^\alpha x^\beta}
\bigr)
(\tilde x^{(j)}_\eta)^{\alpha\beta}.
\end{multline}
Substituting
(\ref{Leading transport equations equation 24})
and
(\ref{Leading transport equations equation 22})
into
(\ref{Leading transport equations equation 19})
we see that all the terms with $(\tilde x^{(j)}_\eta)^{\alpha\beta}$ cancel out
and we get
\begin{multline}
\label{Leading transport equations equation 25}
\{
[v^{(j)}]^*,A_1-h^{(j)},v^{(j)}
\}
(\tilde x,\tilde\xi)
=
\\
-[\tilde v^{(j)}]^*
\bigl(
(\tilde A_1)_{x^\alpha\xi_\alpha}
-
\tilde h^{(j)}_{x^\alpha\xi_\alpha}
\bigr)
\tilde v^{(j)}
-2
[\tilde v^{(j)}]^*
\bigl(
(\tilde A_1)_{x^\alpha}
-\tilde h^{(j)}_{x^\alpha}
\bigr)
\tilde v^{(j)}_{\xi_\alpha}.
\end{multline}
Thus, the proof of the identity
(\ref{Leading transport equations equation 13})
has been reduced to the proof of the identity~(\ref{Leading transport equations equation 25}).

Observe now that formula
(\ref{Leading transport equations equation 25})
no longer has Hamiltonian trajectories present in it.
This means that we can drop all the tildes and
rewrite
(\ref{Leading transport equations equation 25})
as
\begin{multline}
\label{Leading transport equations equation 26}
\{
[v^{(j)}]^*,A_1-h^{(j)},v^{(j)}
\}
=
\\
-[v^{(j)}]^*
\bigl(
A_1
-h^{(j)}
\bigr)_{x^\alpha\xi_\alpha}
v^{(j)}
-2
[v^{(j)}]^*
\bigl(
A_1
-h^{(j)}
\bigr)_{x^\alpha}
v^{(j)}_{\xi_\alpha}\,,
\end{multline}
where the arguments are $(x,\xi)$.
We no longer need to restrict our consideration to the particular
point $(x,\xi)=(\tilde x,\tilde\xi)$:
if we prove
(\ref{Leading transport equations equation 26})
for an arbitrary $(x,\xi)\in T'M$
we will prove it for a particular
$(\tilde x,\tilde\xi)\in T'M$.

The proof of the identity
(\ref{Leading transport equations equation 26})
is straightforward. We note that
\begin{multline}
\label{Leading transport equations equation 27}
[v^{(j)}]^*
(A_1-h^{(j)})_{x^\alpha\xi_\alpha}
v^{(j)}=
\\
-
[v^{(j)}]^*
(A_1-h^{(j)})_{x^\alpha}
v^{(j)}_{\xi_\alpha}
-
[v^{(j)}]^*
(A_1-h^{(j)})_{\xi_\alpha}
v^{(j)}_{x^\alpha}
\end{multline}
and substituting
(\ref{Leading transport equations equation 27})
into
(\ref{Leading transport equations equation 26})
reduce the latter to the form
\begin{multline}
\label{Leading transport equations equation 28}
\{
[v^{(j)}]^*,A_1-h^{(j)},v^{(j)}
\}
=
\\
[v^{(j)}]^*
\bigl(
A_1
-h^{(j)}
\bigr)_{\xi_\alpha}
v^{(j)}_{x^\alpha}
-
[v^{(j)}]^*
\bigl(
A_1
-h^{(j)}
\bigr)_{x^\alpha}
v^{(j)}_{\xi_\alpha}.
\end{multline}
But
\begin{equation}
\label{Leading transport equations equation 29}
[v^{(j)}]^*
\bigl(
A_1
-h^{(j)}
\bigr)_{x^\alpha}
=
-
[v^{(j)}_{x^\alpha}]^*
\bigl(
A_1
-h^{(j)}
\bigr),
\end{equation}
\begin{equation}
\label{Leading transport equations equation 30}
[v^{(j)}]^*
\bigl(
A_1
-h^{(j)}
\bigr)_{\xi_\alpha}
=
-
[v^{(j)}_{\xi_\alpha}]^*
\bigl(
A_1
-h^{(j)}
\bigr).
\end{equation}
Substituting
(\ref{Leading transport equations equation 29})
and
(\ref{Leading transport equations equation 30})
into
(\ref{Leading transport equations equation 28})
we get
\[
\{
[v^{(j)}]^*,A_1-h^{(j)},v^{(j)}
\}
=
[v^{(j)}_{x^\alpha}]^*
\bigl(
A_1
-h^{(j)}
\bigr)
v^{(j)}_{\xi_\alpha}
-
[v^{(j)}_{\xi_\alpha}]^*
\bigl(
A_1
-h^{(j)}
\bigr)
v^{(j)}_{x^\alpha}
\]
which agrees with the definition of the generalised Poisson bracket
(\ref{generalised Poisson bracket on matrix-functions}).

\section{Proof of formula (\ref{subprincipal symbol of OI at time zero})}
\label{Proof of formula}

In this section we prove formula
(\ref{subprincipal symbol of OI at time zero}).
Our approach is as follows.

We write down explicitly the transport equations
(\ref{Leading transport equations equation 8}) at $t=0$,
i.e.
\begin{equation}
\label{Proof of formula equation 1}
\bigl[v^{(l)}\bigr]^*
\,
\left.
\bigl[\mathfrak{S}^{(j)}_{-1}f^{(j)}_1+\mathfrak{S}^{(j)}_0f^{(j)}_0\bigr]
\right|_{t=0}
=0,
\qquad l\ne j.
\end{equation}
We use the same local coordinates for $x$ and $y$ and we
assume all our phase functions to be linear, i.e.~we assume
that for each $j$ we have
(\ref{Leading transport equations equation 14}).
Using linear phase functions is justified for small $t$ because we
have
$(\xi^{(j)}_\eta)_\alpha{}^\beta(0;y,\eta)=\delta_\alpha{}^\beta$
and, hence, $\det\varphi^{(j)}_{x^\alpha\eta_\beta}(t,x;y,\eta)\ne0$
for small $t$. Writing down equations
(\ref{Proof of formula equation 1}) for linear phase functions
is much easier than for general phase functions
(\ref{algorithm equation 2}).

Using linear phase functions has the additional
advantage that the initial condition
(\ref{algorithm equation 32}) simplifies and reads now
$\sum_ju^{(j)}(0;y,\eta)=I$.
In view of
(\ref{decomposition of symbol of OI into homogeneous components}),
this implies, in particular, that
\begin{equation}
\label{Proof of formula equation 2}
\sum_j
u^{(j)}_{-1}(0)=0.
\end{equation}
Here and further on in this section we drop,
for the sake of brevity, the arguments
$(y,\eta)$ in $u^{(j)}_{-1}$.

Of course, the formula we are proving,
formula (\ref{subprincipal symbol of OI at time zero}),
does not depend
on our choice of phase functions. It is just easier to carry
out calculations for linear phase functions.

We will show that
(\ref{Proof of formula equation 1})
is a system of complex linear algebraic equations for the unknowns
$u^{(j)}_{-1}(0)$. The total number of equations
(\ref{Proof of formula equation 1}) is $m^2-m$. However, for each
$j$ and $l$
the LHS of (\ref{Proof of formula equation 1}) is a row of $m$
elements, so (\ref{Proof of formula equation 1}) is, effectively,
a system of $m(m^2-m)$ scalar equations.

Equation
(\ref{Proof of formula equation 2})
is a single matrix equation, so it is,
effectively,
a system of $m^2$ scalar equations.

Consequently, the system
(\ref{Proof of formula equation 1}),
(\ref{Proof of formula equation 2}) is, effectively,
a system of $m^3$ scalar equations.
This is exactly the number of unknown scalar elements
in the $m$ matrices $u^{(j)}_{-1}(0)$.

In the remainder of this section we write down explicitly
the LHS of (\ref{Proof of formula equation 1})
and solve the linear algebraic system
(\ref{Proof of formula equation 1}),
(\ref{Proof of formula equation 2})
for the unknowns
$u^{(j)}_{-1}(0)$.
This will allow us to prove formula
(\ref{subprincipal symbol of OI at time zero}).

Before starting explicit calculations we observe that
equations (\ref{Proof of formula equation 1}) can be equivalently rewritten as
\begin{equation}
\label{Proof of formula equation 3}
P^{(l)}
\,
\left.
\bigl[\mathfrak{S}^{(j)}_{-1}f^{(j)}_1+\mathfrak{S}^{(j)}_0f^{(j)}_0\bigr]
\right|_{t=0}
=0,
\qquad l\ne j,
\end{equation}
where $P^{(l)}:=[v^{(l)}(y,\eta)]\,[v^{(l)}(y,\eta)]^*$
is the orthogonal projection onto the eigenspace corresponding to
the (normalised) eigenvector $v^{(l)}(y,\eta)$ of the principal
symbol.
We will deal with
(\ref{Proof of formula equation 3})
rather than with
(\ref{Proof of formula equation 1}).
This is simply a matter of convenience.

\subsection{Part 1 of the proof of formula (\ref{subprincipal symbol of OI at time zero})}
\label{Part 1}

Our task in this subsection is to calculate
the LHS of (\ref{Proof of formula equation 3}).
In our calculations we use the explicit formula
(\ref{formula for principal symbol of oscillatory integral})
for the principal symbol $u^{(j)}_0(t;y,\eta)$
which was proved in Section~\ref{Leading transport equations}.

At $t=0$ formula (\ref{Leading transport equations equation 4}) reads
\[
\left.
\bigl[\mathfrak{S}^{(j)}_{-1}f^{(j)}_1\bigr]
\right|_{t=0}
=
i
\left.
\left[
\frac{\partial^2}{\partial x^\alpha\eta_\alpha}
\bigl(
A_1(x,\eta)
-h^{(j)}(y,\eta)
-(x-y)^\gamma h^{(j)}_{y^\gamma}(y,\eta)
\bigr)
P^{(j)}(y,\eta)
\right]
\right|_{x=y}
\]
which gives us
\begin{equation}
\label{Proof of formula equation 4}
\left.
\bigl[\mathfrak{S}^{(j)}_{-1}f^{(j)}_1\bigr]
\right|_{t=0}
=
i
\left[
(A_1-h^{(j)})_{y^\alpha\eta_\alpha}P^{(j)}
+
(A_1-h^{(j)})_{y^\alpha}P^{(j)}_{\eta_\alpha}
\right].
\end{equation}
In the latter formula we dropped, for the sake of brevity,
the arguments $(y,\eta)$.

At $t=0$ formula (\ref{Leading transport equations equation 5}) reads
\begin{multline}
\label{Proof of formula equation 5}
\left.
\bigl[\mathfrak{S}^{(j)}_0f^{(j)}_0\bigr]
\right|_{t=0}
=
-i\{v^{(j)},h^{(j)}\}[v^{(j)}]^*
+
\left(
A_0
-
q^{(j)}
+
\frac i2h^{(j)}_{y^\alpha\eta_\alpha}
\right)
P^{(j)}
\\
+[A_1-h^{(j)}]u^{(j)}_{-1}(0)\,,
\end{multline}
where $q^{(j)}$ is the function
(\ref{phase appearing in principal symbol})
and we dropped, for the sake of brevity,
the arguments $(y,\eta)$.
Note that in writing down
(\ref{Proof of formula equation 5})
we used the fact that
\[
\left.
\left[
(d_{\varphi^{(j)}})^{-1}
\partial_t
d_{\varphi^{(j)}}
\right]
\right|_{(t,x;y,\eta)=(0,y;y,\eta)}
=-\frac12
h^{(j)}_{y^\alpha\eta_\alpha}(y,\eta)\,,
\]
compare with formula
(\ref{Leading transport equations equation 20}).

Substituting formulae
(\ref{Proof of formula equation 4})
and
(\ref{Proof of formula equation 5})
into
(\ref{Proof of formula equation 3})
we get
\begin{equation}
\label{Proof of formula equation 6}
(h^{(l)}-h^{(j)})P^{(l)}u^{(j)}_{-1}(0)+P^{(l)}B^{(j)}_0=0,
\qquad l\ne j,
\end{equation}
where
\begin{equation}
\label{Part 1 result}
B^{(j)}_0=
\left(
A_0-q^{(j)}-\frac i2h^{(j)}_{y^\alpha\eta_\alpha}+i(A_1)_{y^\alpha\eta_\alpha}
\right)
P^{(j)}
-i
h^{(j)}_{\eta_\alpha}P^{(j)}_{y^\alpha}
+i(A_1)_{y^\alpha}P^{(j)}_{\eta_\alpha}.
\end{equation}
The subscript in $B^{(j)}_0$ indicates the degree of homogeneity in $\eta$.

\subsection{Part 2 of the proof of formula (\ref{subprincipal symbol of OI at time zero})}
\label{Part 2}

Our task in this subsection is to
solve the linear algebraic system
(\ref{Proof of formula equation 6}),
(\ref{Proof of formula equation 2})
for the unknowns
$u^{(j)}_{-1}(0)$.

It is easy to see that
the unique solution to the system
(\ref{Proof of formula equation 6}),
(\ref{Proof of formula equation 2})
is
\begin{equation}
\label{Part 2 result}
u^{(j)}_{-1}(0)
=\sum_{l\ne j}
\frac
{P^{(l)}B^{(j)}_0+P^{(j)}B^{(l)}_0}
{h^{(j)}-h^{(l)}}\,.
\end{equation}
Summation in (\ref{Part 2 result}) is carried out over all $l$
different from $j$.

\subsection{Part 3 of the proof of formula (\ref{subprincipal symbol of OI at time zero})}
\label{Part 3}

Our task in this subsection is to calculate $[U^{(j)}(0)]_\mathrm{sub}$.

We have
\begin{equation}
\label{subprincipal symbol of Uj0 equation 1}
[U^{(j)}(0)]_\mathrm{sub}
=u^{(j)}_{-1}(0)-\frac i2P^{(j)}_{y^\alpha\eta_\alpha}.
\end{equation}
Here the sign in front of $\frac i2$ is opposite to that in
(\ref{definition of subprincipal symbol})
because the way we write $U^{(j)}(0)$ is using the dual symbol.

Substituting
(\ref{Part 2 result})
and
(\ref{Part 1 result})
into (\ref{subprincipal symbol of Uj0 equation 1})
we get
\begin{multline}
\label{subprincipal symbol of Uj0 equation 2}
[U^{(j)}(0)]_\mathrm{sub}
=
-\frac i2P^{(j)}_{y^\alpha\eta_\alpha}
+\sum_{l\ne j}\frac1{h^{(j)}-h^{(l)}}
\\
\times
\bigl(
P^{(l)}
[
(A_0+i(A_1)_{y^\alpha\eta_\alpha})P^{(j)}
-ih^{(j)}_{\eta_\alpha}P^{(j)}_{y^\alpha}
+i(A_1)_{y^\alpha}P^{(j)}_{\eta_\alpha}
]
\\
\qquad\qquad+
P^{(j)}
[
(A_0+i(A_1)_{y^\alpha\eta_\alpha})P^{(l)}
-ih^{(l)}_{\eta_\alpha}P^{(l)}_{y^\alpha}
+i(A_1)_{y^\alpha}P^{(l)}_{\eta_\alpha}
]
\bigr)
\\
=
\sum_{l\ne j}
\frac
{
P^{(l)}A_\mathrm{sub}P^{(j)}
+
P^{(j)}A_\mathrm{sub}P^{(l)}
}
{
h^{(j)}-h^{(l)}
}
+\frac i2
\Bigl(
-P^{(j)}_{y^\alpha\eta_\alpha}
+
\sum_{l\ne j}
\frac
{
G_{jl}
}
{
h^{(j)}-h^{(l)}
}
\Bigr)\,,
\end{multline}
where
\begin{multline*}
G_{jl}:=
P^{(l)}
[
(A_1)_{y^\alpha\eta_\alpha}P^{(j)}
-2h^{(j)}_{\eta_\alpha}P^{(j)}_{y^\alpha}
+2(A_1)_{y^\alpha}P^{(j)}_{\eta_\alpha}
]
\\
+
P^{(j)}
[
(A_1)_{y^\alpha\eta_\alpha}P^{(l)}
-2h^{(l)}_{\eta_\alpha}P^{(l)}_{y^\alpha}
+2(A_1)_{y^\alpha}P^{(l)}_{\eta_\alpha}
]
\,.
\end{multline*}

We have
\begin{multline*}
G_{jl}
=
2P^{(l)}\{A_1,P^{(j)}\}
+
2P^{(j)}\{A_1,P^{(l)}\}
\\
+
P^{(l)}
[
(A_1-h^{(j)})_{y^\alpha\eta_\alpha}P^{(j)}
+2(A_1-h^{(j)})_{\eta_\alpha}P^{(j)}_{y^\alpha}
]
\\
+
P^{(j)}
[
(A_1-h^{(l)})_{y^\alpha\eta_\alpha}P^{(l)}
+2(A_1-h^{(l)})_{\eta_\alpha}P^{(l)}_{y^\alpha}
]
\\
=
2P^{(l)}\{A_1,P^{(j)}\}
+
2P^{(j)}\{A_1,P^{(l)}\}
-
P^{(l)}\{A_1-h^{(j)},P^{(j)}\}
-
P^{(j)}\{A_1-h^{(l)},P^{(l)}\}
\\
+
P^{(l)}
[
(A_1-h^{(j)})_{y^\alpha\eta_\alpha}P^{(j)}
+(A_1-h^{(j)})_{\eta_\alpha}P^{(j)}_{y^\alpha}
+(A_1-h^{(j)})_{y^\alpha}P^{(j)}_{\eta_\alpha}
]
\\
+
P^{(j)}
[
(A_1-h^{(l)})_{y^\alpha\eta_\alpha}P^{(l)}
+(A_1-h^{(l)})_{\eta_\alpha}P^{(l)}_{y^\alpha}
+(A_1-h^{(l)})_{y^\alpha}P^{(l)}_{\eta_\alpha}
]
\\
=
P^{(l)}\{A_1+h^{(j)},P^{(j)}\}
+
P^{(j)}\{A_1+h^{(l)},P^{(l)}\}
\\
-
P^{(l)}
(A_1-h^{(j)})
P^{(j)}_{y^\alpha\eta_\alpha}
-
P^{(j)}
(A_1-h^{(l)})
P^{(l)}_{y^\alpha\eta_\alpha}
\\
=
P^{(l)}\{A_1+h^{(j)},P^{(j)}\}
+
P^{(j)}\{A_1+h^{(l)},P^{(l)}\}
\\
-
P^{(l)}
(h^{(l)}-h^{(j)})
P^{(j)}_{y^\alpha\eta_\alpha}
-
P^{(j)}
(h^{(j)}-h^{(l)})
P^{(l)}_{y^\alpha\eta_\alpha}
\\
=
P^{(l)}\{A_1+h^{(j)},P^{(j)}\}
+
P^{(j)}\{A_1+h^{(l)},P^{(l)}\}
+(h^{(j)}-h^{(l)})
(
P^{(l)}
P^{(j)}_{y^\alpha\eta_\alpha}
-
P^{(j)}
P^{(l)}_{y^\alpha\eta_\alpha}
)\,,
\end{multline*}
so formula (\ref{subprincipal symbol of Uj0 equation 2}) can be rewritten as
\begin{multline}
\label{subprincipal symbol of Uj0 equation 3}
[U^{(j)}(0)]_\mathrm{sub}
=
\frac i2
\Bigl(
-P^{(j)}_{y^\alpha\eta_\alpha}
+
\sum_{l\ne j}
(
P^{(l)}
P^{(j)}_{y^\alpha\eta_\alpha}
-
P^{(j)}
P^{(l)}_{y^\alpha\eta_\alpha}
)
\Bigr)
\\
+
\frac12
\sum_{l\ne j}
\frac
{
P^{(l)}(2A_\mathrm{sub}P^{(j)}+i\{A_1+h^{(j)},P^{(j)}\})
+
P^{(j)}(2A_\mathrm{sub}P^{(l)}+i\{A_1+h^{(l)},P^{(l)}\})
}
{
h^{(j)}-h^{(l)}
}\,.
\end{multline}

But
\begin{multline*}
\sum_{l\ne j}
(
P^{(l)}
P^{(j)}_{y^\alpha\eta_\alpha}
-
P^{(j)}
P^{(l)}_{y^\alpha\eta_\alpha}
)
=
\Bigl(\,
\sum_{l\ne j}
P^{(l)}
\Bigr)
P^{(j)}_{y^\alpha\eta_\alpha}
-
P^{(j)}
\Bigl(\,
\sum_{l\ne j}
P^{(l)}
\Bigr)_{y^\alpha\eta_\alpha}
\\
=(I-P^{(j)})P^{(j)}_{y^\alpha\eta_\alpha}
-P^{(j)}(I-P^{(j)})_{y^\alpha\eta_\alpha}
=P^{(j)}_{y^\alpha\eta_\alpha},
\end{multline*}
so formula (\ref{subprincipal symbol of Uj0 equation 3}) can be simplified to read
\begin{multline}
\label{Part 3 result}
[U^{(j)}(0)]_\mathrm{sub}
\\
=
\frac12
\sum_{l\ne j}
\frac
{
P^{(l)}(2A_\mathrm{sub}P^{(j)}+i\{A_1+h^{(j)},P^{(j)}\})
+
P^{(j)}(2A_\mathrm{sub}P^{(l)}+i\{A_1+h^{(l)},P^{(l)}\})
}
{
h^{(j)}-h^{(l)}
}
\,.
\end{multline}

\subsection{Part 4 of the proof of formula (\ref{subprincipal symbol of OI at time zero})}
\label{Part 4}

Our task in this subsection is to calculate $\operatorname{tr}[U^{(j)}(0)]_\mathrm{sub}$.

Formula (\ref{Part 3 result}) implies
\begin{equation}
\label{trace of subprincipal symbol of Uj0 equation 1}
\operatorname{tr}[U^{(j)}(0)]_\mathrm{sub}
=
\frac i2\operatorname{tr}\sum_{l\ne j}
\frac
{
P^{(l)}\{A_1,P^{(j)}\}
+
P^{(j)}\{A_1,P^{(l)}\}
}
{
h^{(j)}-h^{(l)}
}
\,.
\end{equation}
Put $A_1=\sum_kh^{(k)}P^{(k)}$ and observe that
\begin{itemize}
\item
terms with the derivatives of $h$ vanish and
\item
the only $k$ which may give nonzero contributions are $k=j$ and $k=l$.
\end{itemize}
Thus, formula (\ref{trace of subprincipal symbol of Uj0 equation 1}) becomes
\begin{multline}
\label{trace of subprincipal symbol of Uj0 equation 2}
\operatorname{tr}[U^{(j)}(0)]_\mathrm{sub}
=
\frac i2\operatorname{tr}\sum_{l\ne j}
\frac1
{
h^{(j)}-h^{(l)}
}
\\
\times\bigl(
h^{(j)}
[
P^{(l)}\{P^{(j)},P^{(j)}\}
+
P^{(j)}\{P^{(j)},P^{(l)}\}
]
+
h^{(l)}
[
P^{(l)}\{P^{(l)},P^{(j)}\}
+
P^{(j)}\{P^{(l)},P^{(l)}\}
]
\bigr).
\end{multline}

We claim that
\begin{multline}
\label{Part 4 auxiliary equation 1}
\operatorname{tr}(P^{(l)}\{P^{(j)},P^{(j)}\})
=
\operatorname{tr}(P^{(j)}\{P^{(j)},P^{(l)}\})
\\
=
-\operatorname{tr}(P^{(l)}\{P^{(l)},P^{(j)}\})
=
-\operatorname{tr}(P^{(j)}\{P^{(l)},P^{(l)}\})
\\
=[v^{(l)}]^*\{v^{(j)},[v^{(j)}]^*\}v^{(l)}
\\
=([v^{(l)}]^*v^{(j)}_{y^\alpha})([v^{(j)}_{\eta_\alpha}]^*v^{(l)})
-([v^{(l)}]^*v^{(j)}_{\eta_\alpha})([v^{(j)}_{y^\alpha}]^*v^{(l)}).
\end{multline}
These facts are established by writing the orthogonal projections
in terms of the eigenvectors and using, if required,
the identities
\[
[v^{(l)}_{y^\alpha}]^*v^{(j)}+[v^{(l)}]^*v^{(j)}_{y^\alpha}=0,
\qquad
[v^{(l)}_{\eta_\alpha}]^*v^{(j)}+[v^{(l)}]^*v^{(j)}_{\eta_\alpha}=0,
\]
\[
[v^{(j)}_{y^\alpha}]^*v^{(l)}+[v^{(j)}]^*v^{(l)}_{y^\alpha}=0,
\qquad
[v^{(j)}_{\eta_\alpha}]^*v^{(l)}+[v^{(j)}]^*v^{(l)}_{\eta_\alpha}=0.
\]
In view of the identities (\ref{Part 4 auxiliary equation 1})
formula (\ref{trace of subprincipal symbol of Uj0 equation 2})
can be rewritten as
\begin{multline}
\label{Part 4 auxiliary equation 2}
\operatorname{tr}[U^{(j)}(0)]_\mathrm{sub}
=
i\operatorname{tr}
\sum_{l\ne j}
P^{(l)}\{P^{(j)},P^{(j)}\}
\\
=
i\operatorname{tr}
(\{P^{(j)},P^{(j)}\}-P^{(j)}\{P^{(j)},P^{(j)}\})
=
-i\operatorname{tr}
(P^{(j)}\{P^{(j)},P^{(j)}\}).
\end{multline}

It remains only to simplify the expression in the RHS of (\ref{Part 4 auxiliary equation 2}).
We have
\begin{multline}
\label{Part 4 auxiliary equation 3}
\operatorname{tr}
(P^{(j)}\{P^{(j)},P^{(j)}\})
=
\{[v^{(j)}]^*,v^{(j)}\}
\\
+[([v^{(j)}]^*v^{(j)}_{y^\alpha})([v^{(j)}]^*v^{(j)}_{\eta_\alpha})-([v^{(j)}]^*v^{(j)}_{\eta_\alpha})([v^{(j)}]^*v^{(j)}_{y^\alpha})]
\\
+[([v^{(j)}_{y^\alpha}]^*v^{(j)})([v^{(j)}_{\eta_\alpha}]^*v^{(j)})-([v^{(j)}_{\eta_\alpha}]^*v^{(j)})([v^{(j)}_{y^\alpha}]^*v^{(j)})]
\\
+[([v^{(j)}]^*v^{(j)}_{y^\alpha})([v^{(j)}_{\eta_\alpha}]^*v^{(j)})-([v^{(j)}]^*v^{(j)}_{\eta_\alpha})([v^{(j)}_{y^\alpha}]^*v^{(j)})]
\\
=
\{[v^{(j)}]^*,v^{(j)}\}
+[([v^{(j)}]^*v^{(j)}_{y^\alpha})([v^{(j)}_{\eta_\alpha}]^*v^{(j)})-([v^{(j)}]^*v^{(j)}_{\eta_\alpha})([v^{(j)}_{y^\alpha}]^*v^{(j)})]
\\
=
\{[v^{(j)}]^*,v^{(j)}\}
-[([v^{(j)}]^*v^{(j)}_{y^\alpha})([v^{(j)}]^*v^{(j)}_{\eta_\alpha})-([v^{(j)}]^*v^{(j)}_{\eta_\alpha})([v^{(j)}]^*v^{(j)}_{y^\alpha})]
\\
=
\{[v^{(j)}]^*,v^{(j)}\}.
\end{multline}
Formulae
(\ref{Part 4 auxiliary equation 2})
and
(\ref{Part 4 auxiliary equation 3})
imply formula (\ref{subprincipal symbol of OI at time zero}).

\section{$\mathrm{U}(1)$ connection}
\label{U(1) connection}

In the preceding Sections
\ref{Algorithm for the construction of the wave group}--\ref{Proof of formula}
we presented technical details
of the construction of the propagator. We saw that
the eigenvectors of the principal symbol, $v^{(j)}(x,\xi)$, play a major role
in this construction. As pointed out in Section~\ref{Main results},
each of these eigenvectors is
defined up to a $\mathrm{U}(1)$ gauge transformation
(\ref{gauge transformation of the eigenvector}),
(\ref{phase appearing in gauge transformation}).
In the end, the full symbols
(\ref{decomposition of symbol of OI into homogeneous components})
of our oscillatory integrals $U^{(j)}(t)$
do not depend on the choice of gauge for the eigenvectors $v^{(j)}(x,\xi)$.
However, the effect of the gauge transformations
(\ref{gauge transformation of the eigenvector}),
(\ref{phase appearing in gauge transformation})
is not as trivial as it may appear at first sight.
We will show in this section that the gauge transformations
(\ref{gauge transformation of the eigenvector}),
(\ref{phase appearing in gauge transformation})
show up, in the form of invariantly defined curvature, in the lower
order terms $u^{(j)}_{-1}(t;y,\eta)$ of the symbols of our oscillatory integrals $U^{(j)}(t)$.
More precisely, we will show that the RHS of
formula~(\ref{subprincipal symbol of OI at time zero})
is the scalar curvature of a connection associated with the gauge transformation
(\ref{gauge transformation of the eigenvector}),
(\ref{phase appearing in gauge transformation}).
Further on in this section, until the very last paragraph, the index $j$ enumerating eigenvalues and
eigenvectors of the principal symbol is assumed to be fixed.

Consider a smooth curve $\Gamma\subset T'M$ connecting points $(y,\eta)$ and $(x,\xi)$.
We write this curve in parametric form as $(z(t),\zeta(t))$, $t\in[0,1]$,
so that $(z(0),\zeta(0))=(y,\eta)$ and $(z(1),\zeta(1))=(x,\xi)$.
Put
\begin{equation}
\label{derivative of eigenvector is orthogonal to eigenvector auxiliary}
w(t):=e^{i\phi(t)}v^{(j)}(z(t),\zeta(t))\,,
\end{equation}
where
$\phi:[0,1]\to\mathbb{R}$
is an unknown function which is to be determined from the condition
\begin{equation}
\label{derivative of eigenvector is orthogonal to eigenvector}
iw^*\dot w=0
\end{equation}
with the dot indicating the derivative with respect to the parameter $t$.
Substituting (\ref{derivative of eigenvector is orthogonal to eigenvector auxiliary})
into
(\ref{derivative of eigenvector is orthogonal to eigenvector})
we get an ordinary differential equation for $\phi$ which
is easily solved, giving
\begin{multline}
\label{formula for phi(1)}
\phi(1)
=\phi(0)+\int_0^1(\dot z^\alpha(t)\,P_\alpha(z(t),\zeta(t))+\dot\zeta_\gamma(t)\,Q^\gamma(z(t),\zeta(t)))\,dt
\\
=\phi(0)+\int_\Gamma(P_\alpha dz^\alpha+Q^\gamma d\zeta_\gamma)\,,
\end{multline}
where
\begin{equation}
\label{formula for P and Q}
P_\alpha:=i[v^{(j)}]^*v^{(j)}_{z^\alpha},
\qquad
Q^\gamma:=i[v^{(j)}]^*v^{(j)}_{\zeta_\gamma}.
\end{equation}
Note that the $2n$-component real quantity $(P_\alpha,Q^\gamma)$
is a covector field (1-form) on $T'M$. This quantity already appeared
in Section~\ref{Main results} as formula (\ref{electromagnetic covector potential}).

Put $f(y,\eta):=e^{i\phi(0)}$, $f(x,\xi):=e^{i\phi(1)}$
and rewrite formula (\ref{formula for phi(1)}) as
\begin{equation}
\label{formula for a(1)}
f(x,\xi)
=f(y,\eta)\,e^{i\int_\Gamma(P_\alpha dz^\alpha+Q^\gamma d\zeta_\gamma)}.
\end{equation}
Let us identify the group $\mathrm{U}(1)$ with the unit circle in the complex
plane, i.e. with $f\in\mathbb{C}$, $\|f\|=1$.
We see that formulae (\ref{formula for a(1)}) and (\ref{formula for P and Q})
give us a rule for the parallel transport of elements
of the group $\mathrm{U}(1)$ along curves in $T'M$. This is the natural
$\mathrm{U}(1)$ connection generated by the normalised field of columns of
complex-valued scalars
\begin{equation}
\label{jth eigenvector of the principal symbol}
v^{(j)}(z,\zeta)=
\bigl(
\begin{matrix}v^{(j)}_1(z,\zeta)&\ldots&v^{(j)}_m(z,\zeta)\end{matrix}
\bigr)^T.
\end{equation}
Recall that the $\Gamma$ appearing in formula (\ref{formula for a(1)}) is a curve
connecting points $(y,\eta)$ and $(x,\xi)$, whereas
the $v^{(j)}(z,\zeta)$ appearing in formulae
(\ref{formula for P and Q}) and (\ref{jth eigenvector of the principal symbol})
enters our construction as
an eigenvector of the principal symbol of our $m\times m$ matrix pseudo\-differential
operator $A$.

In practice, dealing with a connection is not as convenient as dealing with
the covariant derivative $\nabla$. The covariant derivative
corresponding to the connection (\ref{formula for a(1)}) is determined as follows.
Let us view the $(x,\xi)$ appearing in formula (\ref{formula for a(1)})
as a variable which takes values close to $(y,\eta)$,
and suppose that the curve $\Gamma$ is a short straight (in local coordinates)
line segment connecting the point $(y,\eta)$ with the point $(x,\xi)$.
We want the covariant derivative of our function
$f(x,\xi)$, evaluated at $(y,\eta)$, to be zero.
Examination of formula (\ref{formula for a(1)}) shows that
the unique covariant derivative satisfying this condition is
\begin{equation}
\label{formula for U(1) covariant derivative}
\nabla_\alpha:=\partial/\partial x^\alpha-iP_\alpha(x,\xi),
\qquad
\nabla^\gamma:=\partial/\partial\xi_\gamma-iQ^\gamma(x,\xi).
\end{equation}

We define the curvature of our $\mathrm{U}(1)$ connection as
\begin{equation}
\label{definition of U(1) curvature}
R:=
-i
\begin{pmatrix}
\nabla_\alpha\nabla_\beta-\nabla_\beta\nabla_\alpha&
\nabla_\alpha\nabla^\delta-\nabla^\delta\nabla_\alpha
\\
\nabla^\gamma\nabla_\beta-\nabla_\beta\nabla^\gamma&
\nabla^\gamma\nabla^\delta-\nabla^\delta\nabla^\gamma
\end{pmatrix}.
\end{equation}
It may seem that the entries of the $(2n)\times(2n)$ matrix (\ref{definition of U(1) curvature})
are differential operators. They are, in fact, operators of multiplication
by ``scalar functions''. Namely, the more explicit form of (\ref{definition of U(1) curvature}) is
\begin{equation}
\label{explicit formula for U(1) curvature}
R=
\begin{pmatrix}
\frac{\partial P_\alpha}{\partial x^\beta}-\frac{\partial P_\beta}{\partial x^\alpha}&
\frac{\partial P_\alpha}{\partial\xi_\delta}-\frac{\partial Q^\delta}{\partial x^\alpha}
\\
\frac{\partial Q^\gamma}{\partial x^\beta}-\frac{\partial P_\beta}{\partial\xi_\gamma}&
\frac{\partial Q^\gamma}{\partial\xi_\delta}-\frac{\partial Q^\delta}{\partial\xi_\gamma}
\end{pmatrix}.
\end{equation}
The $(2n)\times(2n)$\,-\,component real quantity (\ref{explicit formula for U(1) curvature})
is a rank 2 covariant antisymmetric tensor (2-form) on $T'M$.
It is an analogue of the electromagnetic tensor.

Substituting (\ref{formula for P and Q}) into
(\ref{explicit formula for U(1) curvature})
we get an expression for curvature in terms of the eigenvector
of the principal symbol
\begin{equation}
\label{more explicit formula for U(1) curvature}
R=i
\begin{pmatrix}
[v^{(j)}_{x^\beta}]^*v^{(j)}_{x^\alpha}-[v^{(j)}_{x^\alpha}]^*v^{(j)}_{x^\beta}&
[v^{(j)}_{\xi_\delta}]^*v^{(j)}_{x^\alpha}-[v^{(j)}_{x^\alpha}]^*v^{(j)}_{\xi_\delta}
\\
[v^{(j)}_{x^\beta}]^*v^{(j)}_{\xi_\gamma}-[v^{(j)}_{\xi_\gamma}]^*v^{(j)}_{x^\beta}&
[v^{(j)}_{\xi_\delta}]^*v^{(j)}_{\xi_\gamma}-[v^{(j)}_{\xi_\gamma}]^*v^{(j)}_{\xi_\delta}
\end{pmatrix}.
\end{equation}
Examination of formula (\ref{more explicit formula for U(1) curvature}) shows that,
as expected, curvature is invariant under the gauge transformation
(\ref{gauge transformation of the eigenvector}),
(\ref{phase appearing in gauge transformation}).

It is natural to take the trace of the upper right block
in (\ref{definition of U(1) curvature}) which,
in the notation (\ref{Poisson bracket on matrix-functions}), gives us
\begin{equation}
\label{scalar curvature of U(1) connection}
-i(\nabla_\alpha\nabla^\alpha-\nabla^\alpha\nabla_\alpha)
=-i\{[v^{(j)}]^*,v^{(j)}\}.
\end{equation}
Thus, we have shown that the RHS of
formula~(\ref{subprincipal symbol of OI at time zero})
is the scalar curvature of our $\mathrm{U}(1)$ connection.

\

We end this section by proving, as promised in Section~\ref{Main results},
formula (\ref{sum of curvatures is zero}) without referring to microlocal analysis.
In the following arguments we use our standard notation for the orthogonal
projections onto the eigenspaces of the principal symbol,
i.e.~we write $P^{(k)}:=v^{(k)}[v^{(k)}]^*$.
We have $\operatorname{tr}\{P^{(j)},P^{(j)}\}=0$
and $\sum_lP^{(l)}=I$
which implies
\begin{multline}
\label{sum of curvatures is zero proof equation 1}
0=\sum_{l,j}\operatorname{tr}(P^{(l)}\{P^{(j)},P^{(j)}\})
\\
=\sum_j\operatorname{tr}(P^{(j)}\{P^{(j)},P^{(j)}\})
+\sum_{l,j:\ l\ne j}\operatorname{tr}(P^{(l)}\{P^{(j)},P^{(j)}\}).
\end{multline}
But, according to formula (\ref{Part 4 auxiliary equation 1}),
for $l\ne j$ we have
\[
\operatorname{tr}(P^{(l)}\{P^{(j)},P^{(j)}\})
=-\operatorname{tr}(P^{(j)}\{P^{(l)},P^{(l)}\}),
\]
so formula
(\ref{sum of curvatures is zero proof equation 1}) can be rewritten as
$\sum_j\operatorname{tr}(P^{(j)}\{P^{(j)},P^{(j)}\})=0$.
It remains only to note that,
according to formula (\ref{Part 4 auxiliary equation 3}),
$\operatorname{tr}(P^{(j)}\{P^{(j)},P^{(j)}\})=\{[v^{(j)}]^*,v^{(j)}\}$.

\section{Singularity of the propagator at $t=0$}
\label{Singularity of the wave group at time zero}

Following the notation of \cite{mybook}, we denote by
\[
\mathcal{F}_{\lambda\to t}[f(\lambda)]=\hat f(t)=\int e^{-it\lambda}f(\lambda)\,d\lambda
\]
the one-dimensional Fourier transform and by
\[
\mathcal{F}^{-1}_{t\to\lambda}[\hat f(t)]=f(\lambda)=(2\pi)^{-1}\int e^{it\lambda}\hat f(t)\,dt
\]
its inverse.

Suppose that we have a Hamiltonian trajectory
$(x^{(j)}(t;y,\eta),\xi^{(j)}(t;y,\eta))$
and a real number $T>0$ such that
$x^{(j)}(T;y,\eta)=y$. We will say in this case
that we have a loop of length $T$ originating
from the point $y\in M$.

\begin{remark}
\label{remark on reversibility}
There is no need to consider loops of negative length $T$ because,
given a $T>0$, we have
$x^{(j)}(T;y,\eta^+)=y$
for some $\eta^+\in T'_yM$ if and only if we have
$x^{(j)}(-T;y,\eta^-)=y$
for some $\eta^-\in T'_yM$. Indeed,
it suffices to relate the
$\eta^\pm$ in accordance with
$\eta^\mp=\xi^{(j)}(\pm T;y,\eta^\pm)$.
\end{remark}

Denote by $\mathcal{T}^{(j)}\subset\mathbb{R}$ the set of lengths $T>0$
of all possible loops generated by the Hamiltonian $h^{(j)}$.
Here ``all possible'' refers to all possible starting points
$(y,\eta)\in T'M$ of Hamiltonian trajectories.
It is easy to see that $0\not\in\overline{\mathcal{T}^{(j)}}$.
We put
\[
\mathbf{T}^{(j)}:=
\begin{cases}
\inf\mathcal{T}^{(j)}\quad&\text{if}\quad\mathcal{T}^{(j)}\ne\emptyset,
\\
+\infty\quad&\text{if}\quad\mathcal{T}^{(j)}=\emptyset.
\end{cases}
\]

In the Riemannian case (i.e.~the case when the Hamiltonian
is a square root of a quadratic polynomial in $\xi$) it is known \cite{sabourau,rotman}
that there is a loop originating from every point of the
manifold $M$ and, moreover, there is an explicit estimate from above for
the number $\mathbf{T}^{(j)}$.
We are not aware of similar results for general Hamiltonians.

We also define
$\mathbf{T}:=\min\limits_{j=1,\ldots,m^+}\mathbf{T}^{(j)}$.

\begin{remark}
\label{remark on negative Hamiltonians}
Note that negative eigenvalues of the principal symbol,
i.e.~Hamiltonians $h^{(j)}(x,\xi)$ with negative index
$j=-1,\ldots,-m^-$,
do not affect the asymptotic formulae
we are about to derive. This is because we are dealing
with the case $\lambda\to+\infty$ rather than $\lambda\to-\infty$.
\end{remark}

Denote by
\begin{equation}
\label{definition of integral kernel of wave group}
u(t,x,y):=
\sum_k e^{-it\lambda_k}v_k(x)[v_k(y)]^*
\end{equation}
the integral kernel of the propagator (\ref{definition of wave group}).
The quantity (\ref{definition of integral kernel of wave group})
can be understood as a distribution in the variable
$t\in\mathbb{R}$ depending on the parameters $x,y\in M$.

The main result of this section is the following
\begin{lemma}
\label{Singularity of the wave group at time zero lemma}
Let $\hat\rho:\mathbb{R}\to\mathbb{C}$ be an infinitely smooth function such that
\begin{equation}
\label{condition on hat rho 1}
\operatorname{supp}\hat\rho\subset(-\mathbf{T},\mathbf{T}),
\end{equation}
\begin{equation}
\label{condition on hat rho 2}
\hat\rho(0)=1,
\end{equation}
\begin{equation}
\label{condition on hat rho 3}
\hat\rho'(0)=0.
\end{equation}
Then, uniformly over $y\in M$, we have
\begin{equation}
\label{Singularity of the wave group at time zero lemma formula}
\mathcal{F}^{-1}_{t\to\lambda}[\hat\rho(t)\operatorname{tr}u(t,y,y)]=
n\,a(y)\,\lambda^{n-1}+(n-1)\,b(y)\,\lambda^{n-2}+O(\lambda^{n-3})
\end{equation}
as $\lambda\to+\infty$.
The densities $a(y)$ and $b(y)$ appearing in the RHS of formula
(\ref{Singularity of the wave group at time zero lemma formula})
are defined in accordance with formulae
(\ref{formula for a(x)}) and (\ref{formula for b(x)}).
\end{lemma}

\emph{Proof\ }
Denote by $(S^*_yM)^{(j)}$ the $(n-1)$-dimensional unit cosphere in the cotangent fibre
defined by the equation $h^{(j)}(y,\eta)=1$
and denote by $d(S^*_yM)^{(j)}$ the
surface area element on $(S^*_yM)^{(j)}$
defined by the condition
$d\eta=d(S^*_yM)^{(j)}\,dh^{(j)}$.
The latter means that we use spherical coordinates in the cotangent fibre
with the Hamiltonian $h^{(j)}$
playing the role of the radial coordinate, see subsection 1.1.10 of \cite{mybook} for details.
In particular, as explained in subsection 1.1.10 of \cite{mybook},
our surface area element $d(S^*_yM)^{(j)}$ is expressed via the Euclidean surface area element as
\[
d(S^*_yM)^{(j)}=
\biggl(\,\sum_{\alpha=1}^n\bigl(h^{(j)}_{\eta_\alpha}(y,\eta)\bigr)^2\biggr)^{-1/2}
\times\,
\text{Euclidean surface area element}
\,.
\]
Denote also
$\,{d{\hskip-1pt\bar{}}\hskip1pt}(S^*_yM)^{(j)}:=(2\pi)^{-n}\,d(S^*_yM)^{(j)}\,$.

According to Corollary 4.1.5 from \cite{mybook} we have
uniformly over $y\in M$
\begin{multline}
\label{Singularity of the wave group at time zero lemma equation 1}
\mathcal{F}^{-1}_{t\to\lambda}[\hat\rho(t)\operatorname{tr}u(t,y,y)]=
\\
\sum_{j=1}^{m^+}
\left(c^{(j)}(y)\,\lambda^{n-1}+d^{(j)}(y)\,\lambda^{n-2}+e^{(j)}(y)\,\lambda^{n-2}\right)
+O(\lambda^{n-3})\,,
\end{multline}
where
\begin{equation}
\label{Singularity of the wave group at time zero lemma equation 2}
c^{(j)}(y)=\int\limits_{(S^*_yM)^{(j)}}
\operatorname{tr}u^{(j)}_0(0;y,\eta)
\,{d{\hskip-1pt\bar{}}\hskip1pt}(S^*_yM)^{(j)}\,,
\end{equation}
\begin{multline}
\label{Singularity of the wave group at time zero lemma equation 3}
d^{(j)}(y)=
\\
(n-1)\int\limits_{(S^*_yM)^{(j)}}
\operatorname{tr}
\left(
-\,i\,\dot u^{(j)}_0(0;y,\eta)
+\frac i2\bigl\{u^{(j)}_0\bigr|_{t=0}\,,h^{(j)}\bigr\}(y,\eta)
\right)
{d{\hskip-1pt\bar{}}\hskip1pt}(S^*_yM)^{(j)}\,,
\end{multline}
\begin{equation}
\label{Singularity of the wave group at time zero lemma equation 4}
e^{(j)}(y)=\int\limits_{(S^*_yM)^{(j)}}
\operatorname{tr}[U^{(j)}(0)]_\mathrm{sub}(y,\eta)
\,{d{\hskip-1pt\bar{}}\hskip1pt}(S^*_yM)^{(j)}\,.
\end{equation}
Here $u^{(j)}_0(t;y,\eta)$ is the principal symbol of the oscillatory integral
(\ref{algorithm equation 1}) and $\dot u^{(j)}_0(t;y,\eta)$ is its time derivative.
Note that in writing the term with the Poisson bracket in
(\ref{Singularity of the wave group at time zero lemma equation 3})
we took account of the fact that Poisson brackets in \cite{mybook}
and in the current paper have opposite signs.

Observe that the integrands in formulae
(\ref{Singularity of the wave group at time zero lemma equation 2})
and
(\ref{Singularity of the wave group at time zero lemma equation 3})
are positively homogeneous in $\eta$ of degree 0,
whereas the integrand in formula
(\ref{Singularity of the wave group at time zero lemma equation 4})
is positively homogeneous in $\eta$ of degree $-1$.
In order to have the same degree of homogeneity,  we rewrite
formula
(\ref{Singularity of the wave group at time zero lemma equation 4})
in equivalent form
\begin{equation}
\label{Singularity of the wave group at time zero lemma equation 5}
e^{(j)}(y)=\int\limits_{(S^*_yM)^{(j)}}
\bigl(
h^{(j)}\operatorname{tr}[U^{(j)}(0)]_\mathrm{sub}
\bigr)
(y,\eta)
\,{d{\hskip-1pt\bar{}}\hskip1pt}(S^*_yM)^{(j)}\,.
\end{equation}

Switching from surface integrals to volume integrals with the help of formula (1.1.15) from \cite{mybook},
we rewrite formulae
(\ref{Singularity of the wave group at time zero lemma equation 2}),
(\ref{Singularity of the wave group at time zero lemma equation 3})
and
(\ref{Singularity of the wave group at time zero lemma equation 5})
as
\begin{equation}
\label{Singularity of the wave group at time zero lemma equation 6}
c^{(j)}(y)=n\int\limits_{h^{(j)}(y,\eta)<1}
\operatorname{tr}u^{(j)}_0(0;y,\eta)
\,{d{\hskip-1pt\bar{}}\hskip1pt}\eta\,,
\end{equation}
\begin{multline}
\label{Singularity of the wave group at time zero lemma equation 7}
d^{(j)}(y)=n(n-1)\times
\\
\int\limits_{h^{(j)}(y,\eta)<1}
\operatorname{tr}
\left(
-\,i\,\dot u^{(j)}_0(0;y,\eta)
+\frac i2\bigl\{u^{(j)}_0\bigr|_{t=0}\,,h^{(j)}\bigr\}(y,\eta)
\right)
{d{\hskip-1pt\bar{}}\hskip1pt}\eta\,,
\end{multline}
\begin{equation}
\label{Singularity of the wave group at time zero lemma equation 8}
e^{(j)}(y)=n\int\limits_{h^{(j)}(y,\eta)<1}
\bigl(
h^{(j)}\operatorname{tr}[U^{(j)}(0)]_\mathrm{sub}
\bigr)
(y,\eta)
\,{d{\hskip-1pt\bar{}}\hskip1pt}\eta\,.
\end{equation}

Substituting formulae
(\ref{formula for principal symbol of oscillatory integral})
and
(\ref{phase appearing in principal symbol})
into formulae
(\ref{Singularity of the wave group at time zero lemma equation 6})
and
(\ref{Singularity of the wave group at time zero lemma equation 7})
we get
\begin{equation}
\label{Singularity of the wave group at time zero lemma equation 9}
c^{(j)}(y)=n\int\limits_{h^{(j)}(y,\eta)<1}
{d{\hskip-1pt\bar{}}\hskip1pt}\eta\,,
\end{equation}
\begin{multline}
\label{Singularity of the wave group at time zero lemma equation 10}
d^{(j)}(y)=-n(n-1)\times
\\
\int\limits_{h^{(j)}(y,\eta)<1}
\left(
[v^{(j)}]^*A_\mathrm{sub}v^{(j)}
-\frac i2
\{
[v^{(j)}]^*,A_1-h^{(j)},v^{(j)}
\}
\right)(y,\eta)
\,{d{\hskip-1pt\bar{}}\hskip1pt}\eta\,.
\end{multline}
Substituting formula
(\ref{subprincipal symbol of OI at time zero})
into formula
(\ref{Singularity of the wave group at time zero lemma equation 8})
we get
\begin{equation}
\label{Singularity of the wave group at time zero lemma equation 11}
e^{(j)}(y)=-n\,i\int\limits_{h^{(j)}(y,\eta)<1}
\bigl(
h^{(j)}\{[v^{(j)}]^*,v^{(j)}\}
\bigr)
(y,\eta)
\,{d{\hskip-1pt\bar{}}\hskip1pt}\eta\,.
\end{equation}

Substituting formulae
(\ref{Singularity of the wave group at time zero lemma equation 9})--(\ref{Singularity of the wave group at time zero lemma equation 11})
into formula
(\ref{Singularity of the wave group at time zero lemma equation 1})
we arrive
at (\ref{Singularity of the wave group at time zero lemma formula}).~$\square$

\begin{remark}
The proof of Lemma~\ref{Singularity of the wave group at time zero lemma}
given above was based on the use of Corollary 4.1.5 from \cite{mybook}.
In the actual statement of Corollary 4.1.5 in \cite{mybook}
uniformity in $y\in M$ was not mentioned because the authors were dealing with
a manifold with a boundary. Uniformity reappeared in the subsequent
Theorem 4.2.1 which involved pseudodifferential cut-offs
separating the point $\,y\,$ from the boundary.
\end{remark}

\section{Mollified spectral asymptotics}
\label{Mollified spectral asymptotics}

\begin{theorem}
\label{theorem spectral function mollified}
Let $\rho:\mathbb{R}\to\mathbb{C}$ be a function from Schwartz space $\mathcal{S}(\mathbb{R})$
whose Fourier transform $\hat\rho$ satisfies conditions
(\ref{condition on hat rho 1})--(\ref{condition on hat rho 3}).
Then, uniformly over $x\in M$, we have
\begin{equation}
\label{theorem spectral function mollified formula}
\int e(\lambda-\mu,x,x)\,\rho(\mu)\,d\mu=
a(x)\,\lambda^n+b(x)\,\lambda^{n-1}+
\begin{cases}
O(\lambda^{n-2})\quad&\text{if}\quad{n\ge3},
\\
O(\ln\lambda)\quad&\text{if}\quad{n=2},
\end{cases}
\end{equation}
as $\lambda\to+\infty$.
The densities $a(x)$ and $b(x)$ appearing in the RHS of formula
(\ref{theorem spectral function mollified formula})
are defined in accordance with formulae
(\ref{formula for a(x)}) and (\ref{formula for b(x)}).
\end{theorem}

\emph{Proof\ }
Our spectral function $e(\lambda,x,x)$ was initially defined only for $\lambda>0$,
see formula (\ref{definition of spectral function}). We extend the definition
to the whole real line by setting
\[
e(\lambda,x,x):=0\quad\text{for}\quad\lambda\le0.
\]

Denote by $e'(\lambda,x,x)$ the derivative, with respect to the spectral
parameter, of the spectral function. Here ``derivative'' is understood in the
sense of distributions. The explicit formula for $e'(\lambda,x,x)$ is
\begin{equation}
\label{theorem spectral function mollified equation 2}
e'(\lambda,x,x):=\sum_{k=1}^{+\infty}\|v_k(x)\|^2\,\delta(\lambda-\lambda_k).
\end{equation}

Formula (\ref{theorem spectral function mollified equation 2}) gives us
\begin{equation}
\label{theorem spectral function mollified equation 3}
\int e'(\lambda-\mu,x,x)\,\rho(\mu)\,d\mu=
\sum_{k=1}^{+\infty}\|v_k(x)\|^2\,\rho(\lambda-\lambda_k).
\end{equation}
Formula (\ref{theorem spectral function mollified equation 3})
implies, in particular, that, uniformly over $x\in M$, we have
\begin{equation}
\label{theorem spectral function mollified equation 4}
\int e'(\lambda-\mu,x,x)\,\rho(\mu)\,d\mu=O(|\lambda|^{-\infty})
\quad\text{as}\quad\lambda\to-\infty\,,
\end{equation}
where $O(|\lambda|^{-\infty})$ is shorthand for ``tends to zero faster
than any given inverse power of $|\lambda|$''.

Formula (\ref{theorem spectral function mollified equation 3})
can also be rewritten as
\begin{equation}
\label{theorem spectral function mollified equation 5}
\int e'(\lambda-\mu,x,x)\,\rho(\mu)\,d\mu=
\mathcal{F}^{-1}_{t\to\lambda}[\hat\rho(t)\operatorname{tr}u(t,x,x)]
-\sum_{k\le0}\|v_k(x)\|^2\,\rho(\lambda-\lambda_k)\,,
\end{equation}
where the distribution $u(t,x,y)$ is defined in accordance with
formula (\ref{definition of integral kernel of wave group}).
Clearly, we have
\begin{equation}
\label{theorem spectral function mollified equation 6}
\sum_{k\le0}\|v_k(x)\|^2\,\rho(\lambda-\lambda_k)=O(\lambda^{-\infty})
\quad\text{as}\quad\lambda\to+\infty\,.
\end{equation}
Formulae
(\ref{theorem spectral function mollified equation 5}),
(\ref{theorem spectral function mollified equation 6})
and Lemma~\ref{Singularity of the wave group at time zero lemma}
imply that, uniformly over $x\in M$, we have
\begin{multline}
\label{theorem spectral function mollified equation 7}
\int e'(\lambda-\mu,x,x)\,\rho(\mu)\,d\mu=
\\
n\,a(x)\,\lambda^{n-1}+(n-1)\,b(x)\,\lambda^{n-2}+O(\lambda^{n-3})
\quad\text{as}\quad\lambda\to+\infty\,.
\end{multline}

It remains to note that
\begin{equation}
\label{theorem spectral function mollified equation 8}
\frac d{d\lambda}\int e(\lambda-\mu,x,x)\,\rho(\mu)\,d\mu
=
\int e'(\lambda-\mu,x,x)\,\rho(\mu)\,d\mu\,.
\end{equation}
Formulae
(\ref{theorem spectral function mollified equation 8}),
(\ref{theorem spectral function mollified equation 4})
and
(\ref{theorem spectral function mollified equation 7})
imply
(\ref{theorem spectral function mollified formula}).~$\square$

\begin{theorem}
\label{theorem counting function mollified}
Let $\rho:\mathbb{R}\to\mathbb{C}$ be a function from Schwartz space $\mathcal{S}(\mathbb{R})$
whose Fourier transform $\hat\rho$ satisfies conditions
(\ref{condition on hat rho 1})--(\ref{condition on hat rho 3}).
Then we have
\begin{equation}
\label{theorem counting function mollified formula}
\int N(\lambda-\mu)\,\rho(\mu)\,d\mu=
a\,\lambda^n+b\,\lambda^{n-1}+
\begin{cases}
O(\lambda^{n-2})\quad&\text{if}\quad{n\ge3},
\\
O(\ln\lambda)\quad&\text{if}\quad{n=2},
\end{cases}
\end{equation}
as $\lambda\to+\infty$.
The constants $a$ and $b$ appearing in the RHS of formula
(\ref{theorem counting function mollified formula})
are defined in accordance with formulae
(\ref{a via a(x)}),
(\ref{formula for a(x)}),
(\ref{b via b(x)})
and
(\ref{formula for b(x)}).
\end{theorem}

\emph{Proof\ }
Formula
(\ref{theorem counting function mollified formula})
follows from formula
(\ref{theorem spectral function mollified formula})
by integration over $M$,
see also formula (\ref{definition of counting function}).~$\square$

\

In stating Theorems
\ref{theorem spectral function mollified}
and
\ref{theorem counting function mollified}
we assumed the mollifier $\rho$ to be complex-valued.
This was done for the sake of generality but may seem
unnatural when mollifying real-valued functions
$e(\lambda,x,x)$ and $N(\lambda)$. One can make our
construction look more natural by dealing only with
real-valued mollifiers $\rho$. Note that if the function $\rho$
is real-valued and even then its Fourier transform
$\hat\rho$ is also real-valued and even and, moreover,
condition (\ref{condition on hat rho 3}) is automatically satisfied.

\section{Unmollified spectral asymptotics}
\label{Unmollified spectral asymptotics}

In this section we derive asymptotic formulae for
the spectral function $e(\lambda,x,x)$ and the
counting function $N(\lambda)$ without mollification.
The section is split into two subsections: in the first
we derive one-term asymptotic formulae and
in the second --- two-term asymptotic formulae.

\subsection{One-term spectral asymptotics}
\label{One-term spectral asymptotics}

\begin{theorem}
\label{theorem spectral function unmollified one term}
We have, uniformly over $x\in M$,
\begin{equation}
\label{theorem spectral function unmollified one term formula}
e(\lambda,x,x)=a(x)\,\lambda^n+O(\lambda^{n-1})
\end{equation}
as $\lambda\to+\infty$.
\end{theorem}

\emph{Proof\ }
The result in question is an immediate consequence of
formulae
(\ref{theorem spectral function mollified equation 8}),
(\ref{theorem spectral function mollified equation 7})
and
Theorem~\ref{theorem spectral function mollified}
from the current paper
and Corollary~B.2.2 from \cite{mybook}.~$\square$

\begin{theorem}
\label{theorem counting function unmollified one term}
We have
\begin{equation}
\label{theorem counting function unmollified one term formula}
N(\lambda)=a\lambda^n+O(\lambda^{n-1})
\end{equation}
as $\lambda\to+\infty$.
\end{theorem}

\emph{Proof\ }
Formula
(\ref{theorem counting function unmollified one term formula})
follows from formula
(\ref{theorem spectral function unmollified one term formula})
by integration over $M$,
see also formula (\ref{definition of counting function}).~$\square$

\subsection{Two-term spectral asymptotics}
\label{Two-term spectral asymptotics}

Up till now, in Section~\ref{Mollified spectral asymptotics}
and subsection~\ref{One-term spectral asymptotics},
our logic was to derive asymptotic formulae for the spectral
function $e(\lambda,x,x)$ first and then obtain corresponding
asymptotic formulae for the counting function $N(\lambda)$
by integration over $M$. Such an approach will not work
for two-term asymptotics because
the geometric conditions required for the existence of
two-term asymptotics of $e(\lambda,x,x)$ and $N(\lambda)$
will be different:
for $e(\lambda,x,x)$ the appropriate geometric conditions
will be formulated in terms of \emph{loops},
whereas
for $N(\lambda)$ the appropriate geometric conditions
will be formulated in terms of \emph{periodic trajectories}.

Hence, in this subsection we deal with
the spectral function $e(\lambda,x,x)$
and the counting function $N(\lambda)$ separately.

In what follows the point $y\in M$ is assumed to be fixed.

Denote by $\Pi_y^{(j)}$ the set of normalised ($h^{(j)}(y,\eta)=1$)
covectors $\eta$ which serve as starting points for loops generated by the Hamiltonian
$h^{(j)}$. Here ``starting point'' refers to the starting point
of a Hamiltonian trajectory
$(x^{(j)}(t;y,\eta),\xi^{(j)}(t;y,\eta))$
moving forward in time ($t>0$),
see also Remark~\ref{remark on reversibility}.

The reason we are not interested in large negative $t$ is that the refined
Fourier Tauberian theorem we will be applying,
Theorem~B.5.1 from \cite{mybook},
does not require information regarding large negative $t$.
And the underlying reason for the latter is the fact that the function
we are studying, $e(\lambda,x,x)$ (and, later, $N(\lambda)$), is real-valued.
The real-valuedness of the function $e(\lambda,x,x)$ implies that its
Fourier transform, $\hat e(t,x,x)$, possesses the symmetry
$\hat e(-t,x,x)=\overline{\hat e(t,x,x)}$.

The set $\Pi_y^{(j)}$ is a subset of the $(n-1)$-dimensional unit cosphere
$(S^*_yM)^{(j)}$ and the latter is equipped with
a natural Lebesgue measure, see proof of
Lemma~\ref{Singularity of the wave group at time zero lemma}.
It is known, see Lemma 1.8.2 in \cite{mybook}, that the
set $\Pi_y^{(j)}$ is measurable.

\begin{definition}
\label{definition of nonfocal point 1}
A point $y\in M$ is said to be \emph{nonfocal} if for each
$j=1,\ldots,m^+$ the set $\Pi_y^{(j)}$ has measure zero.
\end{definition}

With regards to the range of the index $j$ in
Definition~\ref{definition of nonfocal point 1},
as well as in sub\-sequent
Definitions~\ref{definition of nonfocal point 2}--\ref{definition of nonperiodicity condition 2},
see Remark~\ref{remark on negative Hamiltonians}.

We call a loop of length $T>0$ \emph{absolutely focused} if
the function
\[
|x^{(j)}(T;y,\eta)-y|^2
\]
has an infinite order zero in the variable $\eta$, and we denote
by $(\Pi_y^a)^{(j)}$ the set of normalised ($h^{(j)}(y,\eta)=1$)
covectors $\eta$ which serve as starting points for absolutely focused loops
generated by the Hamiltonian $h^{(j)}$.
It is known, see Lemma 1.8.3 in \cite{mybook}, that the
set $(\Pi_y^a)^{(j)}$ is measurable and,
moreover, the set $\Pi_y^{(j)}\setminus(\Pi_y^a)^{(j)}$ has measure zero.
This allows us to reformulate Definition~\ref{definition of nonfocal point 1}
as follows.

\begin{definition}
\label{definition of nonfocal point 2}
A point $y\in M$ is said to be \emph{nonfocal} if for each
$j=1,\ldots,m^+$ the set $(\Pi_y^a)^{(j)}$ has measure zero.
\end{definition}

In practical applications it is easier to work with
Definition~\ref{definition of nonfocal point 2}
because the set $(\Pi_y^a)^{(j)}$ is usually much
thinner than the set $\Pi_y^{(j)}$.

In order to derive a two-term asymptotic formula for the
spectral function $e(\lambda,x,x)$ we need the following
lemma (compare with Lemma~\ref{Singularity of the wave group at time zero lemma}).

\begin{lemma}
\label{Singularity of the wave group at time nonzero lemma pointwise}
Suppose that the point $y\in M$ is nonfocal.
Then for any complex-valued function $\hat\gamma\in C_0^\infty(\mathbb{R})$
with $\operatorname{supp}\hat\gamma\subset(0,+\infty)$ we have
\begin{equation}
\label{Singularity of the wave group at time nonzero lemma pointwise formula}
\mathcal{F}^{-1}_{t\to\lambda}[\hat\gamma(t)\operatorname{tr}u(t,y,y)]=
o(\lambda^{n-1})
\end{equation}
as $\lambda\to+\infty$.
\end{lemma}

\emph{Proof\ }
The result in question is a special case of
Theorem~4.4.9 from \cite{mybook}.~$\square$

\

The following theorem is our main result regarding the spectral function $e(\lambda,x,x)$.

\begin{theorem}
\label{theorem spectral function unmollified two term}
If the point $x\in M$ is nonfocal then
the spectral function $e(\lambda,x,x)$ admits
the two-term asymptotic expansion
(\ref{two-term asymptotic formula for spectral function})
as $\lambda\to+\infty$.
\end{theorem}

\emph{Proof\ }
The result in question is an immediate consequence of
formulae
(\ref{theorem spectral function mollified equation 7}),
Theorem~\ref{theorem spectral function mollified}
and
Lemma~\ref{Singularity of the wave group at time nonzero lemma pointwise}
from the current paper
and Theorem~B.5.1 from \cite{mybook}.~$\square$

\

We now deal with the counting function $N(\lambda)$.

Suppose that we have a Hamiltonian trajectory
$(x^{(j)}(t;y,\eta),\xi^{(j)}(t;y,\eta))$
and a real number $T>0$ such that
$(x^{(j)}(T;y,\eta),\xi^{(j)}(T;y,\eta))=(y,\eta)$.
We will say in this case
that we have a $T$-periodic trajectory originating
from the point $(y,\eta)\in T'M$.

Denote by $(S^*M)^{(j)}$ the unit cosphere bundle,
i.e.~the $(2n-1)$-dimensional surface in the cotangent
bundle defined by the equation $h^{(j)}(y,\eta)=1$.
The unit cosphere bundle is equipped with a natural Lebesgue measure:
the $(2n-1)$-dimensional surface area element on  $(S^*M)^{(j)}$ is
$dy\,d(S^*_yM)^{(j)}$ where $d(S^*_yM)^{(j)}$ is the
$(n-1)$-dimensional surface area
element on the unit cosphere $(S^*_yM)^{(j)}$, see proof of
Lemma~\ref{Singularity of the wave group at time zero lemma}.

Denote by $\Pi^{(j)}$ the set of points in $(S^*M)^{(j)}$
which serve as starting points for periodic trajectories generated by the Hamiltonian
$h^{(j)}$.
It is known, see Lemma 1.3.4 in \cite{mybook}, that the
set $\Pi^{(j)}$ is measurable.

\begin{definition}
\label{definition of nonperiodicity condition 1}
We say that the nonperiodicity condition is fulfilled
if for each
$j=1,\ldots,m^+$ the set $\Pi^{(j)}$ has measure zero.
\end{definition}

We call a $T$-periodic trajectory \emph{absolutely periodic} if
the function
\[
|x^{(j)}(T;y,\eta)-y|^2+|\xi^{(j)}(T;y,\eta)-\eta|^2
\]
has an infinite order zero in the variables $(y,\eta)$, and we denote
by $(\Pi^a)^{(j)}$ the set of points in $(S^*M)^{(j)}$
which serve as starting points for absolutely periodic trajectories
generated by the Hamiltonian $h^{(j)}$.
It is known, see Corollary 1.3.6 in \cite{mybook}, that the
set $(\Pi^a)^{(j)}$ is measurable and,
moreover, the set $\Pi^{(j)}\setminus(\Pi^a)^{(j)}$ has measure zero.
This allows us to reformulate Definition~\ref{definition of nonperiodicity condition 1}
as follows.

\begin{definition}
\label{definition of nonperiodicity condition 2}
We say that the nonperiodicity condition is fulfilled
if for each
$j=1,\ldots,m^+$ the set $(\Pi^a)^{(j)}$ has measure zero.
\end{definition}

In practical applications it is easier to work with
Definition~\ref{definition of nonperiodicity condition 2}
because the set $(\Pi^a)^{(j)}$ is usually much
thinner than the set $\Pi^{(j)}$.

In order to derive a two-term asymptotic formula for the
counting function $N(\lambda)$ we need the following
lemma.

\begin{lemma}
\label{Singularity of the wave group at time nonzero lemma integrated}
Suppose that the nonperiodicity condition is fulfilled.
Then for any complex-valued function $\hat\gamma\in C_0^\infty(\mathbb{R})$
with $\operatorname{supp}\hat\gamma\subset(0,+\infty)$ we have
\begin{equation}
\label{Singularity of the wave group at time nonzero lemma integrated formula}
\int_M
\mathcal{F}^{-1}_{t\to\lambda}[\hat\gamma(t)\operatorname{tr}u(t,y,y)]\,dy=
o(\lambda^{n-1})
\end{equation}
as $\lambda\to+\infty$.
\end{lemma}

\emph{Proof\ }
The result in question is a special case of
Theorem~4.4.1 from \cite{mybook}.~$\square$

\

The following theorem is our main result regarding the counting function $N(\lambda)$.

\begin{theorem}
\label{theorem counting function unmollified two term}
If the nonperiodicity condition is fulfilled then
the counting function $N(\lambda)$ admits
the two-term asymptotic expansion
(\ref{two-term asymptotic formula for counting function})
as $\lambda\to+\infty$.
\end{theorem}

\emph{Proof\ }
The result in question is an immediate consequence of
formulae
(\ref{definition of counting function}),
(\ref{theorem spectral function mollified equation 7}),
Theorem~\ref{theorem spectral function mollified}
and
Lemma~\ref{Singularity of the wave group at time nonzero lemma integrated}
from the current paper
and Theorem~B.5.1 from \cite{mybook}.~$\square$

\section{$\mathrm{U}(m)$ invariance of the second asymptotic coefficient}
\label{U(m) invariance}

We prove in this section that the RHS of formula
(\ref{formula for b(x)})
is invariant under unitary transformations
(\ref{unitary transformation of operator A}),
(\ref{matrix appearing in unitary transformation of operator})
of our operator $A$.
The arguments presented in this section bear some
similarity to those from Section~\ref{U(1) connection},
the main difference being that the unitary matrix-function in question
is now a function on the base manifold~$M$ rather than on $T'M$.

Fix a point $x\in M$ and an index $j$ (index enumerating the eigenvalues
and eigenvectors of the principal symbol) and consider the expression
\begin{multline}
\label{formula for bj(x)}
\int\limits_{h^{(j)}(x,\xi)<1}
\biggl(
[v^{(j)}]^*A_\mathrm{sub}v^{(j)}
\\
-\frac i2
\bigl\{
[v^{(j)}]^*,A_1-h^{(j)},v^{(j)}
\bigr\}
+\frac i{n-1}h^{(j)}\bigl\{[v^{(j)}]^*,v^{(j)}\bigr\}
\biggr)(x,\xi)\,
d\xi\,,
\end{multline}
compare with (\ref{formula for b(x)}).
We will show that this expression
is invariant under the transformation
(\ref{unitary transformation of operator A}),
(\ref{matrix appearing in unitary transformation of operator}).

The transformation
(\ref{unitary transformation of operator A}),
(\ref{matrix appearing in unitary transformation of operator})
induces the following transformation of the principal
and subprincipal symbols of the operator $A$:
\begin{equation}
\label{transformation of the principal symbol}
A_1\mapsto RA_1R^*,
\end{equation}
\begin{equation}
\label{transformation of the subprincipal symbol}
A_\mathrm{sub}\mapsto
RA_\mathrm{sub}R^*
+\frac i2
\left(
R_{x^\alpha}(A_1)_{\xi_\alpha}R^*
-
R(A_1)_{\xi_\alpha}R^*_{x^\alpha}
\right).
\end{equation}
The eigenvalues of the principal symbol remain unchanged,
whereas the eigen\-vectors transform as
\begin{equation}
\label{transformation of the eigenvectors of the principal symbol}
v^{(j)}\mapsto Rv^{(j)}.
\end{equation}
Substituting formulae
(\ref{transformation of the principal symbol})--(\ref{transformation of the eigenvectors of the principal symbol})
into the RHS of
(\ref{formula for bj(x)})
we conclude that the increment of
the expression (\ref{formula for bj(x)}) is
\begin{multline*}
\int\limits_{h^{(j)}(x,\xi)<1}
\biggl(\,
\frac i2[v^{(j)}]^*
\left(
R^*R_{x^\alpha}(A_1)_{\xi_\alpha}-(A_1)_{\xi_\alpha}R^*_{x^\alpha}R
\right)
v^{(j)}
\\
-
\frac i2
\left(
[v^{(j)}]^*R^*_{x^\alpha}R(A_1-h^{(j)})v^{(j)}_{\xi_\alpha}
-
[v^{(j)}_{\xi_\alpha}]^*(A_1-h^{(j)})R^*R_{x^\alpha}v^{(j)}
\right)
\\
+
\frac i{n-1}h^{(j)}
\left(
[v^{(j)}]^*R^*_{x^\alpha}Rv^{(j)}_{\xi_\alpha}
-
[v^{(j)}_{\xi_\alpha}]^*R^*R_{x^\alpha}v^{(j)}
\right)
\biggr)(x,\xi)\,
d\xi\,,
\end{multline*}
which can be rewritten as
\begin{multline*}
-\frac i2\int\limits_{h^{(j)}(x,\xi)<1}
\biggl(
h^{(j)}_{\xi_\alpha}
\left(
[v^{(j)}]^*R^*_{x^\alpha}Rv^{(j)}
-
[v^{(j)}]^*R^*R_{x^\alpha}v^{(j)}
\right)
\\
-\frac 2{n-1}h^{(j)}
\left(
[v^{(j)}]^*R^*_{x^\alpha}Rv^{(j)}_{\xi_\alpha}
-
[v^{(j)}_{\xi_\alpha}]^*R^*R_{x^\alpha}v^{(j)}
\right)
\biggr)(x,\xi)\,
d\xi\,.
\end{multline*}
In view of the identity $R^*R=I$ the above expression can be further simplified,
so that it reads now
\begin{multline}
\label{increment of bj(x) first iteration}
i\int\limits_{h^{(j)}(x,\xi)<1}
\biggl(
h^{(j)}_{\xi_\alpha}[v^{(j)}]^*R^*R_{x^\alpha}v^{(j)}
\\
-\frac1{n-1}h^{(j)}
\left(
[v^{(j)}]^*R^*R_{x^\alpha}v^{(j)}_{\xi_\alpha}
+
[v^{(j)}_{\xi_\alpha}]^*R^*R_{x^\alpha}v^{(j)}
\right)
\biggr)(x,\xi)\,
d\xi\,.
\end{multline}

Denote
\begin{equation}
\label{defiition of B alpha}
B_\alpha(x):=-iR^*R_{x^\alpha}
\end{equation}
and observe that this set of matrices,
enumerated by the tensor index $\alpha$ running through the values $1,\ldots,n$,
is Hermitian.
Denote also $b_\alpha(x,\xi):=[v^{(j)}]^*B_\alpha v^{(j)}$
and observe that these $b_\alpha$ are positively homogeneous in $\xi$ of degree 0.
Then the expression
(\ref{increment of bj(x) first iteration})
can be rewritten as
\begin{equation}
\label{increment of bj(x) second iteration}
-
\int\limits_{h^{(j)}(x,\xi)<1}
\left(
h^{(j)}_{\xi_\alpha}\,b_\alpha
-\frac 1{n-1}\,h^{(j)}\,\frac{\partial b_\alpha}{\partial\xi_\alpha}
\right)\!(x,\xi)\,
d\xi\,.
\end{equation}
Lemma 4.1.4 and formula (1.1.15) from \cite{mybook} tell us
that the expression (\ref{increment of bj(x) second iteration}) is zero.

\section{Teleparallel connection}
\label{Teleparallel connection}

In this section we work under the additional assumptions
(\ref{Assumption 1}), (\ref{Assumption 2}) and (\ref{Assumption 3}),
i.e.~we study a $2\times2$ matrix differential operator in dimension 3
with trace-free principal symbol. Our aim is to show that in this
case the principal symbol generates additional geometric structures
which allow us to reformulate the results of our spectral analysis
in a much clearer geometric language.

Let us show first that the manifold $M$ in this case is parallelizable.
The principal symbol $A_1(x,\xi)$ is linear in $\xi$ so it can be written as
\begin{equation}
\label{principal symbol via Pauli matrices}
A_1(x,\xi)=\sigma^\alpha(x)\,\xi_\alpha\,,
\end{equation}
where $\sigma^\alpha(x)$, $\alpha=1,2,3$, are some trace-free Hermitian
$2\times2$ matrix-functions.
Let us denote the elements of the matrices $\sigma^\alpha$ as $\sigma^\alpha{}_{\dot ab}$
where the dotted index, running through the values $\dot1,\dot2$, enumerates the rows
and the undotted index, running through the values $1,2$, enumerates the columns;
this notation is taken from \cite{MR2670535}. Put
\begin{equation}
\label{frame via Pauli matrices}
V_1{}^\alpha(x):=\operatorname{Re}\sigma^\alpha{}_{\dot12}(x),
\quad
V_2{}^\alpha(x):=-\operatorname{Im}\sigma^\alpha{}_{\dot12}(x),
\quad
V_3{}^\alpha(x):=\sigma^\alpha{}_{\dot11}(x).
\end{equation}
Formula (\ref{frame via Pauli matrices}) defines a triple of smooth real vector fields
$V_j(x)$, $j=1,2,3$, on the manifold $M$.
These vector fields are linearly independent at every point $x$ of the manifold:
this follows from the fact that $\det A_1(x,\xi)\ne0$, $\forall(x,\xi)\in T'M$
(ellipticity). Thus, the triple of vector fields $V_j$
is a \emph{frame}. The existence of a frame means that the manifold $M$
is parallelizable.

Conversely, given a frame $V_j$
we uniquely recover the elliptic principal symbol
$A_1(x,\xi)$ via formulae
(\ref{principal symbol via Pauli matrices}),
(\ref{Pauli matrices 1})
and
(\ref{Pauli matrices 2}).
Thus, a principal symbol is equivalent to a frame.

It is easy to see that the frame elements $V_j$ are
orthonormal with respect to the metric
(\ref{definition of metric}). Moreover, the
metric can be defined directly from the frame as
\begin{equation}
\label{definition of metric via frame}
g^{\alpha\beta}=
V_j{}^\alpha\,V_j{}^\beta\,,
\end{equation}
where the repeated frame index $j$ indicates summation over $j=1,2,3$.
The two definitions of the metric,
(\ref{definition of metric})
and
(\ref{definition of metric via frame}),
are equivalent.

Parallelizability implies orientability.
Having chosen a particular orientation we define the Hodge star
in the standard way. We will use the Hodge star later on in this
section in order to simplify calculations involving the torsion tensor.

Note that the topological invariant $\mathbf{c}$ introduced
in Section~\ref{Main results} in accordance with formula
(\ref{definition of relative orientation})
can be equivalently (and more naturally)
defined in terms of the frame as
\begin{equation}
\label{definition of relative orientation more natural}
\mathbf{c}:=\operatorname{sgn}\det V_j{}^\alpha.
\end{equation}

The crucial new geometric structure is the teleparallel connection.
We already defined it in Section~\ref{Main results} in accordance with
formula (\ref{definition of parallel transport}), i.e.~via the principal
symbol. This connection can be equivalently defined via the frame as follows.
Suppose we have a vector $v$ based at
the point $y\in M$ and we want to construct a parallel vector
$u$ based at the point $x\in M$. We decompose the vector $v$
with respect to the frame at the point $y$, $v=c^jV_j(y)$,
and reassemble it with the same coefficients $c^j$ at the point $x$,
defining $u:=c^jV_j(x)$.

We now define the covariant derivative corresponding to the teleparallel
connection. Our teleparallel connection is a special case of an affine
connection, so we are looking at a covariant derivative acting on
vectors/covectors in the usual manner
\begin{equation}
\label{definition of affine covariant derivative}
\nabla_\mu v^\alpha=\partial v^\alpha/\partial x^\mu+\Gamma^\alpha{}_{\mu\beta}\,v^\beta\,,
\qquad
\nabla_\mu w_\beta=\partial w_\beta/\partial x^\mu-\Gamma^\alpha{}_{\mu\beta}\,w_\alpha\,.
\end{equation}
Of course, the above $\nabla$ should not be confused with the $\nabla$ from
Section~\ref{U(1) connection}.
The teleparallel connection coefficients are defined from the conditions
\begin{equation}
\label{conditions for teleparallel connection coefficients}
\nabla_\mu V_j{}^\alpha=0\,,
\end{equation}
where the $V_j$ are elements of our frame.
Formula (\ref{conditions for teleparallel connection coefficients}) gives a system
of 27 linear algebraic equations for the determination of 27 unknown
connection coefficients. It is known
(see, for example, formula (A2) in \cite{MR2573111}),
that the unique solution of this system is
\begin{equation}
\label{formula for teleparallel connection coefficients}
\Gamma^\alpha{}_{\mu\beta}=V_k{}^\alpha(\partial\mathbf{V}_{k\beta}/\partial x^\mu)\,,
\end{equation}
where
\begin{equation}
\label{definition of coframe}
\mathbf{V}_{k\beta}:=g_{\beta\gamma}V_k{}^\gamma\,.
\end{equation}
The triple of covector fields $\mathbf{V}_k$, $k=1,2,3$, is called the \emph{coframe}.
The frame and coframe uniquely determine each other via the relation
\begin{equation}
\label{alternative definition of coframe}
V_j{}^\alpha\mathbf{V}_{k\alpha}=\delta_{jk}.
\end{equation}

One can check by performing explicit calculations that the teleparallel connection
has the following two important properties:
\begin{itemize}
\item
$\nabla_\alpha g_{\beta\gamma}=0$, which means that the connection is metric compatible and
\item
$\nabla_\alpha\nabla_\beta-\nabla_\beta\nabla_\alpha=0$,
which means that the Riemann curvature tensor is zero.
\end{itemize}

The tensor characterising the ``strength'' of the teleparallel
connection is not the Riemann curvature tensor but the torsion
tensor (\ref{definition of torsion}).
The teleparallel connection is, in a sense, the opposite of the
more common Levi-Civita connection:
the Levi-Civita connection has zero torsion but nonzero curvature,
whereas the teleparallel connection has nonzero torsion but zero curvature.
In our paper we distinguish these two affine connections by using
different notation for connection coefficients:
we write the teleparallel connection coefficients
as $\Gamma^\alpha{}_{\beta\gamma}$
and the Levi-Civita connection coefficients (Christoffel symbols)
as $\left\{{{\alpha}\atop{\beta\gamma}}\right\}$,
see formula (\ref{definition of Christoffel symbols}).

Substituting
(\ref{formula for teleparallel connection coefficients})
into
(\ref{definition of torsion})
we arrive at the following explicit formula for the torsion tensor
of the teleparallel connection
\begin{equation}
\label{explicit formula for torsion}
T=V_j\otimes d\mathbf{V}_j\,,
\end{equation}
where the $d$ stands for the exterior derivative.
For the sake of clarity we rewrite formula
(\ref{explicit formula for torsion}) in more detailed form, retaining all tensor indices,
\begin{equation}
\label{more explicit formula for torsion}
T^\alpha{}_{\beta\gamma}
=V_j{}^\alpha
(
\partial\mathbf{V}_{j\gamma}/\partial x^\beta
-
\partial\mathbf{V}_{j\beta}/\partial x^\gamma
)\,.
\end{equation}
As always, the repeated index $j$ appearing in formulae
(\ref{explicit formula for torsion})
and
(\ref{more explicit formula for torsion})
indicates summation over $j=1,2,3$.

As pointed out in Section~\ref{Main results},
it is more convenient to work with the rank two tensor $\overset{*}T$
defined by formula (\ref{definition of torsion with a star})
rather than with the rank three tensor $T$.
Substituting
(\ref{explicit formula for torsion})
into
(\ref{definition of torsion with a star})
we get
\begin{equation}
\label{explicit formula for torsion with a star}
\overset{*}T=V_j\otimes\operatorname{curl}\mathbf{V}_j\,,
\end{equation}
where
\begin{equation}
\label{definition of curl}
(\operatorname{curl}\mathbf{V}_j)_\beta:=(*d\mathbf{V}_j)_\beta
=\frac12\,(d\mathbf{V}_j)^{\gamma\delta}\,\varepsilon_{\gamma\delta\beta}
\,\sqrt{\det g_{\mu\nu}}\,.
\end{equation}

The remainder of this section is devoted to the proof of formula
(\ref{curvature via torsion and metric})
expressing the scalar curvature of the $\mathrm{U}(1)$ connection via
the torsion of the teleparallel connection and the metric.

We fix an arbitrary point $Q\in T'M$ and prove formula
(\ref{curvature via torsion and metric}) at this point.
As the LHS and RHS of (\ref{curvature via torsion and metric}) are invariant under
changes of local coordinates~$x$,
it is sufficient to prove formula
(\ref{curvature via torsion and metric}) in
Riemann normal coordinates, i.e.~local coordinates
such that $x=0$ corresponds to the projection of the point $Q$ onto the base manifold,
$g_{\mu\nu}(0)=\delta_{\mu\nu}$ and $\frac{\partial g_{\mu\nu}}{\partial x^\lambda}(0)=0$.
Moreover, as the formula we are proving involves only
first partial derivatives, we may assume, without loss of generality,
that $g_{\mu\nu}(x)=\delta_{\mu\nu}$ for all $x$ in some neighbourhood of the origin.
Thus, it is sufficient to prove formula (\ref{curvature via torsion and metric})
for the case of Euclidean metric.

As both the LHS and RHS of (\ref{curvature via torsion and metric})
have the same degree of homogeneity in~$\xi$, namely, $-1$, it is sufficient
to prove formula (\ref{curvature via torsion and metric}) for $\xi$ of norm 1.
Moreover, by rotating our Cartesian coordinate system we can reduce the
case of general $\xi$ of norm 1 to the case
\begin{equation}
\label{special momentum}
\xi=
\begin{pmatrix}
0&0&1
\end{pmatrix}.
\end{equation}

There is one further simplification that can be made: we claim
that it is sufficient to prove
formula (\ref{curvature via torsion and metric}) for the case when
\begin{equation}
\label{special frame}
V_j{}^\alpha(0)=\mathbf{c}\delta_j{}^\alpha,
\end{equation}
i.e.~for the case when at the point $x=0$
the elements of the frame are aligned with the coordinate axes;
here $\mathbf{c}=\pm1$ is the topological invariant defined
in accordance with formula (\ref{definition of relative orientation})
or, equivalently, in accordance with
formula (\ref{definition of relative orientation more natural}).
This claim follows from the observation that the
LHS of formula (\ref{curvature via torsion and metric})
is invariant under rigid special unitary transformations of
the column-function $v^+(x,\xi)$,
\[
v^+\mapsto Rv^+,
\]
where ``rigid'' refers to the fact that the matrix $R\in\mathrm{SU}(2)$ is constant.
Of course, the column-function $Rv^+$ is no longer an eigenvector of the original
principal symbol, but a new principal symbol obtained from the old one by the rigid
special orthogonal transformation of the frame
(\ref{orthogonal transformation of frame 2})
with the $3\times3$ special orthogonal matrix $O$ expressed in terms of the
$2\times2$ special unitary matrix $R$ in accordance with
(\ref{orthogonal transformation of frame 3}).
One can always choose the special unitary matrix $R$ so that
at the point $x=0$
the elements of the new frame are aligned with the coordinate axes
(in fact, there are two possible choices of $R$ which differ by sign).
It remains only to note the well known fact that the tensor
$\overset{*}T$
appearing in the RHS of formula (\ref{curvature via torsion and metric})
is also invariant under rigid special orthogonal transformations of the frame.

Having made all the simplifying assumptions listed above, we are now in a
position to prove formula (\ref{curvature via torsion and metric}).
We give the proof for the case
\begin{equation}
\label{topological invariant equals one}
\mathbf{c}=+1\,.
\end{equation}
There is no need to give a separate proof for
the case $\mathbf{c}=-1$ as the two cases reduce to one another
by means of the identity
(\ref{sum of curvatures is zero})
and the observation that torsion
(\ref{explicit formula for torsion})
is invariant under inversion of the frame.

Let us calculate the RHS of (\ref{curvature via torsion and metric}) first.
In view of (\ref{special frame})
we have, in the linear approximation in $x$,
\begin{equation}
\label{frame in terms of microrotations}
\begin{pmatrix}
V_1{}^1(x)&V_1{}^2(x)&V_1{}^3(x)\\
V_2{}^1(x)&V_2{}^2(x)&V_2{}^3(x)\\
V_3{}^1(x)&V_3{}^2(x)&V_3{}^3(x)
\end{pmatrix}
=
\begin{pmatrix}
1&w^3(x)&-w^2(x)\\
-w^3(x)&1&w^1(x)\\
w^2(x)&-w^1(x)&1
\end{pmatrix},
\end{equation}
where $w$ is some smooth vector-function which vanishes at $x=0$.
Formula (\ref{frame in terms of microrotations}) is the standard formula
for the linearisation of an orthogonal matrix about the identity;
see also formula (10.1) in \cite{rotational_elasticity}.
Note that in Cosserat elasticity literature the vector-function $w$ is called the
\emph{vector of microrotations}. Substituting
(\ref{frame in terms of microrotations}) into
(\ref{explicit formula for torsion with a star}) we get, at $x=0$,
\begin{equation}
\label{torsion with a star in terms of microrotations}
\overset{*}T_{\alpha\beta}=\partial w_\beta/\partial x^\alpha-\delta_{\alpha\beta}\operatorname{div}w,
\end{equation}
which is formula (10.5) from \cite{rotational_elasticity}.
Here we freely lower and raise tensor indices
using the fact that the metric is Euclidean
(in the Euclidean case it does not matter whether a tensor index
comes as a subscript or a superscript).
Substituting
(\ref{torsion with a star in terms of microrotations})
and
(\ref{special momentum})
into the RHS of (\ref{curvature via torsion and metric}) we get,
at our point $Q\in T'M$,
\begin{equation}
\label{RHS of curvature via torsion and metric}
\frac12\,
\frac
{\overset{*}T{}^{\alpha\beta}\xi_\alpha\xi_\beta}
{(g^{\mu\nu}\xi_\mu\xi_\nu)^{3/2}}\,
=
-\frac12(\partial w^1/\partial x^1+\partial w^2/\partial x^2)\,.
\end{equation}

Let us now calculate the LHS of (\ref{curvature via torsion and metric}).
The equation for
the eigenvector $v^+(x,\xi)$
of the principal symbol is
\begin{equation}
\label{equation for v plus}
\begin{pmatrix}
V_3{}^\alpha\xi_\alpha-\|\xi\|&(V_1-iV_2)^\alpha\xi_\alpha
\\
(V_1+iV_2)^\alpha\xi_\alpha&-V_3{}^\alpha\xi_\alpha-\|\xi\|
\end{pmatrix}
\begin{pmatrix}
v^+_1\\ v^+_2
\end{pmatrix}=0\,.
\end{equation}
In view of
(\ref{special momentum}),
(\ref{special frame})
and
(\ref{topological invariant equals one})
the (normalised) solution of
(\ref{equation for v plus})
at our point $Q\in T'M$ is
\[
v^+=
\begin{pmatrix}
1\\0
\end{pmatrix}.
\]
Of course, our $v^+(x,\xi)$
is defined up to the gauge transformation
(\ref{gauge transformation of the eigenvector}),
(\ref{phase appearing in gauge transformation}),
however the LHS of (\ref{curvature via torsion and metric})
is invariant under this gauge transformation.
We now perturb equation
(\ref{equation for v plus})
about the point $Q\in T'M$,
that is, about
\[
x=0,
\qquad
\xi=
\begin{pmatrix}
0&0&1
\end{pmatrix},
\]
making use of formula
(\ref{frame in terms of microrotations}),
which gives us the following equation for the
increment $\delta v^+$ of
the eigenvector $v^+(x,\xi)$
of the principal symbol:
\begin{multline*}
\begin{pmatrix}
0&0
\\
0&-2
\end{pmatrix}
\begin{pmatrix}
\delta v^+_1
\\
\delta v^+_2
\end{pmatrix}
+
\begin{pmatrix}
0&-w^2(x)-iw^1(x)
\\
-w^2(x)+iw^1(x)&0
\end{pmatrix}
\begin{pmatrix}
1
\\
0
\end{pmatrix}
\\
+
\begin{pmatrix}
0&\delta\xi_1-i\delta\xi_2
\\
\delta\xi_1+i\delta\xi_2&-2\delta\xi_3
\end{pmatrix}
\begin{pmatrix}
1
\\
0
\end{pmatrix}
=0,
\end{multline*}
or, equivalently,
\begin{equation}
\label{equation for increment of v plus component 2}
\delta v^+_2
=
\frac12
(-w^2(x)+iw^1(x)+\delta\xi_1+i\delta\xi_2).
\end{equation}
Formula
(\ref{equation for increment of v plus component 2})
has to be supplemented by the normalisation condition
\linebreak
$\|v^+(x,\xi)\|=1$, which in its linearised form reads
\begin{equation}
\label{equation for increment of v plus component 1}
\operatorname{Re}
\delta v^+_1=0.
\end{equation}
Formulae
(\ref{equation for increment of v plus component 1})
and
(\ref{equation for increment of v plus component 2})
define $\delta v^+$
modulo an arbitrary
$\operatorname{Im}
\delta v^+_1$, with this degree of freedom being associated
with the gauge transformation
(\ref{gauge transformation of the eigenvector}),
(\ref{phase appearing in gauge transformation}).
Without loss of generality we may assume that the gauge is
chosen so that
\begin{equation}
\label{equation for increment of v plus component 1 choice of gauge}
\operatorname{Im}
\delta v^+_1=0.
\end{equation}

Combining formulae
(\ref{equation for increment of v plus component 1}),
(\ref{equation for increment of v plus component 1 choice of gauge}) and
(\ref{equation for increment of v plus component 2})
we get
\begin{equation}
\label{formula for increment of v plus}
\delta v^+=
\frac12
\begin{pmatrix}
0
\\
-w^2(x)+iw^1(x)+\delta\xi_1+i\delta\xi_2
\end{pmatrix}.
\end{equation}
Recall that the $w$ appearing in this formula
is some smooth vector-function which vanishes at $x=0$.

Differentiation of
(\ref{formula for increment of v plus})
gives us
\begin{equation}
\label{derivative of v plus in x}
\frac{\partial v^+}{\partial x^\alpha}=
\frac12
\begin{pmatrix}
0
\\
-\partial w^2/\partial x^\alpha+i\partial w^1/\partial x^\alpha
\end{pmatrix},
\end{equation}
\begin{equation}
\label{derivative of v plus in xi}
\frac{\partial v^+}{\partial\xi_1}=
\frac12
\begin{pmatrix}
0
\\
1
\end{pmatrix},
\qquad
\frac{\partial v^+}{\partial\xi_2}=
\frac12
\begin{pmatrix}
0
\\
i
\end{pmatrix},
\qquad
\frac{\partial v^+}{\partial\xi_3}=
0.
\end{equation}
Formulae
(\ref{derivative of v plus in x})
and
(\ref{derivative of v plus in xi})
imply that at our point $Q\in T'M$
\begin{equation}
\label{LHS of curvature via torsion and metric}
-i\{[v^+]^*,v^+\}
=
-\frac12(\partial w^1/\partial x^1+\partial w^2/\partial x^2).
\end{equation}

Comparing formulae
(\ref{RHS of curvature via torsion and metric})
and
(\ref{LHS of curvature via torsion and metric})
and recalling (\ref{topological invariant equals one}),
we arrive at the required result
(\ref{curvature via torsion and metric}).

\

We end this section by writing down an explicit self-contained formula for
the trace of the tensor $\overset{*}T$. Note that according to formula
(\ref{formula for b(x) with assumption 1, 2 and 3}),
it is only the trace of $\overset{*}T$ that we need for our
spectral asymptotics.
Formulae
(\ref{explicit formula for torsion with a star})
and
(\ref{definition of curl})
imply
\begin{multline}
\label{explicit formula for the trace of torsion with a star}
\operatorname{tr}\overset{*}T
=\sqrt{\det g^{\alpha\beta}}\,
\bigl[
\mathbf{V}_{j1}
\,\partial\mathbf{V}_{j3}/\partial x^2
+
\mathbf{V}_{j2}
\,\partial\mathbf{V}_{j1}/\partial x^3
+
\mathbf{V}_{j3}
\,\partial\mathbf{V}_{j2}/\partial x^1
\\
-
\mathbf{V}_{j1}
\,\partial\mathbf{V}_{j2}/\partial x^3
-
\mathbf{V}_{j2}
\,\partial\mathbf{V}_{j3}/\partial x^1
-
\mathbf{V}_{j3}
\,\partial\mathbf{V}_{j1}/\partial x^2
\bigr].
\end{multline}
Here the coframe $\mathbf{V}_j$
is determined from the principal symbol $A_1(x,\xi)$
in accordance with formulae
(\ref{principal symbol via Pauli matrices}),
(\ref{frame via Pauli matrices})
and
(\ref{definition of coframe})
or
(\ref{alternative definition of coframe}),
whereas the metric $g$
is determined from the principal symbol $A_1(x,\xi)$
in accordance with formula
(\ref{definition of metric})
or
(\ref{definition of metric via frame}).

\section{Proof of Theorem \ref{main theorem}}
\label{Proof of Theorem}

As Theorem \ref{main theorem} is an if and only if theorem,
our proof comes in two parts.

\subsection{Part 1 of the proof of Theorem \ref{main theorem}}
\label{Part 1 of the proof of Theorem}

Let $A$ be a massless Dirac operator on half-densities. We need to
prove that
a)~the subprincipal symbol of this operator,
$A_\mathrm{sub}(x)$, is proportional to the identity matrix and b)
the second asymptotic coefficient of the spectral function, $b(x)$,
is zero.

As we have already established the formula for $b(x)$, see
(\ref{formula for b(x) with assumption 1, 2 and 3}), this part of
the proof of Theorem \ref{main theorem} reduces to proving that
the explicit formula for the subprincipal symbol of the
massless Dirac operator on half-densities is
\begin{equation}
\label{Part 1 of the proof of Theorem eq 1}
A_\mathrm{sub}(x)=\frac{\mathbf{c}}4
\,\bigl(\operatorname{tr}\overset{*}T(x)\bigr)\,I\,,
\end{equation}
where $I$ is the $2\times2$ identity matrix.

We give the proof of (\ref{Part 1 of the proof of Theorem eq 1}) for the case
(\ref{topological invariant equals one}).
There is no need to give a separate proof for
the case $\mathbf{c}=-1$ as the two cases reduce to one another
by inversion of the frame:
the full symbol of the massless Dirac operator on half-densities
changes sign under inversion of the frame
and hence its subprincipal symbol
changes sign under inversion of the frame,
whereas torsion
(\ref{explicit formula for torsion})
is invariant under inversion of the frame.

We fix an arbitrary point $P\in M$ and prove the identity
(\ref{Part 1 of the proof of Theorem eq 1}) at this point.
As the LHS and RHS of (\ref{Part 1 of the proof of Theorem eq 1}) are invariant
under changes of local coordinates~$x$,
it is sufficient to check the identity (\ref{Part 1 of the proof of Theorem eq 1}) in
Riemann normal coordinates, i.e.~local coordinates
such that $x=0$ corresponds to the point $P$,
$g_{\mu\nu}(0)=\delta_{\mu\nu}$ and $\frac{\partial g_{\mu\nu}}{\partial x^\lambda}(0)=0$.
Moreover, as the identity we are proving involves only
first partial derivatives, we may assume, without loss of generality,
that $g_{\mu\nu}(x)=\delta_{\mu\nu}$ for all $x$ in some neighbourhood of the origin.
Furthermore, by rotating our Cartesian coordinate system we can
achieve (\ref{special frame}), which opens the way to the use,
in the linear approximation in $x$,
of formula (\ref{frame in terms of microrotations}).

Substituting (\ref{frame in terms of microrotations}) into
(\ref{Pauli matrices 1}), we get,
in the linear approximation in $x$,
\begin{multline}
\label{Part 1 of the proof of Theorem eq 2}
\sigma^1=
\begin{pmatrix}
w^2&1+iw^3\\
1-iw^3&-w^2
\end{pmatrix}
=\sigma_1\,,
\\
\sigma^2=
\begin{pmatrix}
-w^1&-i+w^3\\
i+w^3&w^1
\end{pmatrix}
=\sigma_2\,,
\\
\sigma^3=
\begin{pmatrix}
1&-iw^1-w^2\\
iw^1-w^2&-1
\end{pmatrix}
=\sigma_3\,.
\end{multline}
Recall that the $w$ appearing in this formula
is some smooth vector-function which vanishes at $x=0$.

Substitution of
(\ref{Part 1 of the proof of Theorem eq 2})
into
(\ref{definition of Weyl operator})
(which coincides with
(\ref{definition of Weyl operator on half-densities})
because we assumed the metric to be Euclidean, $g_{\mu\nu}(x)=\delta_{\mu\nu}$)
allows us to evaluate the full symbol
$A(x,\xi)=A_1(x,\xi)+A_0(x)$
of the massless Dirac operator on half-densities:
\begin{multline}
\label{Part 1 of the proof of Theorem eq 3}
A_1(x,\xi)=
\begin{pmatrix}
\xi_3&\xi_1-i\xi_2\\
\xi_1+i\xi_2&-\xi_3
\end{pmatrix}
\\
+
\begin{pmatrix}
w^2\xi_1-w^1\xi_2&iw^3\xi_1+w^3\xi_2+(-iw^1-w^2)\xi_3\\
-iw^3\xi_1+w^3\xi_2+(iw^1-w^2)\xi_3&-w^2\xi_1+w^1\xi_2
\end{pmatrix},
\end{multline}
\begin{equation}
\label{Part 1 of the proof of Theorem eq 4}
A_0(0)=-\frac i4
\begin{pmatrix}
0&1\\
1&0
\end{pmatrix}
\begin{pmatrix}
0&1\\
1&0
\end{pmatrix}
\begin{pmatrix}
\partial w^2/\partial x^1&i\partial w^3/\partial x^1\\
-i\partial w^3/\partial x^1&-\partial w^2/\partial x^1
\end{pmatrix}+\ldots.
\end{equation}
Here formula
(\ref{Part 1 of the proof of Theorem eq 3})
is written in the linear approximation in $x$,
whereas formula
(\ref{Part 1 of the proof of Theorem eq 4})
displays, for the sake of brevity, only one term out of nine
(the one corresponding to $\alpha=\beta=1$ in
(\ref{definition of Weyl operator})) with the remaining
eight terms concealed within the dots $\ldots$.

Substituting
(\ref{Part 1 of the proof of Theorem eq 4})
and
(\ref{Part 1 of the proof of Theorem eq 3})
into
(\ref{definition of subprincipal symbol}),
we get
\begin{equation}
\label{Part 1 of the proof of Theorem eq 5}
A_\mathrm{sub}(0)=
-\frac12\,
(\operatorname{div}w)\,I.
\end{equation}
But, according to
(\ref{torsion with a star in terms of microrotations}),
\begin{equation}
\label{Part 1 of the proof of Theorem eq 6}
\operatorname{tr}\overset{*}T(0)=-2\operatorname{div}w.
\end{equation}
Formulae
(\ref{Part 1 of the proof of Theorem eq 5}),
(\ref{Part 1 of the proof of Theorem eq 6})
and (\ref{topological invariant equals one})
imply formula
(\ref{Part 1 of the proof of Theorem eq 1})
at $x=0$.

\subsection{Part 2 of the proof of Theorem \ref{main theorem}}
\label{Part 2 of the proof of Theorem}

Let $A$ be an operator satisfying assumptions
(\ref{Assumption 1}), (\ref{Assumption 2}) and (\ref{Assumption 3})
and such that
a) the subprincipal symbol of this operator,
$A_\mathrm{sub}(x)$, is proportional to the identity matrix and b)
the second asymptotic coefficient of the spectral function, $b(x)$,
is zero. We need to prove that $A$ is a massless Dirac operator on half-densities.

As we have already established the formula for $b(x)$, see
(\ref{formula for b(x) with assumption 1, 2 and 3}), we have,
for our operator $A$, the identity
(\ref{Part 1 of the proof of Theorem eq 1}).
Let $V_j$ be the frame corresponding to the principal
symbol of the operator $A$,
see formulae
(\ref{principal symbol via Pauli matrices})
and
(\ref{frame via Pauli matrices}).
Now, let $B$ be the massless Dirac operator on half-densities
corresponding to the same frame. Then the principal symbols of the
operators $A$ and $B$ coincide.
But the subprincipal symbols of the
operators $A$ and $B$ coincide as well, as in both cases
these are determined via the frame according to the same formula
(\ref{Part 1 of the proof of Theorem eq 1}) (for the massless
Dirac operator $B$ this is the result from
subsection~\ref{Part 1 of the proof of Theorem}).
A first order differential operator is determined by its
principal and subprincipal symbols, hence, $A=B$.~$\square$

\section{Spectral asymmetry}
\label{Spectral asymmetry}

In this section we deal with the special case when the operator
$A$ is differential (as opposed to pseudodifferential).
No assumptions are made regarding $n$, $m$ or $\operatorname{tr}A_1$.

Our aim is to examine what happens when we change the sign of the operator.
In other words, we compare the original operator $A$ with the operator
$\tilde A:=-A$. In theoretical physics the transformation
$A\mapsto-A$ would be interpreted as time reversal,
see equation (\ref{dynamic equation most basic}).

It is easy to see that for a differential operator the number $m$
(number of equations in our system) has to be even and that the
principal symbol has to have the same number of positive and negative
eigenvalues.
In the notation of Section~\ref{Main results}
this fact can be expressed as $m=2m^+=2m^-$.

It is also easy to see that the
principal symbols of the two operators, $A$ and $\tilde A$,
and the eigenvalues and eigenvectors of the principal symbols
are related as
\begin{equation}
\label{Spectral asymmetry equation 1}
A_1(x,\xi)=\tilde A_1(x,-\xi),
\end{equation}
\begin{equation}
\label{Spectral asymmetry equation 2}
h^{(j)}(x,\xi)=\tilde h^{(j)}(x,-\xi),
\end{equation}
\begin{equation}
\label{Spectral asymmetry equation 3}
v^{(j)}(x,\xi)=\tilde v^{(j)}(x,-\xi),
\end{equation}
whereas the subprincipal symbols are related as
\begin{equation}
\label{Spectral asymmetry equation 4}
A_\mathrm{sub}(x)=-\tilde A_\mathrm{sub}(x).
\end{equation}

Formulae
(\ref{formula for a(x)}),
(\ref{formula for b(x)}),
(\ref{generalised Poisson bracket on matrix-functions}),
(\ref{Poisson bracket on matrix-functions})
and
(\ref{Spectral asymmetry equation 1})--(\ref{Spectral asymmetry equation 4})
imply
\begin{equation}
\label{Spectral asymmetry equation 5}
a(x)=\tilde a(x),
\qquad
b(x)=-\tilde b(x).
\end{equation}
Substituting (\ref{Spectral asymmetry equation 5}) into
(\ref{a via a(x)}) and (\ref{b via b(x)}) we get
\begin{equation}
\label{Spectral asymmetry equation 6}
a=\tilde a,
\qquad
b=-\tilde b.
\end{equation}

Formulae (\ref{two-term asymptotic formula for counting function})
and (\ref{Spectral asymmetry equation 6}) imply that the spectrum
of a generic first order differential operator is asymmetric about $\lambda=0$.
This phenomenon is known in differential geometry as
\emph{spectral asymmetry}
\cite{atiyah_short_paper,atiyah_part_1,atiyah_part_2,atiyah_part_3}.

If we square our operator $A$ and consider the spectral problem
$A^2v=\lambda^2v$,
then the terms $\pm b\lambda^{n-1}$ cancel
out and the second asymptotic coefficient of the counting function
(as well as the spectral function) of the operator $A^2$ turns to zero.
This is in agreement with the known fact that for an even order semi-bounded
matrix differential operator acting on a manifold without boundary
the second asymptotic coefficient of the counting function is zero, see
Section 6 of \cite{VassilievFuncAn1984} and \cite{SafarovIzv1989}.

The case of the massless Dirac operator is special because,
according to Theorem~\ref{main theorem},
the spectrum (as well as the spectral function) of this operator
is asymptotically symmetric about $\lambda=0$ in the two leading terms.
However, despite this asymptotic symmetry,
we believe that for a generic Riemannian 3-manifold the spectrum
of the massless Dirac operator is asymmetric. In stating this belief
we are in agreement with the discussion presented on page 1298 of
\cite{trautman}; note that in the case of an odd-dimensional manifold
the author of \cite{trautman} refers
to the massless Dirac operator as the \emph{Pauli} operator.
And, of course, our belief that for a generic Riemannian 3-manifold the spectrum
of the massless Dirac operator is asymmetric is closely related to the fact that
in dimension 3 the massless Dirac operator commutes with the operator
of charge conjugation, see
formulae
(\ref{commutes})
and
(\ref{antilinear}).

\section{Bibliographic review}
\label{Bibliographic review}

To our knowledge, the first publication on the subject
of two-term spectral asymptotics for systems
was Ivrii's 1980 paper \cite{IvriiDoklady1980}
in Section 2 of
which the author stated, without proof, a formula for the second
asymptotic coefficient of the counting function.
In a subsequent 1982 paper \cite{IvriiFuncAn1982}
Ivrii acknowledged that the formula from
\cite{IvriiDoklady1980} was incorrect and gave a new formula, labelled (0.6), followed by a ``proof''.
In his 1984 Springer Lecture Notes \cite{ivrii_springer_lecture_notes}
Ivrii acknowledged on page 226 that both his
previous formulae for the second asymptotic coefficient were
incorrect and stated, without proof, yet another formula.

Roughly at the same time Rozenblyum \cite{grisha} also stated
a formula for the second asymptotic coefficient of the counting function
of a first order system.

The formulae from \cite{IvriiDoklady1980}, \cite{IvriiFuncAn1982} and \cite{grisha}
are fundamentally flawed because they are proportional to the subprincipal
symbol. As our formulae
(\ref{b via b(x)}) and (\ref{formula for b(x)})
show, the second
asymptotic coefficient of the counting function
may be nonzero even when the subprincipal symbol is zero.
This illustrates, yet again, the difference between scalar
operators and systems.

The formula on page 226 of \cite{ivrii_springer_lecture_notes}
gives an algorithm for the calculation of the correction term
designed to take account of the effect
described in the previous paragraph. This algorithm
requires the evaluation of a limit of a complicated expression
involving the integral, over the cotangent bundle,
of the trace of the symbol of the resolvent of the operator $A$
constructed by means of pseudodifferential calculus. This algorithm
was revisited in Ivrii's 1998 book, see formulae (4.3.39) and (4.2.25)
in \cite{ivrii_book}.

The next contributor to the subject was Safarov
who, in his 1989 DSc Thesis~\cite{SafarovDSc}, wrote down a formula
for the second asymptotic coefficient of the counting function
which was ``almost'' correct.
This formula appears in \cite{SafarovDSc} as formula (2.4).
As explained in Section~\ref{Main results},
Safarov lost only the curvature terms
$\,-\frac{ni}{n-1}\int h^{(j)}\{[v^{(j)}]^*,v^{(j)}\}$.
Safarov's DSc Thesis \cite{SafarovDSc} provides arguments which are sufficiently
detailed and we were able to identify the precise point
(page 163) at which the mistake occurred.

In 1998 Nicoll rederived \cite{NicollPhD} Safarov's formula
(\ref{formula for principal symbol of oscillatory integral})
for the principal symbols of the propagator, using a method
slightly different from \cite{SafarovDSc}, but stopped short
of calculating the second asymptotic coefficient
of the counting function.

In 2007 Kamotski and Ruzhansky \cite{kamotski}
performed an analysis of the
propagator of a first order elliptic system based on the
approach of Rozenblyum \cite{grisha}, but stopped short
of calculating the second asymptotic coefficient
of the counting function.

One of the authors of this paper, Vassiliev, considered systems in Section 6 of
his 1984 paper \cite{VassilievFuncAn1984}. However, that paper dealt with systems of a very special type:
differential (as opposed to pseudodifferential) and of even (as opposed to odd) order.
In this case the second asymptotic coefficients
of the counting function and the spectral function vanish, provided the
manifold does not have a boundary.


\appendix

%
%
%
%
%

\section{The massless Dirac operator}
\label{The massless Dirac operator}

Let $M$ be a 3-dimensional connected compact oriented manifold equipped with a
Riemannian metric $g_{\alpha\beta}$, $\alpha,\beta=1,2,3$ being the
tensor indices. Note that we are more prescriptive in this appendix
than in the main text of the paper: in the main text orientability
and existence of a metric emerged as consequences of the way we
stated the problem, whereas in this appendix they are \emph{a priori}
assumptions.

We work only in local coordinates with prescribed orientation.

It is known \cite{Stiefel,Kirby} that
a 3-dimensional oriented manifold is parallelizable,
i.e.~there exist smooth real vector fields $V_j$, $j=1,2,3$, that
are linearly independent at every point $x$ of the manifold.
(This fact is often referred to as \emph{Steenrod's theorem}.)
Each
vector $V_j(x)$ has coordinate components $V_j{}^\alpha(x)$,
$\alpha=1,2,3$. Note that we use the Latin letter $j$ for
enumerating the vector fields (this is an \emph{anholonomic} or
\emph{frame} index) and the Greek letter $\alpha$ for enumerating
their components (this is a \emph{holonomic} or \emph{tensor}
index). The triple of linearly independent vector fields $V_j$,
$j=1,2,3$, is called a \emph{frame}. Without loss of generality we
assume further on that the vector fields $V_j$ are orthonormal
with respect to our metric: this can always be
achieved by means of the Gram--Schmidt process.

Define Pauli matrices
\begin{equation}
\label{Pauli matrices 1}
\sigma^\alpha(x):=
s^j\,V_j{}^\alpha(x)\,,
\end{equation}
where
\begin{equation}
\label{Pauli matrices 2}
s^1:=
\begin{pmatrix}
0&1\\
1&0
\end{pmatrix}
=s_1\,,
\qquad
s^2:=
\begin{pmatrix}
0&-i\\
i&0
\end{pmatrix}
=s_2\,,
\qquad
s^3:=
\begin{pmatrix}
1&0\\
0&-1
\end{pmatrix}
=s_3\,.
\end{equation}
In formula (\ref{Pauli matrices 1})
summation is carried out over the repeated frame index $j=1,2,3$,
and $\alpha=1,2,3$ is the free tensor index.

The massless Dirac operator is the matrix operator
\begin{equation}
\label{definition of Weyl operator}
W:=-i\sigma^\alpha
\left(
\frac\partial{\partial x^\alpha}
+\frac14\sigma_\beta
\left(
\frac{\partial\sigma^\beta}{\partial x^\alpha}
+\left\{{{\beta}\atop{\alpha\gamma}}\right\}\sigma^\gamma
\right)
\right),
\end{equation}
where summation is carried out over
$\alpha,\beta,\gamma=1,2,3$, and
\begin{equation}
\label{definition of Christoffel symbols}
\left\{{{\beta}\atop{\alpha\gamma}}\right\}:=
\frac12g^{\beta\delta}
\left(
\frac{\partial g_{\gamma\delta}}{\partial x^\alpha}
+
\frac{\partial g_{\alpha\delta}}{\partial x^\gamma}
-
\frac{\partial g_{\alpha\gamma}}{\partial x^\delta}
\right)
\end{equation}
are the Christoffel symbols.
Here and throughout this appendix we raise and lower
tensor indices using the metric.
Note that we chose the letter ``$W$'' for denoting the massless Dirac operator
because in theoretical physics literature it is often referred to as the \emph{Weyl}
operator.

Formula (\ref{definition of Weyl operator}) is the formula from
\cite{MR2670535}, only written in matrix notation (i.e.~without spinor
indices). Note that in the process of transcribing formulae from
\cite{MR2670535} into matrix notation we used the identity
\begin{equation}
\label{raising spinor indices in Pauli matrices}
\epsilon\sigma^\alpha\epsilon=(\sigma^\alpha)^T,
\end{equation}
$\alpha=1,2,3$, where
\begin{equation}
\label{metric spinor}
\epsilon:=
\begin{pmatrix}
0&-1\\
1&0
\end{pmatrix}
\end{equation}
is the `metric spinor'.
The identity (\ref{raising spinor indices in Pauli matrices})
gives a simple way of raising/lowering spinor indices in Pauli matrices in the non-relativistic
($\alpha\ne0$) setting.

Physically, our massless Dirac operator
(\ref{definition of Weyl operator})
describes a single neutrino living in a 3-dimensional compact universe
$M$. The eigenvalues of the massless Dirac operator are the energy levels.

Observe that the sign of $\det V_j{}^\alpha$ is preserved throughout the connected
oriented manifold $M$. Having $\det V_j{}^\alpha>0$ means that our frame has
positive orientation (relative to the prescribed orientation of local coordinates)
and $\det V_j{}^\alpha<0$ means that our frame has negative orientation.
Accordingly,
we say that our massless Dirac operator (\ref{definition of Weyl operator})
has positive/negative orientation depending on the sign of $\det V_j{}^\alpha$.
Of course, the transformation
$W\mapsto-W$ changes the orientation of the massless Dirac operator.

The massless Dirac operator (\ref{definition of Weyl operator}) acts on columns
$v=\begin{pmatrix}v_1&v_2\end{pmatrix}^T$
of complex-valued scalar functions.
In differential geometry this object is referred to as
a (Weyl) spinor so as to emphasise the fact that $v$ transforms
in a particular way under transformations of the orthonormal frame $V$.
However, as in our exposition the frame $V$ is assumed to be chosen
\emph{a priori}, we can treat the components of the spinor as scalars.
This issue will be revisited below when we state Property 4
of the massless Dirac operator.

We now list the main properties of the massless Dirac operator.

\

\textbf{Property 1.}
The massless Dirac operator is invariant under changes of local
coordinates $x$, i.e.~it maps 2-columns of smooth scalar functions $M\to\mathbb{C}^2$
to 2-columns of smooth scalar functions $M\to\mathbb{C}^2$
regardless of the choice of local coordinates.

In order to establish this property we examine separately the two operators
\begin{equation}
\label{definition of Weyl operator part 1}
\sigma^\alpha
\frac\partial{\partial x^\alpha}
\end{equation}
and
\begin{equation}
\label{definition of Weyl operator part 2}
\sigma^\alpha
\sigma_\beta
\left(
\frac{\partial\sigma^\beta}{\partial x^\alpha}
+\left\{{{\beta}\atop{\alpha\gamma}}\right\}\sigma^\gamma
\right)
\end{equation}
appearing in formula (\ref{definition of Weyl operator}).

Let us act with the differential operator (\ref{definition of Weyl operator part 1}) on a
2-column $u:M\to\mathbb{C}^2$ of smooth scalar functions.
Then
$\frac{\partial u}{\partial x^\alpha}$
is a column-valued covector (i.e.~pair of gradients),
$\sigma^\alpha$ is a matrix-valued vector,
so matrix multiplication combined with contraction in $\alpha$
gives a column-valued scalar. Thus, the operator
(\ref{definition of Weyl operator part 1}) is invariant under changes of local coordinates.

As to the multiplication operator (\ref{definition of Weyl operator part 2}), its invariance
follows from the observation that
$
\left(
\frac{\partial\sigma^\beta}{\partial x^\alpha}
+\left\{{{\beta}\atop{\alpha\gamma}}\right\}\sigma^\gamma
\right)
$
is a matrix-valued tensor.

\

\textbf{Property 2.}
The massless Dirac operator
is formally self-adjoint (symmetric)
with respect to the inner product
\begin{equation}
\label{inner product on colums of scalars}
\int_Mv^*w\,\sqrt{\det g_{\alpha\beta}}\,dx
\end{equation}
on 2-columns of smooth scalar functions $v,w:M\to\mathbb{C}^2$.

Indeed, the adjoint operator is
\begin{equation}
\label{definition of Weyl operator adjoint}
W^*=-i\frac1{\sqrt{\det g_{\kappa\lambda}}}
\frac\partial{\partial x^\alpha}\sqrt{\det g_{\mu\nu}}\,\sigma^\alpha
+\frac i4
\left(
\frac{\partial\sigma^\beta}{\partial x^\alpha}
+\left\{{{\beta}\atop{\alpha\gamma}}\right\}\sigma^\gamma
\right)
\sigma_\beta
\sigma^\alpha.
\end{equation}
Comparing formulae
(\ref{definition of Weyl operator})
and
(\ref{definition of Weyl operator adjoint})
we see that in order to prove formal self-adjointness we need to show that
\begin{multline}
\label{formal self-adjointness 1}
\left(
\frac{\partial\sigma^\beta}{\partial x^\alpha}
+\left\{{{\beta}\atop{\alpha\gamma}}\right\}\sigma^\gamma
\right)
\sigma_\beta
\sigma^\alpha
+
\sigma^\alpha
\sigma_\beta
\left(
\frac{\partial\sigma^\beta}{\partial x^\alpha}
+\left\{{{\beta}\atop{\alpha\gamma}}\right\}\sigma^\gamma
\right)
\\
=
\frac4{\sqrt{\det g_{\kappa\lambda}}}
\left(
\frac\partial{\partial x^\alpha}\sqrt{\det g_{\mu\nu}}\,\sigma^\alpha
\right).
\end{multline}
We fix an arbitrary point $P\in M$ and prove the identity
(\ref{formal self-adjointness 1}) at this point.
In view of Property 1, it is sufficient to check the identity (\ref{formal self-adjointness 1}) in
Riemann normal coordinates, i.e.~local coordinates
such that $x=0$ corresponds to the point $P$,
$g_{\mu\nu}(0)=\delta_{\mu\nu}$ and $\frac{\partial g_{\mu\nu}}{\partial x^\lambda}(0)=0$.
Moreover, as the identity we are proving involves only
first partial derivatives, we may assume, without loss of generality,
that $g_{\mu\nu}(x)=\delta_{\mu\nu}$ for all $x$ in some neighbourhood of the origin.
Thus, the problem has been reduced to proving that
variable (i.e.~dependent on $x$) Pauli matrices
in Euclidean space satisfy the identity
\begin{equation}
\label{formal self-adjointness 2}
\left(
\frac{\partial\sigma^\beta}{\partial x^\alpha}
\right)
\sigma^\beta
\sigma^\alpha
+
\sigma^\alpha
\sigma^\beta
\left(
\frac{\partial\sigma^\beta}{\partial x^\alpha}
\right)
=
4
\left(
\frac{\partial\sigma^\alpha}{\partial x^\alpha}
\right).
\end{equation}
Note that in (\ref{formal self-adjointness 2}) we made all the
tensor indices upper, using the fact that the metric is Euclidean
(in the Euclidean case it does not matter whether a tensor index
comes as a subscript or a superscript). Of course, we still retain
the convention of summation over repeated indices.

In order to prove (\ref{formal self-adjointness 2}) we recall
the basic identity for Pauli matrices which in the Euclidean case reads
\begin{equation}
\label{formal self-adjointness 3}
\sigma^\mu\sigma^\nu+\sigma^\nu\sigma^\mu=2I\delta^{\mu\nu},
\end{equation}
where $I$ is the $2\times2$ identity matrix.
(For a general metric one would have written the above formula
with $g^{\mu\nu}$ instead of $\delta^{\mu\nu}$.)
Formula (\ref{formal self-adjointness 3}) implies
\begin{equation}
\label{formal self-adjointness 4}
\sigma^\mu\sigma^\mu=3I,
\end{equation}
\begin{equation}
\label{formal self-adjointness 5}
\sigma^\mu\sigma^\kappa\sigma^\mu=-\sigma^\kappa,
\end{equation}
\begin{equation}
\label{formal self-adjointness 6}
\partial(\sigma^\mu\sigma^\nu+\sigma^\nu\sigma^\mu)/\partial x^\lambda=0.
\end{equation}
Using formulae
(\ref{formal self-adjointness 3})--(\ref{formal self-adjointness 6})
we get
\begin{multline}
\label{formal self-adjointness 7}
\left(
\frac{\partial\sigma^\beta}{\partial x^\alpha}
\right)
\sigma^\beta
\sigma^\alpha
+
\sigma^\alpha
\sigma^\beta
\left(
\frac{\partial\sigma^\beta}{\partial x^\alpha}
\right)
=
-
\sigma^\beta
\left(
\frac{\partial\sigma^\beta}{\partial x^\alpha}
\right)
\sigma^\alpha
-
\sigma^\alpha
\left(
\frac{\partial\sigma^\beta}{\partial x^\alpha}
\right)
\sigma^\beta
\\
=
\sigma^\beta
\sigma^\beta
\left(
\frac{\partial\sigma^\alpha}{\partial x^\alpha}
\right)
+
\left(
\frac{\partial\sigma^\alpha}{\partial x^\alpha}
\right)
\sigma^\beta
\sigma^\beta
+
\sigma^\beta
\left(
\frac{\partial\sigma^\alpha}{\partial x^\alpha}
\right)
\sigma^\beta
+
\sigma^\beta
\left(
\frac{\partial\sigma^\alpha}{\partial x^\alpha}
\right)
\sigma^\beta
\\
+
\sigma^\beta
\sigma^\alpha
\left(
\frac{\partial\sigma^\beta}{\partial x^\alpha}
\right)
+
\left(
\frac{\partial\sigma^\beta}{\partial x^\alpha}
\right)
\sigma^\alpha
\sigma^\beta
\\
=
3
\left(
\frac{\partial\sigma^\alpha}{\partial x^\alpha}
\right)
+
3
\left(
\frac{\partial\sigma^\alpha}{\partial x^\alpha}
\right)
-
\left(
\frac{\partial\sigma^\alpha}{\partial x^\alpha}
\right)
-
\left(
\frac{\partial\sigma^\alpha}{\partial x^\alpha}
\right)
\\
-
\sigma^\alpha
\sigma^\beta
\left(
\frac{\partial\sigma^\beta}{\partial x^\alpha}
\right)
-
\left(
\frac{\partial\sigma^\beta}{\partial x^\alpha}
\right)
\sigma^\beta
\sigma^\alpha
+
2
\delta^{\alpha\beta}
\left(
\frac{\partial\sigma^\beta}{\partial x^\alpha}
\right)
+
2
\delta^{\alpha\beta}
\left(
\frac{\partial\sigma^\beta}{\partial x^\alpha}
\right)
\\
=
-
\left(
\frac{\partial\sigma^\beta}{\partial x^\alpha}
\right)
\sigma^\beta
\sigma^\alpha
-
\sigma^\alpha
\sigma^\beta
\left(
\frac{\partial\sigma^\beta}{\partial x^\alpha}
\right)
+
8
\left(
\frac{\partial\sigma^\alpha}{\partial x^\alpha}
\right).
\end{multline}
Comparing the left- and right-hand sides of
(\ref{formal self-adjointness 7})
we arrive at
(\ref{formal self-adjointness 2}).

\

\textbf{Property 3.}
The massless Dirac operator $W$ commutes
\begin{equation}
\label{commutes}
\mathrm{C}(Wv)=W\mathrm{C}(v)
\end{equation}
with the antilinear map
\begin{equation}
\label{antilinear}
v\mapsto\mathrm{C}(v):=\epsilon\overline{v}.
\end{equation}
Here the map (\ref{antilinear})
acts on columns
$v=\begin{pmatrix}v_1&v_2\end{pmatrix}^T$
of complex-valued scalar functions,
with $\epsilon$ being the `metric spinor'
defined in accordance with (\ref{metric spinor}).
The commutativity property (\ref{commutes}) follows from the
explicit formula for the massless Dirac operator
(\ref{definition of Weyl operator}) and the identity
$\epsilon\sigma^\alpha=-\overline{\sigma^\alpha}\epsilon$,
$\alpha=1,2,3$,
the latter being a consequence of formula
(\ref{raising spinor indices in Pauli matrices}).

Formula (\ref{commutes}) implies that $v$ is an eigenfunction
of the massless Dirac operator corresponding to an eigenvalue $\lambda$
if and only if
$\mathrm{C}(v)$ is an eigenfunction
of the massless Dirac operator corresponding to the same eigenvalue
$\lambda$.
Hence, all eigenvalues of the massless Dirac operator have even
multiplicity. Moreover, any eigenfunction $v$ and its `partner' $\mathrm{C}(v)$
make the same contribution to the spectral function
(\ref{definition of spectral function})
at every point $x$ of the manifold $M$.


We do not use the commutativity property (\ref{commutes}) of the
massless Dirac operator in the current paper.

The antilinear operator (\ref{antilinear}) is, of course, the charge conjugation
operator which we already encountered in Section \ref{Main results},
see formula (\ref{charge conjugation at the level of principal symbol}).
The difference between the arguments presented in this appendix and those in
Section \ref{Main results} is that in this appendix we deal with the differential
operator, whereas in Section \ref{Main results} we dealt with the principal
symbol. This leads to opposite commutation properties: the charge conjugation
operator commutes with the Weyl operator but it anti\-commutes  with its principal
symbol. The source of this difference is the $\,i\,$ appearing in the RHS of
formula (\ref{definition of Weyl operator}).

\

\textbf{Property 4.}
This property has to do with a particular behaviour under $\mathrm{SU}(2)$ transformations.
Let $R:M\to\mathrm{SU}(2)$ be an arbitrary smooth special unitary matrix-function.
Let us introduce new Pauli matrices
\begin{equation}
\label{special unitary transformation of Pauli matrices}
\tilde\sigma^\alpha:=R\sigma^\alpha R^*
\end{equation}
and a new operator $\tilde W$ obtained by replacing the $\sigma$
in (\ref{definition of Weyl operator}) by $\tilde\sigma$.
It turns out (and this is Property 4) that the two operators,
$\tilde W$ and $W$, are related in exactly the same way as the Pauli
matrices, $\tilde\sigma$ and $\sigma$, that is,
\begin{equation}
\label{special unitary transformation of Weyl operator}
\tilde W=RW R^*.
\end{equation}

In order to prove formula
(\ref{special unitary transformation of Weyl operator})
we write down the operator $\tilde W$ explicitly and rearrange
terms:
\begin{multline*}
\tilde W:=-iR\sigma^\alpha R^*
\left(
\frac\partial{\partial x^\alpha}
+\frac14R\sigma_\beta R^*
\left(
\frac{\partial(R\sigma^\beta R^*)}{\partial x^\alpha}
+\left\{{{\beta}\atop{\alpha\gamma}}\right\}R\sigma^\gamma R^*
\right)
\right)
\\
=-iR\sigma^\alpha\frac\partial{\partial x^\alpha}R^*
+iR\sigma^\alpha
\frac{\partial R^*}{\partial x^\alpha}
\\
-\frac i4R\sigma^\alpha\sigma_\beta
\left(
\frac{\partial\sigma^\beta}{\partial x^\alpha}
+\left\{{{\beta}\atop{\alpha\gamma}}\right\}\sigma^\gamma
\right)R^*
-\frac i4R\sigma^\alpha\sigma_\beta R^*
\left(
\frac{\partial R}{\partial x^\alpha}\sigma^\beta R^*
+
R\sigma^\beta\frac{\partial R^*}{\partial x^\alpha}
\right)
\\
=RWR^*
+iR\sigma^\alpha
\frac{\partial R^*}{\partial x^\alpha}
-\frac i4R\sigma^\alpha\sigma_\beta R^*
\left(
\frac{\partial R}{\partial x^\alpha}\sigma^\beta R^*
+
R\sigma^\beta\frac{\partial R^*}{\partial x^\alpha}
\right).
\end{multline*}
Hence, proving (\ref{special unitary transformation of Weyl operator})
reduces to proving that
\begin{equation}
\label{proving property 3 equation 1}
\sigma^\alpha\sigma_\beta R^*
\left(
\frac{\partial R}{\partial x^\alpha}\sigma^\beta R^*
+
R\sigma^\beta\frac{\partial R^*}{\partial x^\alpha}
\right)
=
4\sigma^\alpha
\frac{\partial R^*}{\partial x^\alpha}\,.
\end{equation}

In order to prove formula (\ref{proving property 3 equation 1})
it is sufficient to show that
\[
\sigma_\beta R^*
\frac{\partial R}{\partial x^\alpha}\sigma^\beta R^*
+
\sigma_\beta
\sigma^\beta\frac{\partial R^*}{\partial x^\alpha}
=
4
\frac{\partial R^*}{\partial x^\alpha}
\]
which, in turn, in view of the identity
$\sigma_\beta\sigma^\beta=3I$
(we already used it in the special case of Euclidean metric,
see formula (\ref{formal self-adjointness 4})),
is equivalent to proving that
\begin{equation}
\label{proving property 3 equation 2}
\sigma_\beta R^*\frac{\partial R}{\partial x^\alpha}\sigma^\beta
=
\frac{\partial R^*}{\partial x^\alpha}R\,.
\end{equation}
The fact that the matrix function $R$ is special unitary implies
that at every point $x$ of the manifold $M$ and for every index $\alpha=1,2,3$
the matrix
$R^*\frac{\partial R}{\partial x^\alpha}$
is trace-free anti-Hermitian, which,
in view of the identity
$\sigma_\beta\sigma^\gamma\sigma^\beta=-\sigma^\gamma$
(we already used it in the special case of Euclidean metric,
see formula (\ref{formal self-adjointness 5})), implies
that formula (\ref{proving property 3 equation 2}) can be
equivalently rewritten as
\begin{equation}
\label{proving property 3 equation 3}
-R^*
\frac{\partial R}{\partial x^\alpha}
=
\frac{\partial R^*}{\partial x^\alpha}R\,.
\end{equation}
But formula (\ref{proving property 3 equation 3}) is an immediate
consequence of the identity $R^*R=I$.


\

Having proved Property 4, let us examine the geometric meaning of the
transformation (\ref{special unitary transformation of Pauli matrices}).
Let us expand the new Pauli matrices $\tilde\sigma$ with respect to
the basis (\ref{Pauli matrices 2}):
\begin{equation}
\label{Pauli matrices 1 new}
\tilde\sigma^\alpha(x)=
s^j\,\tilde V_j{}^\alpha(x).
\end{equation}
Formulae
(\ref{Pauli matrices 1}),
(\ref{Pauli matrices 1 new})
and
(\ref{special unitary transformation of Pauli matrices})
give us the following identity relating
the new vector fields $\tilde V^j$
and the old vector fields $V^j$:
\begin{equation}
\label{orthogonal transformation of frame 1}
Rs^kR^*V_k=s^j\,\tilde V_j\,.
\end{equation}
Resolving (\ref{orthogonal transformation of frame 1})
for $\tilde V_j$ we get
\begin{equation}
\label{orthogonal transformation of frame 2}
\tilde V_j=O_j{}^kV_k\,,
\end{equation}
where the real scalars $O_j{}^k$ are given by the formula
\begin{equation}
\label{orthogonal transformation of frame 3}
O_j{}^k=\frac12\operatorname{tr}(s_jRs^kR^*)\,.
\end{equation}
Note that in writing formulae
(\ref{orthogonal transformation of frame 1}) and (\ref{orthogonal transformation of frame 2})
we chose to hide the tensor index, i.e. we chose to hide the coordinate components of our vector fields.
Say, formula (\ref{orthogonal transformation of frame 2}) written
in more detailed form reads $\tilde V_j{}^\alpha=O_j{}^kV_k{}^\alpha$.

The scalars (\ref{orthogonal transformation of frame 2})
can be viewed as elements of a real $3\times3$ matrix-function $O$
with the first index, $j$, enumerating rows and the second, $k$, enumerating columns.
It is easy to check that this matrix-function $O$ is special orthogonal.
Hence, the new vector fields $\tilde V_j$ are orthonormal and have the same orientation
as the old vector fields $V_j$.
We have shown that the
transformation (\ref{special unitary transformation of Pauli matrices})
has the geometric meaning of switching from our original oriented orthonormal frame $V_j$
to a new oriented orthonormal frame~$\tilde V_j$.

Formula (\ref{orthogonal transformation of frame 3}) means that the special
unitary matrix $R$ is, effectively, a square root of the special orthogonal
matrix $O$.
It is easy to see that for a given matrix $O\in\mathrm{SO}(3)$ formula
(\ref{orthogonal transformation of frame 3}) defines the
matrix $R\in\mathrm{SU}(2)$ uniquely up to sign.
This observation allows us to view the issue
of the geometric meaning of the
transformation (\ref{special unitary transformation of Pauli matrices})
the other way round:
given a pair of orthonormal frames, $V_j$ and $\tilde V_j$, with the
same orientation, we can recover the special orthogonal
matrix-function $O(x)$ from formula (\ref{orthogonal transformation of frame 2})
and then attempt finding a smooth special unitary matrix-function $R(x)$
satisfying (\ref{orthogonal transformation of frame 3}).
Unfortunately, this may not always be possible due to topological obstructions.
We can only guarantee the absence of topological obstructions when the two
frames, $V_j$ and $\tilde V_j$, are sufficiently close to each other, which
is equivalent to saying that
we can only guarantee the absence of topological obstructions when
the special orthogonal matrix-function $O(x)$
is sufficiently close to the identity matrix for all $x\in M$.

We illustrate the possibility of a topological obstruction by means of an explicit example.

\begin{example}
\label{example of 3-torus}
Consider the unit torus $\mathbb{T}^3$ parameterized by cyclic coordinates $x^\alpha$,
$\alpha=1,2,3$, of period $2\pi$. The metric is assumed to be Euclidean. Define the orthonormal
frame as
\begin{equation}
\label{frame on torus}
V_1{}^\alpha=
\begin{pmatrix}
\cos k_3x^3\\
\sin k_3x^3\\
0
\end{pmatrix},
\qquad
V_2{}^\alpha=
\begin{pmatrix}
-\sin k_3x^3\\
\cos k_3x^3\\
0
\end{pmatrix},
\qquad
V_3{}^\alpha=
\begin{pmatrix}
0\\
0\\
1
\end{pmatrix},
\end{equation}
where $k_3\in\mathbb{Z}$ is a parameter.
Let $W$ be the massless Dirac operator corresponding to the frame
(\ref{frame on torus}) with some even $k_3$
and let $\tilde W$ be the massless Dirac operator corresponding to the frame
(\ref{frame on torus}) with some odd $k_3$.
We claim that there does not exist
a smooth matrix-function $R:\mathbb{T}^3\to\mathrm{SU}(2)$ which would give
(\ref{orthogonal transformation of frame 3}),
where $O(x)$ is the special orthogonal matrix-function defined by formula
(\ref{orthogonal transformation of frame 2}).
To prove this, it is sufficient to show that the two operators,
$W$ and $\tilde W$, have different spectra.
Straightforward separation of variables shows that any half-even integer
(positive or negative) is an eigenvalue of $\tilde W$ but is not an eigenvalue of $W$.
What happens in this example is that a special unitary matrix-function $R(x)$
satisfying (\ref{orthogonal transformation of frame 3}) can be defined
locally but not globally: if we try to construct $R(x^3)$
moving along the circumference of the torus $x^3\in(-\pi,\pi)$
we end up with a discontinuity, $\lim\limits_{x^3\to-\pi^+}R(x^3)=-\lim\limits_{x^3\to\pi^-}R(x^3)$.
\end{example}

In fact, one can generalise Example \ref{example of 3-torus}
by introducing rotations in three different directions, which leads
to eight genuinely distinct parallelizations. See also \cite{wild} page~524.

Let us emphasise that the topological obstructions we were discussing
have nothing to do with Stiefel--Whitney classes. We are working on
a parallelizable manifold and the Stiefel--Whitney class of such
a manifold is trivial. The topological issue at hand is that our
parallelizable manifold
may be equipped with different spin structures.

We say that two massless Dirac operators, $W$ and $\tilde W$, are
equivalent if there exists a smooth matrix-function
$R:M\to\mathrm{SU}(2)$ such that the corresponding Pauli matrices,
$\sigma^\alpha$ and $\tilde\sigma^\alpha$, are related in accordance
with (\ref{special unitary transformation of Pauli matrices}).
In view of Property 4 (see formula
(\ref{special unitary transformation of Weyl operator}))
all massless Dirac operators from the same
equivalence class generate the same spectral function
(\ref{definition of spectral function})
and the same counting function
(\ref{definition of counting function}),
so for the purposes of our paper viewing such operators as
equivalent is most natural.

As explained above, there may be many distinct equivalence classes
of massless Dirac operators, the difference between
which is topological. Studying the spectral theoretic implications of
these topological differences is beyond the scope of our paper.
The two-term asymptotics
(\ref{two-term asymptotic formula for spectral function})
and
(\ref{two-term asymptotic formula for counting function})
derived in the main text of our paper do not feel this topology.

In theoretical physics the $\mathrm{SU}(2)$ freedom involved in defining the
massless Dirac operator is interpreted as a gauge degree of
freedom. We do not adopt this point of view (at least explicitly)
in order to fit the massless Dirac operator into the
standard spectral theoretic framework.

We defined the massless Dirac operator
(\ref{definition of Weyl operator})
as an operator acting on 2-columns of scalar functions,
i.e.~on 2-columns of quantities which do not change under changes of
local coordinates. This necessitated the introduction of the density
$\sqrt{\det g_{\alpha\beta}}$ in the formula
(\ref{inner product on colums of scalars}) for the inner product.
In spectral theory it is more common to work with half-densities.
Hence, we introduce the operator
\begin{equation}
\label{definition of Weyl operator on half-densities}
W_{1/2}:=
(\det g_{\kappa\lambda})^{1/4}
\,W\,
(\det g_{\mu\nu})^{-1/4}
\end{equation}
which maps half-densities to half-densities.
We call the operator
(\ref{definition of Weyl operator on half-densities})
\emph{the massless Dirac operator on half-densities}.

\end{document}